\documentclass[final,11pt,authoryear]{elsarticle} 
\usepackage[a4paper,left=2cm,right=2cm,top=2cm,bottom=2cm,]{geometry} 
\usepackage{amssymb}
\usepackage{amsthm}

\usepackage{amsmath,amsfonts,amsbsy}
\allowdisplaybreaks 
\usepackage{bm}
\usepackage[usenames,dvipsnames,svgnames,table]{xcolor}
\usepackage{latexsym}
\usepackage{graphicx}
\usepackage{array}
\usepackage{arrayjobx}
\usepackage{tikz}
\usetikzlibrary{arrows.meta,automata,positioning}
\usepackage{subcaption} 

\usetikzlibrary{calc,patterns,positioning,arrows,chains,shapes,matrix,shapes.symbols,shapes.multipart}
\usepackage{pgfplots}
\pgfplotsset{compat=newest}
\usepgfplotslibrary{statistics}
\usepackage{longtable}
\usepackage{tabularx}
\usepackage{booktabs}

\usepackage{hyperref}

\usepackage{cleveref}
\usepackage{algorithm}
\usepackage[noend]{algpseudocode}

\usepackage{adjustbox} 

\usepackage{textcomp}
\usepackage{mleftright} 

\usepackage{cases} 

\usepackage[symbol]{footmisc}

\usepackage{multicol}
\usepackage{multirow}

\usepackage{diagbox}

\usepackage{makecell}

\usepackage{colortbl}
\usepackage{color}
\definecolor{tsp}{gray}{0.9}
\definecolor{sdvrp}{gray}{0.75}
\definecolor{ctp}{gray}{0.6}
\definecolor{mctp}{gray}{0.45}
\definecolor{hist}{gray}{0.6}

\usepackage{threeparttable, caption}

\usepackage{pgfplots, pgfplotstable}
\pgfplotsset{compat=1.16}

\pgfplotstableread[col sep=comma,]{Plots/M1.csv}{\tableA}
\pgfplotstableread[col sep=comma,]{Plots/M2.csv}{\tableB}
\pgfplotstableread[col sep=comma,]{Plots/M3.csv}{\tableAbis}
\pgfplotstableread[col sep=comma,]{Plots/M4.csv}{\tableBbis}

\usepackage{siunitx}

\begin{document}

	
	\newcommand{\ones}{\mathbf 1}
	\newcommand{\reals}{\mathbb{R}} 
	\newcommand{\binaryDomain}{\{0,1\}} 
	\newcommand{\integers}{\mathbb{Z}} 
	\newcommand{\symm}{{\mbox{\bf S}}}  
	
	\newcommand{\nullspace}{{\mathcal N}}
	\newcommand{\range}{{\mathcal R}}
	\newcommand{\Rank}{\mathop{\bf Rank}}
	\newcommand{\Tr}{\mathop{\bf Tr}}
	\newcommand{\diag}{\mathop{\bf diag}}
	\newcommand{\card}{\mathop{\bf card}}
	\newcommand{\rank}{\mathop{\bf rank}}
	\newcommand{\conv}{\mathop{\bf conv}}
	\newcommand{\prox}{\mathbf{prox}}
	
	\newcommand{\Expect}{\mathop{\bf E{}}}
	\newcommand{\Prob}{\mathop{\bf Prob}}
	\newcommand{\Co}{{\mathop {\bf Co}}} 
	\newcommand{\dist}{\mathop{\bf dist{}}}
	\newcommand{\epi}{\mathop{\bf epi}} 
	\newcommand{\Vol}{\mathop{\bf vol}}
	\newcommand{\dom}{\mathop{\bf dom}} 
	\newcommand{\intr}{\mathop{\bf int}}
	\newcommand{\sign}{\mathop{\bf sign}}
	
	\newcommand{\cf}{{cf. }}
	\newcommand{\eg}{{e.g., }}
	\newcommand{\ie}{{i.e., }}
	\newcommand{\etc}{{etc. }}
	
	\newcommand{\stoppingNodes}{V^\textrm{sto}}
	\newcommand{\stoppingNodesSel}{V^\textrm{sel}}
	\newcommand{\stoppingNodesAlt}{V^\textrm{alt}}
	\newcommand{\demandNodes}{V^\textrm{dem}}
	\newcommand{\prio}[1]{V^\textrm{pref}_{#1}}
	\newcommand{\prioList}[2]{\textrm{pref}(#1,#2)}
	\newcommand{\stoppingNodesBefore}{V^\textrm{before}}
	\newcommand{\stoppingNodesAfter}{V^\textrm{after}}
	
	\newcommand{\wasteNodesPlus}{V^{\textrm{was}+}}
	
	\newcommand{\wasteNodesSelPlus}{V^{\textrm{sel}+}}
	
	\newcommand{\depot}{\sigma}
	\newcommand{\wasteDepot}{\tau}
	\newcommand{\stopTime}{t^{\textrm{sto}}}
	\newcommand{\speedColl}{s^{\textrm{col}}}
	\newcommand{\speedDep}{s^{\textrm{dep}}}
	
	\newcommand{\probAbbr}{WCLRP}
	
	\newcommand{\probAbbrVRP}{C$m$-CTP-R}
	\newcommand{\probAbbrFLP}{FLP-PW}
	\newcommand{\probAbbrLRP}{C-CT-LRP}
	\newcommand{\probAbbrLRPE}{C-CT-LRP-2E}
	
	\newcommand{\RNwS}{RN-wS}
	\newcommand{\RNnoS}{RN-noS}
	\newcommand{\CGwS}{CG-wS}
	\newcommand{\CGnoS}{CG-noS}
	
	\newcounter{commentcounter}
	\setcounter{commentcounter}{1}
	\long\def\symbolfootnote[#1]#2{\begingroup \def\thefootnote{\fnsymbol{footnote}}\footnote[#1]{#2}\endgroup}
	\newcommand{\commentV}[1]{{\footnotesize\textbf{\textcolor{green}{(C.\arabic{commentcounter})}}\symbolfootnote[4]{\texttt{\textcolor{green}{(C.\arabic{commentcounter})~#1}}}}\addtocounter{commentcounter}{1}}
	\newcommand{\commentR}[1]{{\footnotesize\textbf{\textcolor{brown}{(C.\arabic{commentcounter})}}\symbolfootnote[4]{\texttt{\textcolor{brown}{(C.\arabic{commentcounter})~#1}}}}\addtocounter{commentcounter}{1}}
	\newcommand{\commentM}[1]{{\footnotesize\textbf{\textcolor{brown}{(C.\arabic{commentcounter})}}\symbolfootnote[4]{\texttt{\textcolor{blue}{(C.\arabic{commentcounter})~#1}}}}\addtocounter{commentcounter}{1}}

\begin{frontmatter}



\title{A capacitated multi-vehicle covering tour problem on a road network and its application to waste collection}

\date{December 21, 2021}


\author[ufr]{Vera Fischer\corref{cor1}}
\ead{vera.fischer@unifr.ch}

\author[ufr]{Meritxell Pacheco Paneque}
\ead{meritxell.pacheco@unifr.ch}

\author[mtl,cirr,ger]{Antoine Legrain}
\ead{antoine.legrain@polymtl.ca}

\author[hslu]{Reinhard~B\"{u}rgy}
\ead{reinhard.buergy@hslu.ch}

\cortext[cor1]{Corresponding author}

\address[ufr]{University of Fribourg, Fribourg, Switzerland}
\address[hslu]{Lucerne University of Applied Sciences and Arts, Lucerne, Switzerland}
\address[mtl]{\'Ecole Polytechnique de Montr\'{e}al, Montr\'{e}al (Qu\'{e}bec), Canada}
\address[cirr]{CIRRELT - Interuniversity Research Centre on Enterprise Networks, Logistics and Transportation, \\ Montr\'{e}al (Qu\'{e}bec), Canada}
\address[ger]{GERAD - Group for Research in Decision Analysis, Montr\'{e}al (Qu\'{e}bec), Canada}

\begin{abstract}

In most Swiss municipalities, a curbside system consisting of heavy trucks stopping at almost each household is used for non-recoverable waste collection. Due to the many stops of the trucks, this strategy causes high fuel consumption, emissions and noise. These effects can be alleviated by reducing the number of stops performed by collection vehicles. One possibility consists of locating collection points throughout the municipality such that residents bring their waste to their most preferred location. The optimization problem consists of selecting a subset of candidate locations to place the points such that each household disposes the waste at the most preferred location. Provided that the underlying road network is available, we refer to this optimization problem as the capacitated multi-vehicle covering tour problem on a road network (\probAbbrVRP{}). We introduce two mixed-integer linear programming (MILP) formulations: a road-network-based formulation that exploits the sparsity of the network and a customer-based formulation typically used in vehicle routing problems (VRP). To solve large instances, we propose a two-phased heuristic approach that addresses the two subproblems the \probAbbrVRP{} is built on: a set covering problem to select the locations and a split-delivery VRP to determine the routes. Computational experiments on both small and real-life instances show that the road-network-based formulation is better suited. Furthermore, the proposed heuristic provides good solutions with optimality gaps below $0.5\%$ and $3.5\%$ for $75\%$ of the small and real-life instances respectively and is able to find better solutions than the exact method for many real-life instances.

\end{abstract}
\begin{keyword}
waste collection \sep multi-vehicle covering tour problem \sep mixed-integer linear programming \sep road network



\end{keyword}

\end{frontmatter}

\section{Introduction}

Waste collection is an important process in waste management. It mainly involves the transportation of waste from collection sites to disposal facilities, and represents one of the primary and most expensive logistical activities performed by any municipality. Indeed, collection costs of municipal solid waste often account for up to 70\% of the total waste management budget (\citealp{tavares2009optimisation}). The design and operation of a waste collection system is a difficult task, since it entails multiple distinguishing features, such as the types of collected waste (\eg dry recyclable, wet food or residual), the collection frequency (\eg weekly, bi-weekly or on-demand), and the pricing (\eg weight- or volume-based). Other key aspects are the containers or bags being used at collection sites and the collection vehicles. Furthermore, residents have high expectations when it comes to waste collection, in the sense that they aim at a frequent collection at a site close to their home but are not willing to spend neither money nor time for it.

From the residents' viewpoint, collection methods are often divided into curbside (pick-up) systems, where the waste is disposed outside their property, and bring (drop-off) systems, where the waste is brought to communal collection sites (\citealp{rodrigues2016waste}). The former method is used in most Swiss municipalities, with heavy trucks stopping at almost each household to collect non-recoverable waste. Curbside systems are the most convenient for residents. However, due to the nature of the trucks and the number of performed stops, this strategy results into large collection times and causes negative effects such as high fuel consumption, emissions and noise. In hopes of designing a more efficient and sustainable residential waste collection system, in this paper we investigate a waste collection concept that consists of the location of collection points such that they remain close to residential buildings and can be easily accessed by foot. In contrast to conventional bring systems, we take into account residents' preferences when placing collection points within a maximum walking distance from households. This concept favors a smoother transition from the existing curbside system and a broader acceptability among residents.

This problem can be formulated as a facility location problem (FLP; \eg \citealp{ghiani2012capacitated, tralhao2010multiobjective}). Nevertheless, since the location of collection points is interlaced with the subsequent collection tours, both decisions should be simultaneously tackled by means of location-routing problems (LRPs; \citealp{prodhon2014survey}). These problems represent an approach to model and solve locational problems while paying special attention to the underlying issues of vehicle routing (\citealp{nagy2007location}). Despite the different nature of the decisions being addressed (location is strategical whereas routing is tactical/operational), the overall system cost may be excessive if they are handled separately. Additionally, the use of LRPs could decrease the total costs over a long planning horizon within which routes might change (\citealp{salhi1999consistency}).

Given the location of a single disposal facility, a set of candidate locations and a homogeneous fleet of capacitated vehicles, we aim at selecting a subset of candidate locations to place the collection points while determining the routes that visit them in order to collect the waste and transport it to the disposal facility. We assume that a given amount of waste is produced at each residential building and that the waste can be split and transported in different vehicles. For each residential building, a given preference list ranks the candidate locations according to some convenience measure (e.g., walking distance, proximity to interesting points). In short, the goal of the problem is to identify a subset of candidate locations that covers all residential buildings with minimum total travel time. We say that a residential building is covered if its waste is collected at a collection point from its preference list. Besides, the waste must be collected at the first collection point in the preference list that is visited by a vehicle.

As defined, this problem is closely related to the multi-vehicle covering tour problem ($m$-CTP;  \citealp{hachicha2000heuristics}). The $m$-CTP can be seen as a LRP that generalizes the vehicle routing problem (VRP). We model our problem as a particular version of the $m$-CTP in which the constraints on the length and number of vertices of each tour are replaced by vehicle capacity constraints. Moreover, each residential building is not only covered by a collection point but also allocated to the most preferred one belonging to a tour according to its preference list. Since we assume that the underlying road network is known, we refer to this problem as the capacitated $m$-CTP on a road network (\probAbbrVRP{}). 

In this paper, we propose a compact mixed-integer linear programming (MILP) formulation for the \probAbbrVRP{}. This formulation relies on a road-network graph, a more detailed approach to represent the road network than the so-called customer-based graph typically used in VRPs. In a customer-based graph, a node is introduced for each customer and each depot, and an arc represents the shortest path between the start node and the end node. Instead, the proposed formulation exploits the sparsity of the road network by introducing decision variables for each road segment. This technique has already been considered for other problems, such as the Steiner traveling salesman problem (STSP) in \cite{letchford2013compact}, and it is extended here for the \probAbbrVRP{}. We also formulate the model that arises from assuming a customer-based graph on the network for a numerical comparison of the two approaches. Furthermore, we evaluate the impact of split collection (\ie the waste at a collection point can be transported in different vehicles) both computationally and from a practical perspective, as this strategy is usually implemented in the context of waste collection.

To handle practically relevant instances of the \probAbbrVRP{} that general-purpose MILP solvers fail to solve, we develop a two-phased heuristic method. It is based upon solution procedures for the two interdependent subproblems the \probAbbrVRP{} is built on: a set covering problem (SCP) to identify the subset of candidate locations and the split-delivery VRP (SDVRP) to generate the tours. In the first phase, we construct set covers and approximate the routing cost associated with each of them by means of a giant tour on the candidate locations of the set cover and their alternatives, which represent additional candidate locations that can cover the same residential buildings. In the second phase, we solve a SDVRP on the candidate locations visited in the giant tour by transforming it into a capacitated VRP (CVRP) with an a priori splitting strategy (\citealp{chen2017novel}) and by solving the resulting problem with the metaheuristic for CVRP proposed by \cite{vidal2020hybrid}. We consider a set of small instances inspired by \cite{letchford2013compact} and a set of real-life instances from various Swiss municipalitie in the numerical experiments performed to evaluate the MILP formulations, validate and assess the performance of the heuristic method and analyze the obtained savings with respect to the state of practice.

The remainder of the paper is organized as follows. \Cref{sec:relatedWork} provides an overview of relevant works in the context of $m$-CTP and road-network representation. \Cref{sec:problem} formally defines the problem and the particularities for its application in waste collection. Sections~\ref{sec:mip} and \ref{sec:heuristic} present the MILP formulations and the two-phased heuristic method, respectively. \Cref{sec:compTests} reports the numerical experiments, and Section~\ref{sec:conclusion} summarizes the main findings and discusses avenues for future research.

\section{Related work} \label{sec:relatedWork}

Despite the increasing attention received by LRPs in the last years (\citealp{schneider2017survey}), relatively few papers on the covering tour problem (CTP) and derivatives thereof have been published. The first appearance of the CTP can be attributed to \cite{current1981multiobjective}. The CTP is formally defined in \cite{gendreau1997covering} as the problem of finding a Hamiltonian tour with minimum length over the subset of nodes to be visited such that each node in the subset of nodes to be covered lies within a prespecified distance from a tour node. A two-index formulation is developed and it is solved exactly with a branch-and-cut algorithm. Similar to \cite{current1989covering}, they further propose a heuristic approach that combines a SCP and a Traveling Salesman Problem (TSP) heuristics. \cite{baldacci2005scatter} formulate the CTP as a two-commodity network flow problem and propose three scatter-search heuristic algorithms. They show that the resulting bound is tighter than the linear programming (LP) relaxation previously used.

The extension of the CTP to multiple vehicles, the $m$-CTP, is defined in \cite{hachicha2000heuristics} as the problem of designing up to $m$ vehicle tours starting and ending at the depot with minimum total length such that the nodes to cover lie within a preset distance of a tour node and both the number of nodes and the length of any tour do not exceed given values. The $m$-CTP reduces to a VRP with unit demands when all nodes must be visited. The authors point out that the $m$-CTP appears to be more difficult than the VRP and therefore heuristics might be the only viable methods to find good solutions for practically relevant instances. They propose three heuristics called modified savings, modified sweep and route-first/cluster-second that are partially based on the corresponding methods for the standard VRP. The performed experiments show that all of them allow to find good solutions for realistic instances within a reasonable computational time, being the modified savings the fastest one and the other two better in terms of solution quality.

\cite{jozefowiez2014branch} introduces a branch-and-price algorithm based on a column generation approach to solve the $m$-CTP exactly. The master problem is a SCP and the subproblem is a variant of the profitable tour problem (PTP; \citealp{dell1995prize}) solved by a branch-and-cut algorithm. This methodology is tested on instances with up to 60 nodes that can be visited and up to 150 nodes to cover. \cite{ha2013exact} propose a two-commodity flow formulation based on the formulation of \citealp{baldacci2005scatter} for the $m$-CTP when the length constraint is relaxed and $m$ is a decision variable. They consider a standard branch-and-cut algorithm to exactly solve the problem. Computational results show that it outperforms the algorithm by \cite{jozefowiez2014branch} in the same context. They also develop a two-phased metaheuristic based on an evolutionary local search (ELS). In the first phase, subsets of nodes that cover all customers are created, and in the second phase, a VRP with unit demands on each subset is solved. The generated solution is within 1.5\% of optimality for the considered test instances with up to 200 nodes. \cite{kammoun2017integration} apply a variable neighborhood search (VNS) heuristic based on the variable neighborhood descent (VND) method for the variant without the length constraint. This method is compared against the one of \cite{ha2013exact} on the same instances and reports better or equal results in a smaller computational time.

The $m$-CTP finds applications in problems that concern the design of bilevel transportation networks. In these problems, only a subset of the nodes is actually visited by vehicles. Some examples include the works of \cite{hachicha2000heuristics} on the location of mobile health-care teams in rural areas; \cite{labbe1986maximizing} on the location of boxes for overnight mail service; and \cite{oliveira2015multi} on the planning of routes for urban patrolling. Additional examples can be found in humanitarian logistics. \cite{naji2012covering} address the location of satellite distribution centers for supplying humanitarian aid throughout a disaster area. Their numerical experiments on randomly-generated data show that only very small instances can be solved efficiently using that mathematical model. The proposed multi-start heuristic produces high-quality solutions for realistic instances in reasonable computational times. More recently, \cite{davoodi2019integrated} deal with a similar problem by developing a hybrid approach that combines an exact solution method (Benders decomposition) and a fast metaheuristic (VNS) to enhance the efficiency of the overall method, as confirmed by the experiments performed on a real-life case study. Another interesting application is school bus routing. In \cite{schittekat2013metaheuristic}, the joint problem of bus stop selection and bus route generation is formulated in its most basic form. The authors characterize a MILP formulation and a parameter-free matheuristic that combines a greedy randomized adaptive search procedure (GRASP) and a VND method. Experiments on randomly-generated instances with up to 80 stops and 800 students show that the matheuristic finds most known optimal solutions much faster than the exact method.

In the context of waste collection, \cite{cubillos2020solution} rely on a bi-objective LRP known as maximal CTP (MCTP) for the location of recycling drop-off stations. As introduced by \cite{current1994median}, in a MCTP a tour must visit $p$ nodes out of $n$ candidate locations with the goals of minimizing the total length of the tour and maximizing the covered demand (i.e., demand located at a certain distance of the visited nodes). The MCTP is a variant of the CTP in which the covering objective is replaced with a constraint that requires complete coverage. In \cite{cubillos2020solution}, the collection costs associated with the routing are approximated with a TSP that is heuristically solved, which means that the capacity on the collection vehicles is disregarded. To handle real-life sized problems, they propose a heuristic method inspired by VNS because its use for location problems has yielded good results in the literature.

Several of the above-mentioned works rely on a customer-based graph to represent the underlying road network. These are complete graphs where a node represents a point of interest (\eg customers, depots) and an arc represents the best path (\eg shortest, fastest) between two nodes (\citealp{huang2017time}). When multiple attributes are defined on road segments (\eg travel cost, distance), this representation can have negative consequences on the solution quality and/or efficiency (\citealp{ben2018vehicle}). To address this issue, a growing number of papers investigate road-network graphs. These graphs mimic the road network by defining the arcs as the road segments and the nodes as the extremities of these segments. \cite{letchford2014pricing} note that it seems both more simple and more natural to work with the original road network. They also show that significant computing time savings can be achieved with respect to customer-based graphs.

There is a collection of papers relying on road-network graphs in the context of STSP. Independently introduced by \cite{orloff1974fundamental}, \cite{fleischmann1985cutting} and \cite{cornuejols1985traveling}, the aim of this problem is to find a minimum-cost cycle that visits a set of required nodes at least once in a road-network graph. \cite{letchford2013compact} propose compact formulations for the STSP with a linear number of variables and constraints where decision variables are introduced for each road segment. The authors show that instances with up to 500 nodes can be solved with the developed exact branch-and-cut algorithms. Note that the STSP has been transformed into the classical TSP in \cite{alvarez2019note}. Clearly, the TSP is a well studied problem and can be solved quickly with the state-of-the-art solver CONCORDE (\citealp{applegate2003implementing}) which enables to solve all instances from the literature to optimality within 20 seconds (most of them within a second). Road-network graphs are as well receiving increasing attention in time-dependent VRP (e.g., \citealp{ben2019branch}).

In conclusion, the review of the literature shows that a variety of real-life problems can be modeled as $m$-CTPs. We show that this is also the case for the \probAbbrVRP{}. To the best of our knowledge, the location of waste collection sites has not yet been modeled within the $m$-CTP framework. Our problem differs from the standard $m$-CTP in that vehicles are capacitated, the nodes that must be covered are allocated to specific visited nodes and the demand associated with a visited node can be split among vehicles. Some of the reviewed applications have already addressed some of these features. In particular, \cite{schittekat2013metaheuristic} and \cite{davoodi2019integrated} include capacitated vehicles and allocation but do not allow for splits. \cite{naji2012covering} consider in addition to multiple commodities and a heterogeneous fleet of vehicles also capacitated vehicles and split delivery, but covered nodes are not allocated to visited nodes. In our case, to model allocation, we take into account the preferences of the residents on where to bring their waste. Indeed, the ordering given by the preference lists of residents is used in the decision-making process of determining the location of collection points. Note that in other applications (\eg \citealp{schittekat2013metaheuristic}), covered nodes are simply assigned to visited nodes following a criterion set by the decision-maker (\eg proximity). Furthermore, we rely on a road-network graph to represent the underlying road network.

\section{Problem definition}
\label{sec:problem}

We provide a formal definition of the \probAbbrVRP{} in Section~\ref{sec:problemStatement} and discuss its specificities in the context of waste collection in Section~\ref{sec:wasteCollApp}.

\subsection{Notation, data, and problem statement}
\label{sec:problemStatement}

Let $G=(V,A)$ be a directed graph with node set $V$ and arc set $A$ representing the (directed) road network and $W$ the set of nodes that have a positive demand that needs to be satisfied, \ie picked up or delivered (no mix). For each demand node $i \in W$, its demand $d_i$ must be satisfied at one and only one node from its preference list $\prio{i} \subseteq V$. Note that node $i \in W$ may or may not be in $\prio{i}$. We assume that $\prio{i}$ is totally ordered. The ordering reflects the preference associated with demand node $i$, which means that $d_i$ must be satisfied at the first node in $\prio{i}$ at which a vehicle stops. We denote by $\prioList{i}{j}$ the index of node $j$ in $\prio{i}$. Then, for two nodes $j, j' \in \prio{i}$, $\prioList{i}{j'} < \prioList{i}{j}$ indicates that node $j'$ is preferred over node $j$ by demand node $i$. We define $\stoppingNodes = \cup_{i \in W} \prio{i}$ as the subset of potential stopping nodes, i.e., the candidate locations that can be visited by the vehicles. The remaining nodes in $V \setminus \stoppingNodes$ might represent, for instance, road intersections.

The arc set $A$ represents directed road segments. Let $c_{hh'}$ be the non-negative length associated with arc $(h,h') \in A$. We assume that these lengths satisfy the triangle inequality. For the customer-based graph representation, we define $\ell_{hh'}$ as the shortest path length from any node $h \in V$ to any node $h' \in V$ in $G$. We assume that $G$ is strongly connected, \ie there exists a path from each node $h \in V$ to each node $h' \in V$, so that $\ell_{hh'}$ is well-defined for each ordered pair of nodes $(h,h') \in V \times V$.

The demand is satisfied by $m$ identical vehicles that are located at a depot $\sigma \in V$. Each vehicle has a capacity of $Q$ and drives exactly one tour starting and ending at the depot. We allow for splits, \ie the total demand that needs to be satisfied at a stopping node may be split up between multiple vehicles. We assume that the demands are arbitrarily divisible.

A solution of the \probAbbrVRP{} is specified by exactly $m$ tours. For each tour, we usually only record the nodes at which the vehicle stops to satisfy the demand by assuming that a vehicle travels on a shortest path from one stop to the next. Hence, tour $k$ with $\textrm{st}(k)$ stops can be specified by the sequence \begin{equation*}
	\left((\sigma{}, 0),(j_1^k,q_1^k)),\dots, (j_{\textrm{st}(k)}^k,q^k_{\textrm{st}(k)}), (\sigma{}, 0)\right),
\end{equation*} where each tuple $(j_n^k,q_n^k), n\in \{1,\dots,\textrm{st}(k)\},$ represents the amount of demand $q^k_n\in\reals_{\geq 0}$ that is satisfied at node $j_n^k\in \stoppingNodes$. Note that given a solution (or at least the nodes at which the $m$ vehicles stop), it is possible to deduce the total demand that needs to be satisfied at each stopping node. Indeed, for each node $i \in W$, we simply need to assign its demand to the first node in $\prio{i}$ at which a vehicle actually stops.

We say that a solution covers $i \in W$ if at least one vehicle stops at some node in $\prio{i}$. A solution is feasible if it covers all demand nodes, the demand of any node $i \in W$ is satisfied at the first node in $\prio{i}$ at which a vehicle stops and the capacity of the vehicles is not exceeded. The objective of the \probAbbrVRP{} is to find a feasible solution with minimum total cost. The total cost is calculated as the sum of the total cost associated with the traveled distances in the tours plus $r$ times the number of stops performed by all vehicles, where $r$ is a given stop penalty value.

\subsection{Waste collection application}
\label{sec:wasteCollApp}

In this section, we characterize the \probAbbrVRP{} for our application in waste collection. The road network of the geographical area under consideration can be extracted from a mapping service (in our case, \citealp{OpenStreetMap}). For every road intersection, we add a node to $V$. The set of potential collection points $\stoppingNodes$ must be specified by the decision-makers (\eg municipality). We place potential collection points on road segments such that the distance between two points is at most 50 m. Then, any road segment longer than 50 m is split into equal-length stretches of less or equal than 50 m. The resulting splitting points are added to $\stoppingNodes$. 

Each residential building is mapped to the node in $V$ whose location is the closest to the building's location. This node becomes a demand node, and therefore belongs to $W$. Notice that multiple buildings might be represented by a single demand node. For the experiments performed in \Cref{sec:compTests}, the demand $d_i$ of node $i \in W$ is obtained by aggregating the average waste production per inhabitant for a given time horizon across the number of inhabitants represented by the node. The preference list $\prio{i}$ of node $i \in W$ contains all potential collection points that are located at a walking distance shorter than $\gamma$. This value represents the maximum walking distance that is assumed for residents to bring their waste to a collection point. We order the points in $\prio{i}$  according to their increasing distance from node $i$. After the preference lists of all demand nodes have been determined, the potential collection points that do not belong to any preference list are deleted from $\stoppingNodes{}$.

We consider a single vehicle with capacity $Q$ to collect the waste. The depot $\sigma{}$ corresponds to the waste disposal facility where vehicles depart and dump the collected waste. Whenever needed (\eg the vehicle capacity has been reached), the vehicle can go to the disposal facility, which is connected to each potential collection point. Note that visiting the disposal facility requires a considerable driving and dumping time. It should therefore only be visited when necessary. In our deterministic setting, the number of visits is simply $m=\lceil d^{\textrm{tot}}/Q\rceil$, where $d^{\textrm{tot}}=\sum_{i\in W} d_i$ is the total amount of waste. Note that a buffer on the vehicle capacity is typically assumed (\eg the vehicle can be filled up to 80\% of its total capacity). From a planning perspective, we say that the vehicle executes $m$ tours starting and ending at the disposal facility.

\begin{figure}[!htbp]
	\centering
	\begin{subfigure}[b]{0.3\textwidth}
		\centering
		\includegraphics[width=1\textwidth]{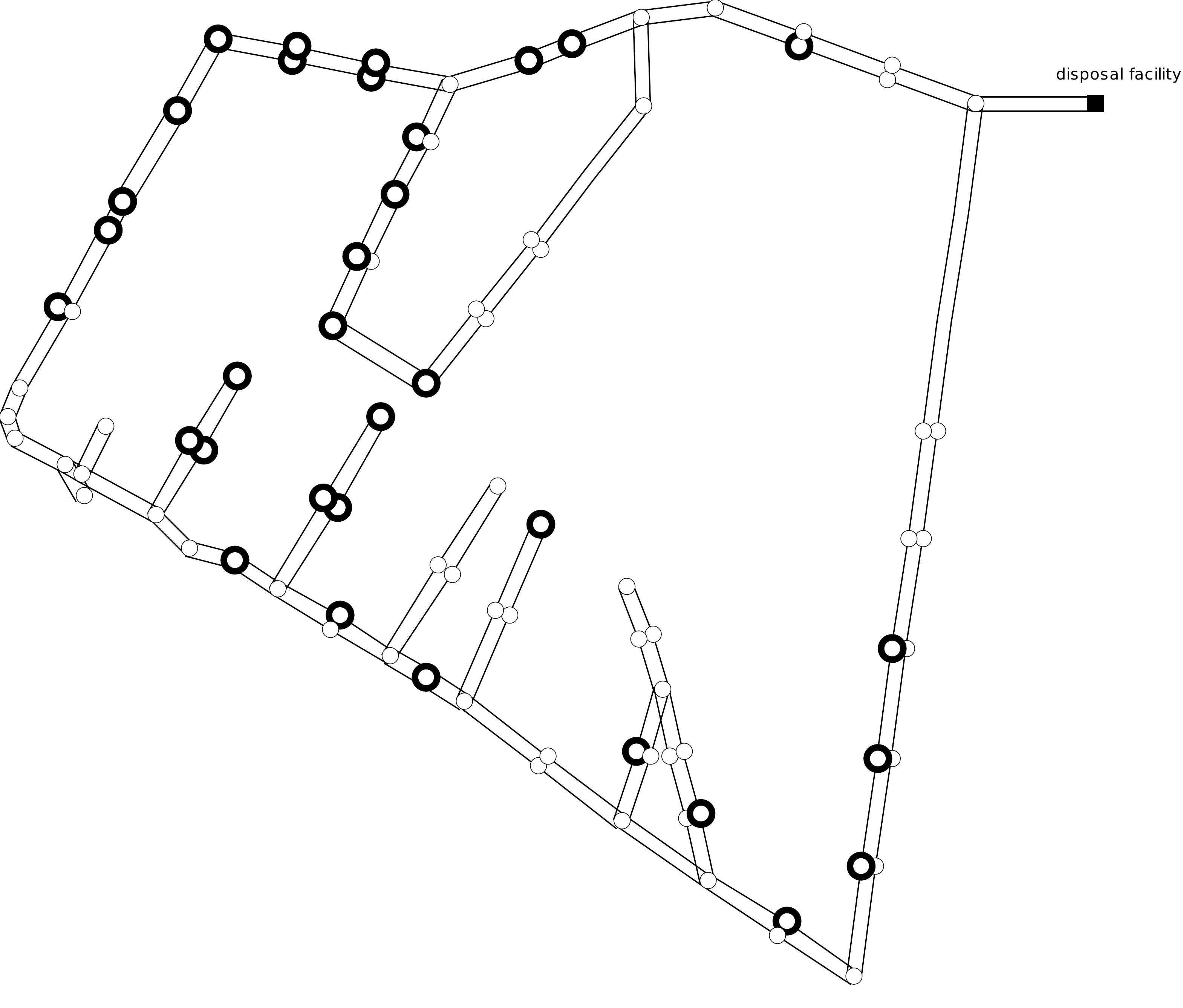}
		\caption{Example}
		\label{fig:graphm1}
	\end{subfigure}
	\hfill
	\begin{subfigure}[b]{0.3\textwidth}
		\centering
		\includegraphics[width=1\textwidth]{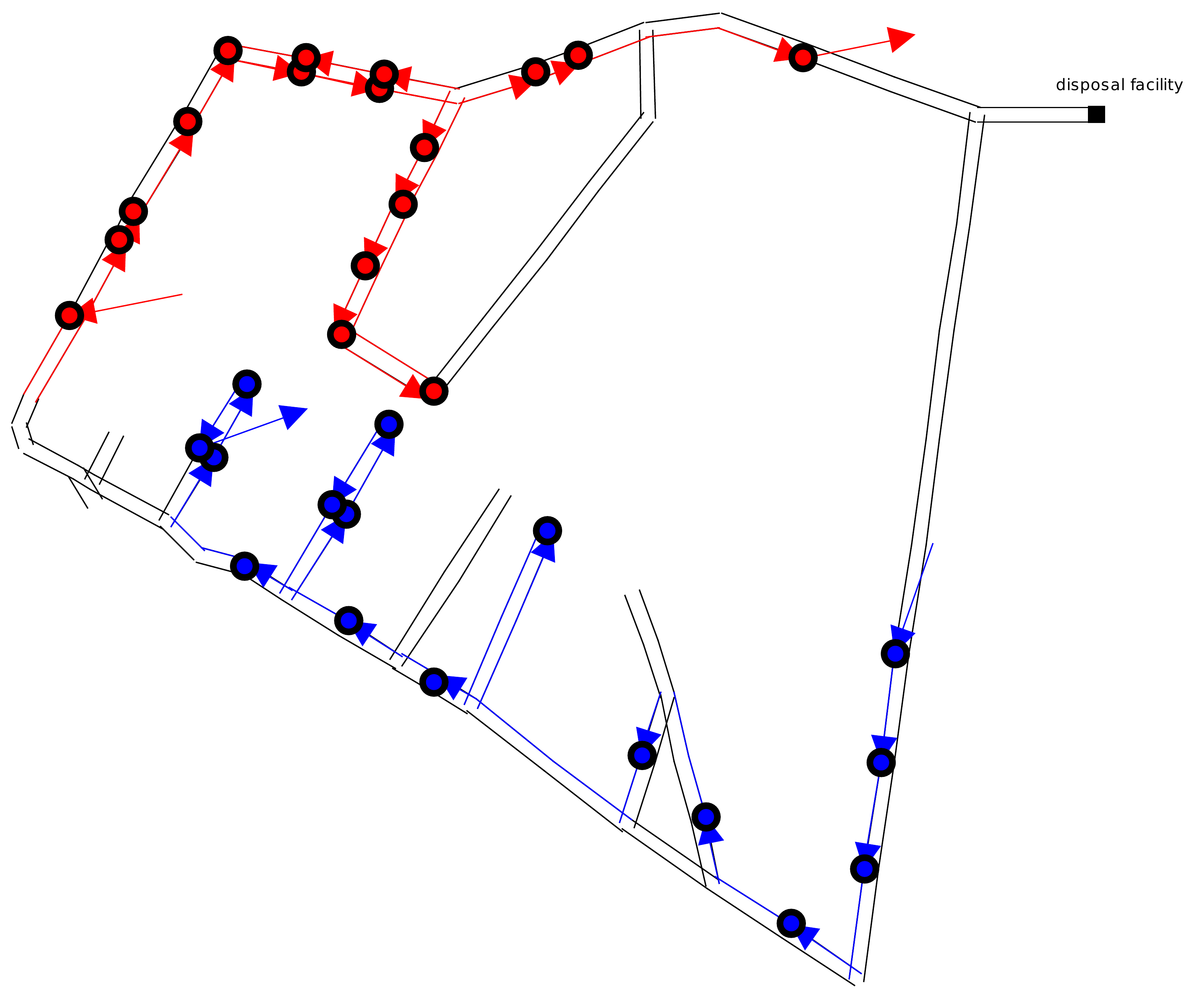}
		\caption{$\gamma = 0$}
		\label{fig:gamma0}
	\end{subfigure}
	\hfill
	\begin{subfigure}[b]{0.3\textwidth}
		\centering
		\includegraphics[width=1\textwidth]{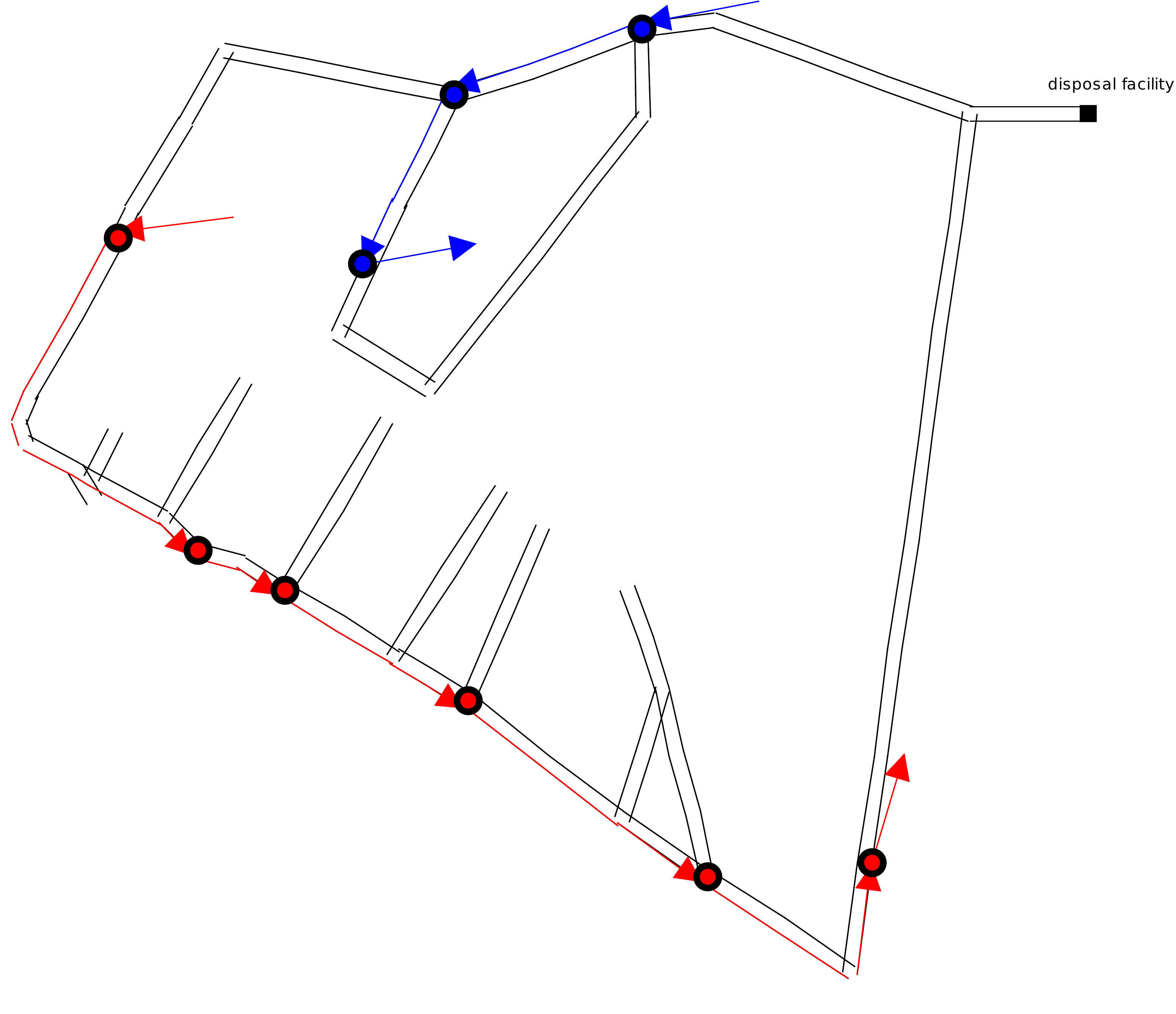}
		\caption{$\gamma = 100$}
		\label{fig:gamma100}
	\end{subfigure}
	\caption{Visualization of our illustrative example and tours 1 (red) and 2 (blue) of an optimal solution with maximum walking distance $\gamma=0$ and $\gamma=100$.}
	\label{fig:graphm1_opt_sol}
\end{figure}

For the sake of illustration, \Cref{fig:graphm1} presents the graph of a small neighborhood of one of the municipalities considered in \Cref{sec:compTests} with 57~residential buildings, 411~inhabitants and an area of 0.13~km$^2$. Graph~$G$ contains in total 97~nodes, out of which 33 are demand nodes, and 307~arcs, out of which 172 are incident to the disposal facility. For readability purposes, we only show the underlying undirected graph of $G$ with left and right road sides explicitly represented. Vehicles are only allowed to turn at the so-called intersections, which are depicted in the graph as nodes in the middle of a road segment (neither on the left nor on the right road side). Furthermore, the waste disposal facility (black square), the demand nodes (bold circles) and the potential collection points (all circles) are depicted. For the disposal facility, only the closest road segments that link it to the municipality road network are drawn.  \Cref{fig:gamma0} and \ref{fig:gamma100} present two solutions for $m=2$ tours to collect the total waste. The solutions have been obtained with the MILP formulation that relies on the road-network graph (see \Cref{section:networkbasedForm}). Two values for the maximum walking distance in meters are considered: $\gamma=0$ (\Cref{fig:gamma0}) and $\gamma=100$ (\Cref{fig:gamma100}). In~\Cref{sec:data_prob_inst} we detail the assumptions on the remaining input data. For $\gamma=0$, tours 1 (red) and 2 (blue) stop at 17 and 16 collection points, respectively. The total travel time is 4411 s and the total stop penalty is 165 s, so the total time is 4576 s. For $\gamma=100$, tours 1 and 2 only stop at 6 and 3 collection points, respectively. The total travel time is 3750 s and the total stop penalty is 45 s, so the total time is 3795 s. This represents almost a 20\% decrease in the total time, which corresponds to a large gain from a practical perspective. 

\section{Compact MILP formulations}
\label{sec:mip}
We propose two compact MILP formulations for the \probAbbrVRP{}. The first formulation (\Cref{section:networkbasedForm}) relies on a road-network graph. For each vehicle, integer decision variables capture the number of times that the vehicle traverses  a road segment. The second one (\Cref{subsection:customerbasedForm}) considers the customer-graph representation of the road network. In this case, binary decision variables for each pair of potential stopping nodes determine whether or not a vehicle successively stops at both nodes. Both formulations allow for splits. This feature may contribute to reduce the overall cost, but increases the search space of the problem and might slow down the solution process. In typical applications of the \probAbbrVRP{}, the demands are small relative to the vehicle capacity. As a consequence, splits might not help much in reducing costs (\citealp{archetti2008, archetti2008split}). We therefore also propose variants to both formulations that forbid splits.

\subsection{Road network-based formulations}
\label{section:networkbasedForm}
Let $\mathcal{M}=\{1,\dots,m\}$ be the set of vehicles enumerated from 1 to $m$. For each arc $(h, h') \in A$ and vehicle $k \in \mathcal{M}$, we introduce an integer variable $x_{hh'k}$ that indicates the number of traversals of vehicle $k$ on arc $(h,h')$. For each node $j \in \stoppingNodes$ and vehicle $k \in \mathcal{M}$, we define a binary variable $y_{jk}$ that takes value 1 if vehicle $k$ stops at node $j$. For each node $i \in W$ and node $j\in\stoppingNodes$, we add a binary variable $z_{ij}$ that takes value 1 if a demand of node $i$ is satisfied at node $j$. For each node $j\in\stoppingNodes$ and vehicle $k \in \mathcal{M}$, we introduce a non-negative continuous variable $q_{jk}$ that indicates the quantity of demand satisfied at node $j$ by vehicle $k$.

We briefly discuss how the tour performed by vehicle $k \in \mathcal{M}$ is constructed with the variables $x$. To this end, we build an auxiliary directed multigraph $G_k^\textrm{T}=(V_k^\textrm{T},A_k^\textrm{T})$, where $V_k^\textrm{T}\subseteq V$. We add a node $h \in V$ to $V_k^\textrm{T}$ if an arc incident to $h$, \ie $(h,h')$ or $(h',h)$, has a positive value $x_{hh'k}$ or $x_{h'hk}$. For any arc $(h,h')\in A$, we add $x_{hh'k}$ copies of $(h,h')$ to $A_k^\textrm{T}$. We then search for an Eulerian cycle in $G_k^\textrm{T}$. Clearly, not all graphs admit an Eulerian cycle. A directed multigraph $D$ admits an Eulerian cycle if and only if $D$ is connected and the in-degree equals the out-degree at each node of $D$ (see, \eg \citealp[Section 1.6]{bang2008digraphs}).

We must therefore ensure these conditions in $G_k^\textrm{T}$. The degree constraints can be directly specified on the variables $x$. For the connectedness condition, we adapt a single-commodity flow-based formulation proposed by \cite{letchford2013compact} for the STSP. For this purpose, we introduce a non-negative continuous variable $f_{hh'k}$ for each vehicle $k \in \mathcal{M}$ to capture the flow passing through $(h, h') \in A$. This flow is determined as follows. For each vehicle $k$, there is a demand equal to $q_{jk}$ associated with node $j\in\stoppingNodes$. This demand must be satisfied by transporting $q_{jk}$ units from node $j$ to the depot by using only the arcs that are traversed by vehicle $k$. Hence, the variable $f_{hh'k}$ can be positive only if $x_{hh'k} \geq 1$. This ensures that all stopping nodes with positive demands visited by vehicle $k$ are connected to the depot in $G_k^\textrm{T}$. Thus, they are all in the same connected component of $G_k^\textrm{T}$. We can efficiently find an Eulerian cycle in this component, for example, with the Hierholzer’s algorithm. Note that in sub-optimal solutions, we may have other connected components in $G_k^\textrm{T}$. These reflect unnecessary traversals and can simply be deleted. Additionally, we may just record the solution in the format given in \Cref{sec:problemStatement} by assuming that a shortest path is followed to travel from one stopping node to the next. These changes can only improve the quality of the solution.

We construct the road-network-based formulation with splits (\RNwS{} for short) as follows:

\footnotesize
\begin{subequations}\label{mip:RNwS}
	\begin{align}
		\text{min } &&& \sum_{k \in \mathcal{M}} \sum_{(h,h') \in A}  c_{hh'} x_{hh'k} +  \sum_{k \in \mathcal{M}} \sum_{j \in \stoppingNodes{}}r y_{jk} \label{eq:OF1}, && \\
		\text{s.t. } &&& \sum_{j \in \prio{i}} z_{ij} =  1 &&\forall i \in W, \label{eq:C1} \\
		&&& \sum_{\substack{j' \in \prio{i}: \\ \prioList{i}{j'} > \prioList{i}{j}}} z_{ij'} \leq 1 - y_{jk} && \forall i \in W, j \in \prio{i}, k \in \mathcal{M}, \label{eq:C2} \\
		&&& \sum_{i \in W: j \in \prio{i}} d_{i} z_{ij} = \sum_{k \in \mathcal{M}} q_{jk} &&\forall j \in \stoppingNodes{}, \label{eq:C3} \\
		&&& \sum_{j \in \stoppingNodes{}} q_{jk} \leq Q &&\forall k \in \mathcal{M}, \label{eq:C3-2} \\
		&&& \sum_{h \in V: (h, j) \in A} x_{hjk} \geq y_{jk} &&\forall j \in \stoppingNodes{}, k \in \mathcal{M}, \label{eq:C4} \\
		&&& \sum_{h' \in V: (h', h) \in A} x_{h'hk} - \sum_{h' \in V: (h, h') \in A} x_{hh'k}= 0 &&\forall h \in V, k \in \mathcal{M}, \label{eq:C5} \\
		&&& \sum_{h' \in V: (h, h') \in A} f_{hh'k} - \sum_{h' \in V: (h', h) \in A} f_{h'hk}  = \begin{cases}
			q_{hk}, & \forall h \in \stoppingNodes{} \label{eq:C6} \\
			0, & \forall h \in V \setminus (\stoppingNodes{} \cup \{\depot{}\})
		\end{cases}  && \forall k \in \mathcal{M}, \\
		&&&\sum_{h \in V: (h,\depot) \in A} f_{h \depot k} = \sum_{j \in \stoppingNodes{}} q_{jk} && \forall  k \in \mathcal{M}, \label{eq:C7} \\
		&&& x_{hh'k} \in \integers_{\geq 0} &&\forall (h, h') \in A, k \in \mathcal{M}, \label{eq:C8} \\ 
		&&& y_{jk} \in \binaryDomain &&\forall  j \in \stoppingNodes{}, k \in \mathcal{M}, \label{eq:C9} \\ 
		&&& z_{ij}  \in \binaryDomain &&\forall  i \in W, j \in \stoppingNodes{}, \label{eq:C10} \\ 
		&&& 0 \leq q_{jk} \leq Q y_{jk}  &&\forall  j \in \stoppingNodes{}, k \in \mathcal{M}, \label{eq:C11} \\ 
		&&& 0 \leq f_{hh'k} \leq Q x_{hh'k} &&\forall (h, h') \in A, k \in \mathcal{M}. \label{eq:C12} 
	\end{align}
\end{subequations}
\normalsize

The objective function~\eqref{eq:OF1} expresses the total cost, which is computed as the sum of the total travel distances plus $r$ times the total number of stops. Constraints~\eqref{eq:C1} ensure that the demand of node $i \in W$ is satisfied at exactly one stopping node from its preference list. Constraints~\eqref{eq:C2} state that this demand is satisfied at the first stopping node in $\prio{i}$ at which some vehicle stops. Constraints~\eqref{eq:C3} guarantee that the total demand that is satisfied at node $j \in \stoppingNodes{}$ is equal to the total quantity that is satisfied by all vehicles stopping at this node. Constraints~\eqref{eq:C3-2} limit the demand satisfied by vehicle $k$ to its capacity $Q$.

Constraints~\eqref{eq:C4} to \eqref{eq:C12} force the variables $x$ to take values that make up valid tours. More precisely, constraints~\eqref{eq:C4} specify that vehicle $k$ can only stop at node $j\in \stoppingNodes{}$ if it traverses an outgoing arc from node $j$ at least once. The degree constraints~\eqref{eq:C5} define that vehicle $k$ enters and leaves any node $h \in V$ the same number of times. Constraints~\eqref{eq:C6} and \eqref{eq:C7} define the flow on the variables $f$ as discussed before. Constraints~\eqref{eq:C6} ensure that the net outflow out of any stopping node $j \in \stoppingNodes{}$ must be $q_{jk}$, which is the demand quantity satisfied at node $j$ by vehicle $k$. For any other node, these constraints impose a 0 net outflow. Constraints~\eqref{eq:C7} enforce that the total quantity that must go into the depot node equals the total quantity satisfied by vehicle $k$.

Finally, constraints~\eqref{eq:C8} to \eqref{eq:C12} define the domain of the decision variables. Additionally, constraints \eqref{eq:C11} link the variables $q$ with the variables $y$ by stating that a positive quantity can be satisfied at node $j \in \stoppingNodes{}$ by vehicle $k$ only if this vehicle stops at this node. Note that $Q$ is a trivial upper bound on the quantity that can be satisfied at any node by one vehicle. Constraints~\eqref{eq:C12} link the variables $f$ with the variables $x$. If $x_{hh'k}=0$, \ie arc $(h,h')$ is not traversed by vehicle $k$, then the flow $f_{hh'k}$ does not pass through it, and therefore $f_{hh'k}=0$. If $x_{hh'k}=1$, the constraint is trivially satisfied as the total flow cannot be larger than $Q$.

As shown in \cite{archetti2006worst}, if the cost matrix satisfies the triangle inequality, then there exists an optimal solution to the SDVRP where the number of splits is less than the number of tours. The number of splits is the sum of the number of splits at each customer, that is defined as the number of tours that visit the customer minus one. We can derive a valid inequality to \RNwS{} from this property with respect to the potential stopping nodes actually visited by the vehicles. This prevents to obtain solutions with many unnecessary splits early in the optimization process. To calculate the number of splits, we introduce the binary variable $s_{j}$ for each node $j \in \stoppingNodes{}$, which takes value 1 if node $j$ is visited at least once. This allows not to take into account the potential stopping nodes that are not visited by any vehicle. These variables are linearly characterized by the following constraints: \begin{subequations}\label{mip:valid_ineq}
	\begin{align}
		s_j &\leq \sum_{k \in \mathcal{M}} y_{jk} && \forall j \in \stoppingNodes{}, \label{eq:valid_ineq_1} \\
		s_j &\geq y_{jk} && \forall j \in \stoppingNodes{}, k \in  \mathcal{M}. \label{eq:valid_ineq_2} 
	\end{align} Constraints~\eqref{eq:valid_ineq_1} force $s_j$ to be equal to 0 if node $j$ is not visited, whereas constraints~\eqref{eq:valid_ineq_2} set $s_j$ to 1 if the node is visited by at least one vehicle. The valid inequality is then written as follows:
	\begin{align}
		\sum_{j \in \stoppingNodes{}} \mleft( \sum_{k \in \mathcal{M}} y_{jk}   - s_j \mright) \leq m-1. && \label{eq:valid_ineq_3}
	\end{align} 
\end{subequations} Note that when node $j \in \stoppingNodes{}$ is not visited by any vehicle, the term in brackets in the left-hand side of constraint~\eqref{eq:valid_ineq_3} is equal to 0. If the node is visited at least once, this term corresponds to the number of splits, \ie number of tours (vehicles) that visit the node minus one.

In order to forbid splits in formulation \eqref{mip:RNwS}, we simply need to add the following constraints: 
\begin{subequations}\label{mip:RNnoS}
	\begin{align}
		\sum_{k \in \mathcal{M}} y_{jk} \leq 1 && \forall j \in \stoppingNodes{}, \label{eq:ctpSplit1}
	\end{align} which impose that at most one vehicle stops at any potential stopping node. These constraints allow to replace constraints~\eqref{eq:C2} with the following stronger variant:
	\begin{align}
		\sum_{\substack{j' \in \prio{i}: \\ \prioList{i}{j'} > \prioList{i}{j}}} z_{ij'} \leq 1 - \sum_{k \in \mathcal{M}} y_{jk} && \forall i \in W, j \in \prio{i}. \label{eq:ctpSplit2}
	\end{align}
\end{subequations} Notice that the valid inequality defined by constraints~\eqref{mip:valid_ineq} is not necessary in this case. We refer to the resulting formulation as \RNnoS{}.

\subsection{Customer‐based graph formulation} \label{subsection:customerbasedForm}
As discussed in \Cref{sec:relatedWork}, VRPs are typically formulated using a so-called customer-based graph. In this section, we show how to construct a customer-based-graph formulation for the \probAbbrVRP{}, called \CGwS{} for short, by adapting \RNwS{}, and then develop the variant that forbids splits, called \CGnoS{}.

Let $G'=(V',A')$ be the complete directed graph made up by the node set $V'=\{\depot{}\} \cup \stoppingNodes{}$ and the arc set $A'$ such that the arc $(j,j')\in A'$ represents a shortest path from $j \in V'$ to $j' \in V'$ of length $\ell_{jj'}$. It is then easy to derive \CGwS{} from \RNwS{}. We only need to replace $V$ with $V'$ and $A$ with $A'$ in formulation~\eqref{mip:RNwS}. We can strengthen constraints~\eqref{eq:C4}. In the customer-based graph, if a vehicle $k$ goes through node $j$, then it also stops there, whereas in the road-network graph it could just pass by without stopping. Thus, we can replace constraints~\eqref{eq:C4} with the following stronger condition:
\begin{equation} \label{mip:strongerC4}
	\sum_{j' \in V'} x_{jj'k} = y_{jk},  \forall j \in \stoppingNodes{}, k \in \mathcal{K}.
\end{equation} In this formulation, the variables $x_{jj'k}$ determine whether vehicle $k$ stops at node $j'$ right after stopping at node $j$. Notice that constraints~\eqref{mip:valid_ineq} are also valid inequalities to \CGwS{}.

When splits are not allowed, \CGnoS{} is a two-index formulation on $G'$. It can be obtained from \CGwS{} by dropping the vehicle index $k$ in the variables and adapting the objective function and constraints accordingly. For each arc $(j,j')\in A'$, we introduce a binary variable $x_{jj'}$ with the same meaning as in \CGwS{}. For each node $j \in \stoppingNodes{}$, we define a binary variable $y_{j}$ that is equal to 1 if node $j$ is visited by a vehicle. For each node $i \in W$ and node $j \in \stoppingNodes{}$, we add a binary variable $z_{ij}$ that is equal to 1 if the demand of node $i$ is satisfied at node $j$. For each node $j \in \stoppingNodes{}$, we introduce a non-negative continuous variable $q_{j}$ that captures the quantity of waste satisfied at node $j$. To prevent subtours, we associate a non-negative continuous variable $f_{jj'}$ with each arc $(j,j') \in A'$ that indicates the number of units of the flow that traverses arc $(j,j')$. \CGnoS{} can be constructed as follows (the components of this formulation have a similar interpretation to their counterparts in \RNnoS{}):

\footnotesize
\begin{subequations}\label{mip:cgnos}
	\begin{align}
		\text{min } &&&  \sum_{j \in \stoppingNodes{}} r y_{j} + \sum_{(j,j')\in A'} \ell_{jj'} x_{jj'} \label{eq:cg0} &&\\
		\text{s.t. } &&& \sum_{j \in \prio{i}} z_{ij} =  1 && \forall i \in W, \label{eq:cg1} \\
		&&& \sum_{\substack{j' \in \prio{i}: \\ \prioList{i}{j'} > \prioList{i}{j}}} z_{ij'} \leq 1 - y_{j} &&\forall i \in W, j \in \prio{i} \label{eq:cg2} \\
		&&& \sum_{i \in W: j \in \prio{i}} d_{i} z_{ij} = q_{j} &&\forall j \in \stoppingNodes{}, \label{eq:cg3} \\
		&&& \sum_{j' \in V'} x_{jj'} = y_{j} && \forall j \in \stoppingNodes{}, \label{eq:cg4} \\
		&&& \sum_{j' \in V'} x_{j'j} - \sum_{j' \in V'} x_{jj'} = 0 && \forall j \in V', \label{eq:cg5} \\
		&&& \sum_{j' \in V' } x_{\depot j'}  = m, \label{eq:cg52} \\
		&&& \sum_{j' \in V'} f_{jj'} - \sum_{j' \in V'} f_{j'j} = q_{j} && \forall j \in \stoppingNodes{},\label{eq:cg6} \\
		&&& \sum_{j \in \stoppingNodes{}}  f_{j\depot} - \sum_{j \in \stoppingNodes{}} f_{\depot j} = \sum_{j \in \stoppingNodes{}} q_{j}, \label{eq:cg7} && \\
		&&& x_{jj'} \in \binaryDomain && \forall (j, j') \in A', \label{eq:cg11} \\
		&&& y_{j} \in \binaryDomain &&\forall j \in \stoppingNodes{}, \label{eq:cg12} \\
		&&& z_{ij}  \in \binaryDomain &&\forall i \in W, j \in \stoppingNodes{}, \label{eq:cg13} \\
		&&& 0 \leq q_{j} \leq Q y_{j}  &&\forall j \in \stoppingNodes{}, \label{eq:cg14}\\
		&&& 0 \leq f_{jj'} \leq Q x_{jj'} &&\forall (j, j') \in A'. \label{eq:cg10}
	\end{align}
\end{subequations} 
\normalsize

\section{A two-phased heuristic method for the \probAbbrVRP{}}
\label{sec:heuristic}

The compact formulations developed in \Cref{sec:mip} become very large for instances of relevant size, which makes general-purpose MILP solvers to fail at finding solutions. In this section, we propose a heuristic method that relies on the two interdependent subproblems the \probAbbrVRP{} is built on. The first subproblem comprises the selection of a subset of potential stopping nodes covering all demand nodes. We call such a subset a set cover. Note that this subproblem is a SCP. Given a set cover, the second subproblem involves the generation of tours visiting its stopping nodes. In principle, this problem is a SDVRP, as the same node might be visited in multiple routes. To increase the flexibility in the routing, we determine a set of alternative nodes, or simply alternatives, for each stopping node in the set cover. We approximate the routing cost associated with a set cover by means of a giant tour that visits either the stopping node or one of its alternatives such that the resulting visited nodes make up a set cover.

The overall method is composed of two phases. In the first phase, we generate a collection of set covers (\Cref{subsec:first}). In the second phase, we build tours on the potential stopping nodes of a given set cover and their alternatives (\Cref{subsec:second}).

\subsection{First phase}
\label{subsec:first}

The aim of this phase is to generate a collection of set covers while taking into account the routing that will be performed in the second phase. To this end, we construct a giant tour using the set cover nodes and their alternatives and redefine the set cover with the nodes "really" visited in such tour. We store the resulting set covers in a list called bestSetCovers according to their giant tour cost. This list has a maximum length (maxLength) and is implemented as a min-max priority queue to efficiently retrieve or delete its minimum or maximum element.

\Cref{algo:first_phase} shows the pseudocode of the first phase. The procedure constructSet constructs a set cover denoted by $\stoppingNodesSel{}$ and calculates the cost of the giant tour associated with it (see \Cref{subsubsec:set_gen}). If this set cover has already been treated in the second phase (see \Cref{subsec:second}), we penalize this fact by multiplying the cost of the associated giant tour (cost($\stoppingNodesSel{}$)) by 1.5. The set cover is added to bestSetCovers if it is different from the ones already included in bestSetCovers, and either (i) cost($\stoppingNodesSel{}$) is lower than the largest giant tour cost (cost($\text{bestSetCovers.max}()$)) or (ii) bestSetCovers has not reached its maximum length (maxLength).

\begin{algorithm}[!ht]
	\caption{First phase of the two-phased heuristic method for the \probAbbrVRP{}}
	\label{algo:first_phase}
	\small
	\begin{algorithmic}[1]
		\State Define $\text{bestSetCovers}$ as an empty min-max priority queue
		\For{$t = 0$ to $\text{maxIterations}$}
		\State Obtain $\stoppingNodesSel{}$ with $\text{constructSet}()$
		\If{$\stoppingNodesSel{} \in$ TreatedSetCovers }
		\State Penalize cost($\stoppingNodesSel{}$)
		\EndIf
		\If{$\stoppingNodesSel{} \notin \text{bestSetCovers}$ and ($\text{cost}(\stoppingNodesSel{}) < \text{cost}(\text{bestSetCovers.max}())$  or $|\text{bestSetCovers}|<\text{maxLength}$)}
		\State $t \leftarrow 0$ 
		\State Add $\stoppingNodesSel{}$ to $\text{bestSetCovers}$
		\EndIf
		
		\EndFor
	\end{algorithmic}
\end{algorithm}

\subsubsection{Set construction} \label{subsubsec:set_gen}

Starting from an empty set, the procedure constructSet generates a set cover $\stoppingNodesSel{}$ by iteratively adding a next node until the resulting set is a set cover. The next node to be added is chosen as follows. Randomly select a demand node $i \in W$ that is not yet covered by $\stoppingNodesSel{}$. To cover node $i$, add one of the stopping nodes in $\prio{i}$ to $\stoppingNodesSel{}$. To intensify the covering of demand nodes, for each node $j \in \prio{i}$, compute the total number of covered demand nodes if node $j$ would be added to $\stoppingNodesSel{}$. Then, select and add to $\stoppingNodesSel{}$ a node maximizing this number. Once $\stoppingNodesSel{}$ is constructed, we remove its redundant nodes. A node $j \in \stoppingNodesSel{}$ is redundant if $\stoppingNodesSel{} \setminus \{j\}$ is a set cover.

Concerning the alternatives, for each node $j \in \stoppingNodesSel{}$, let $W_j \subseteq W$ be the set of nodes whose demand is satisfied at node $j$. It contains all demand nodes $i$ for which node $j$ is the first potential stopping node in $\prio{i}$ among the ones in $\stoppingNodesSel{}$, \ie for any other node $ j' \in \stoppingNodesSel{} \cap \prio{i}$, $\prioList{i}{j} < \prioList{i}{j'}$. Let $\stoppingNodesAlt_{j} \subseteq \stoppingNodes{}$ be the set of alternatives to node $j \in \stoppingNodesSel{}$. It contains the nodes $j' \in \stoppingNodes{} \setminus \{j\}$ that can cover the demand nodes in $W_j$. We denote by $\stoppingNodesAlt = \cup_{j \in \stoppingNodesSel{}} \stoppingNodesAlt_{j}$ the set of alternatives.

We now construct a giant tour using the nodes $\stoppingNodesSel{} \cup \stoppingNodesAlt \subseteq \stoppingNodes{}$. The construction of the giant tour is based on the savings obtained from the gradual insertion of stopping nodes (inspired by \cite{clarke1964scheduling}). To start the giant tour, a node from the given set cover $\stoppingNodesSel{}$ is randomly selected. The subsequent nodes are inserted either before the first node or after the last node in the current giant tour. This is also decided at random. Once the position is determined (either before the first node or after the last one), the nodes are selected according to the largest reported savings among the available groups of potential stopping nodes. A group of potential stopping nodes is defined by a set cover node and its alternatives. As soon as a node from a group is inserted in the giant tour, the group is no longer considered for insertion. We keep on adding nodes to the giant tour until all the demand is satisfied. We finally redefine $\stoppingNodesSel{}$ as the set cover that contains the stopping nodes visited in the giant tour.

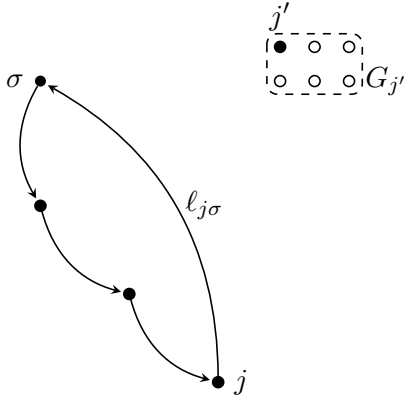
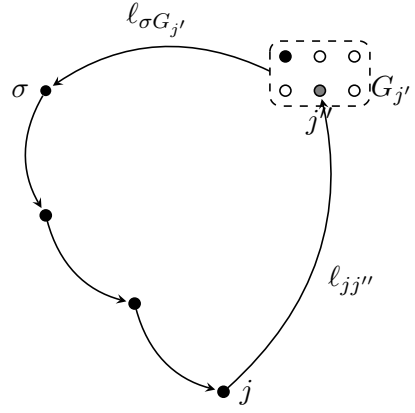
\begin{figure}[!htbp]
	\centering
	\begin{subfigure}{0.4\textwidth}
		\centering
		\begin{tikzpicture}[
			> = stealth, 
			shorten > = 1pt, 
			auto,
			semithick 
			]
			\node[circle, fill, scale=0.4, label=left:$\depot{}$] (depot){} ;
			\node[circle, draw, scale=0.4, right = 3cm of depot] (jndem4){} ;
			\node[circle, draw, scale=0.4, right = 0.3cm of jndem4 ] (jndem5){} ;
			\node[circle, draw, scale=0.4, right = 0.3cm of jndem5,  label=right:$G_{j'}$] (jndem6){} ;
			\node[circle, fill, draw, scale=0.4, above = 0.3cm of jndem4, label = above: $j'$] (jndem1){} ;
			\node[circle, draw, scale=0.4, above = 0.3cm of jndem5] (jndem2){} ;
			\node[circle, draw, scale=0.4, above = 0.3cm of jndem6] (jndem3){} ;
			\draw[dashed, rounded corners]  ($(jndem1.north west)+(-0.12,0.12)$) rectangle ($(jndem6.south east)+(0.12,-0.12)$) {};
			\node[circle, draw, fill, scale=0.4, below = 1.5cm of depot] (j){} ;
			\node[circle, draw, fill, scale=0.4, below right = 1.5cm of j] (gt2){} ;
			\node[circle, draw, fill, scale=0.4, below right = 1.5cm of gt2, label=right:$j$] (gt3){} ;
			\path[bend right, ->, left] (depot) edge  node[] {} (j);
			\path[bend right, ->, left] (j) edge  node[] {} (gt2);
			\path[bend right, ->, left] (gt2) edge  node[] {} (gt3);
			\path[bend right, ->, right] (gt3) edge  node[] {$\ell_{j \depot}$} (depot);	
			
		\end{tikzpicture} 
		\subcaption{Current giant tour ($j$ denotes the last visited node) and group of potential stopping nodes $G_j'$.}
		\label{fig:savings_a}
	\end{subfigure}%
	\begin{subfigure} {0.15\textwidth}
		$\quad$
	\end{subfigure}
	\begin{subfigure}{0.4\textwidth}
		\centering
		\begin{tikzpicture}[
			> = stealth, 
			shorten > = 1pt, 
			auto,
			semithick 
			]
			\node[circle, fill, scale=0.4, label=left:$\depot{}$] (depot){} ;
			\node[circle, draw, scale=0.4, right = 3cm of depot] (jndem4){} ;
			\node[circle,  fill = gray,  draw, scale=0.4, right = 0.3cm of jndem4, label=below:$j''$] (jndem5){} ;
			\node[circle, draw, scale=0.4, right = 0.3cm of jndem5, label=right:$G_{j'}$] (jndem6){} ;
			\node[circle, fill, draw, scale=0.4, above = 0.3cm of jndem4] (jndem1){} ;
			\node[circle, draw, scale=0.4, above = 0.3cm of jndem5] (jndem2){} ;
			\node[circle, draw, scale=0.4, above = 0.3cm of jndem6] (jndem3){} ;
			\draw[dashed, rounded corners ]  ($(jndem1.north west)+(-0.15,0.15)$) rectangle ($(jndem6.south east)+(0.15,-0.15)$) {};
			\node[circle, draw, fill, scale=0.4, below = 1.5cm of depot] (j){} ;
			\node[circle, draw, fill, scale=0.4, below right = 1.5cm of j] (gt2){} ;
			\node[circle, draw, fill, scale=0.4, below right = 1.5cm of gt2, label=right:$j$] (gt3){} ;
			\path[bend right, ->, left] (depot) edge  node[] {} (j);
			\path[bend right, ->, left] (j) edge  node[] {} (gt2);
			\path[bend right, ->, left] (gt2) edge  node[] {} (gt3);
			\path[bend right, ->, above] ($(jndem1.north west)+(-0.15,-0.25)$) edge  node[] {$\ell_{\depot G_{j'}}$} (depot);	
			\path[bend right, ->, below right] (gt3) edge  node[] {$\ell_{j j''}$} (jndem5);	
			
		\end{tikzpicture}
		\subcaption{Savings obtained from inserting node $j'' \in G_{j'}$ after node $j$: $s_{jj''} = \ell_{\depot{}G_{j'}} + \ell_{j \depot} - \ell_{j j''}$.}
		\label{fig:savings_b}
	\end{subfigure}
	\caption{Savings calculation for the insertion of a group of potential stopping nodes after the last node of the current giant tour.}
	\label{fig:savings}
\end{figure}

\Cref{fig:savings} illustrates the calculation of savings for a given group of potential stopping nodes. Imagine the next node is to be inserted after the last node in the giant tour, which is denoted by $j$. We denote by $G_{j'}$ the group defined by the set cover node $j'$ (black node) and its alternatives (white nodes). Given a group, the candidate node to be inserted in the giant tour corresponds to the closest node to node $j$, which is denoted by $j''$ (gray node) in \Cref{fig:savings_b}. The savings obtained from inserting node $j''$ after node $j$ in the giant tour are calculated with respect to the simple tour that would connect the depot with the group, i.e., $s_{jj''} = \ell_{\depot{}G_{j'}} + \ell_{j \depot} - \ell_{j j''}$. As shown in \Cref{fig:savings_b}, $\ell_{\depot{}G_{j'}}$ corresponds to the shortest path connecting the group with the depot that can be saved from being included in the giant tour, $\ell_{j \depot}$ is the link removed from the tour and $\ell_{j j''}$ is the link added to the tour.

\subsection{Second phase} \label{subsec:second}

The goal of the second phase is to solve a SDVRP on $\stoppingNodesSel{}$. The tours performed by the vehicles are constructed in two consecutive steps that enable a simplified treatment of potential split nodes (\ie potential stopping nodes that are visited by multiple vehicles) and vehicle capacities. In the first step, we transform the SDVRP into a CVRP by means of an a priori splitting strategy of the nodes (see \Cref{subsubsec:sdvrp_cvrp}). This allows us to use a state-of-the-art algorithm to solve the resulting CVRP in the second step (see \Cref{subsubsec:hgs_cvrp}).  

\Cref{algo:second_phase} shows the pseudocode of this phase. For each set cover $\stoppingNodesSel{}  \in$ BestSets, we first add it to the list TreatedSetCovers so that the feedback mechanism with the first phase via the penalization of the cost of the associated giant tour can be applied. We then transform the SDVRP into a CVRP and solve it with the selected algorithm. We keep track of the solution with the lowest overall cost ($\text{BestSol}$) and update it if applicable.

\begin{algorithm}[!htbp]
	\caption{Second phase of the two-phased heuristic method for the \probAbbrVRP{}}
	\label{algo:second_phase}
	\small
	\begin{algorithmic}[1]
		
		\State $\text{BestSol} \leftarrow \emptyset$
		\State $\text{totalCost(BestSol)} \leftarrow \infty$ 
		\For{$\stoppingNodesSel{} \in \text{BestSets}$}
		\State $\text{Add } \stoppingNodesSel{}  \text{ to TreatedSetCovers}$ 
		\State $\text{Transform SDVRP associated with } \stoppingNodesSel{}  \text{ into a CVRP}$ 
		\State $\text{CVRPSol } \leftarrow \text{solve resulting CVRP with HGS-CVRP}$ 
		\State $\text{SDVRPSol } \leftarrow \text{transform CVRPSol into a SDVRP solution}$ 
		\If{$\text{totalCost(SDVRPSol)} < \text{totalCost}(\text{BestSol})$}
		\State  $\text{BestSol} \leftarrow \text{SDVRPSol} $
		\EndIf
		\EndFor
	\end{algorithmic}
\end{algorithm}

\subsubsection{From SDVRP to CVRP} \label{subsubsec:sdvrp_cvrp}

To solve the SDVRP, we consider the method developed by \cite{chen2017novel}. The idea is to transform the SDVRP into a CVRP at the expense of an increased number of customers. To do so, an a priori splitting strategy is used, \ie each customer demand is split in advance and not when the problem is being solved. This approach allows to use any CVRP solver instead of developing tailored algorithms for SDVRP. Finally, the solution to the CVRP is transformed into a solution to the original SDVRP.

\cite{chen2017novel} propose a 20/10/5/1 rule that splits each customer demand into $m_{20}$ pieces of $0.2Q$, $m_{10}$ pieces of $0.1Q$, $m_{5}$ pieces of $0.05Q$, $m_{1}$ pieces of $0.01Q$, and at most one remaining piece of less than $0.01Q$, where $Q$ is the vehicle capacity. The number of pieces for each customer $j$ is calculated as follows: $m_{20} = \max \{m \in \mathbb{Z}^+ \cup \{0\} | 0.2Qm \leq d_j\}$, $m_{10} = \max \{m \in \mathbb{Z}^+ \cup \{0\} | 0.1Qm \leq d_j - 0.2Qm_{20}\}$, $m_{5} = \max \{m \in \mathbb{Z}^+ \cup \{0\} | 0.05Qm \leq d_j - 0.2Qm_{20} - 0.1Qm_{10} \}$ and $m_{1} = \max \{m \in \mathbb{Z}^+ \cup \{0\} | 0.05Qm \leq d_j - 0.2Qm_{20} - 0.1Qm_{10} - 0.05Qm_{5} \}$, where $d_j$ is the demand of customer $j$. For instance, if $Q= 100$ and $d_j=76$, this demand is split into six nodes with demands 20, 20, 20, 10, 5, and 1. Notice that the only way not to eliminate the optimal SDVRP solution when both the vehicle capacity and the customer demands are integers is to split the demands into unit demands. Nevertheless, this strategy will dramatically increase the size of the resulting CVRP, and consequently the running time. As pointed out by the authors, the split rule might not eliminate the optimal solution if the vehicles (routes in our case) that serve split demands are not fully loaded, since this involves many different ways to split the demand in an optimal solution.

For the \probAbbrVRP{}, the amount of demand to be satisfied at each stopping node is determined during the first phase, when the nodes that form the giant tour are selected. For these nodes, we only split the ones whose demand exceeds 10\% of the vehicle capacity, i.e., the nodes $j$ such that $d_j \geq 0.1Q$, with the 20/10/5/1 split rule. In this way, we prevent to split small demands (with respect to vehicle capacity), and therefore limit the number of nodes that are created when defining the associated CVRP.

\subsubsection{Hybrid genetic search for the CVRP } \label{subsubsec:hgs_cvrp}

We rely on the hybrid genetic search for the CVRP (HGS-CVRP) developed by \cite{vidal2020hybrid}. This state-of-the-art algorithm uses the same general methodology as \cite{vidal2012hybrid} but includes an additional neighborhood called SWAP* that consists in exchanging two customers between different routes without an insertion in place. Computational experiments have shown that HGS-CVRP stands as a leading metaheuristic regarding solution quality and convergence speed.

The performance of HGS-CVRP comes from the combination of three main strategies. First, a synergistic combination of crossover-based and neighborhood-based search. The former allows a diversified search in the solution space, while the latter enables aggressive solution improvement. Second, a controlled exploration of infeasible solutions in which any excess load in the routes is linearly penalized. Third, population diversity management strategies that allow to maintain a diversified and high-quality set of solutions and counterbalance the loss of diversity due to the neighborhood search. We refer the reader to \cite{vidal2020hybrid} for further details on the HGS-CVRP algorithm.

We apply the open-source implementation of HGS-CVRP (\url{https://github.com/vidalt/HGS-CVRP}) with a lower number of iterations without improvement (10000 instead of 20000 to speed up the process) and the rest of the parameters set to their default values. The algorithm returns a solution file with the nodes visited in each route and the total cost. We then simply need to transform this solution into a solution to the original SDVRP.

\section{Computational experiments}
\label{sec:compTests}

In this section, we present the tests conducted on a set of small instances inspired by \cite{letchford2013compact} and a set of real-life instances from four different municipalities in Switzerland. In~\Cref{sec:data_prob_inst}, we describe the various datasets and problem instances considered for such tests. \Cref{sec:comp_milp_form} compares the MILP formulations against each other, \Cref{sec:perf_heur} assesses the performance of the heuristic method and \Cref{sec:practical_aspects} discusses some practical aspects of the proposed waste collection concept with respect to the state of practice. We refer the reader to the \ref{sec:app_results} for the exhaustive results on the tested datasets and problem instances.

\subsection{Datasets and problem instances} \label{sec:data_prob_inst}

Each problem instance is generated with a road-network dataset containing information on the nodes and the arcs of the underlying graph and by providing values to the problem parameters. The associated customer-based graphs are constructed according to \Cref{subsection:customerbasedForm}. The two main parameters that characterize an instance are the maximum walking distance ($\gamma$) and the number of tours ($m$). The former is considered to determine the preference lists, \ie the potential stopping nodes within a distance less or equal than $\gamma$ ordered in increasing distance, as described in \Cref{sec:wasteCollApp}, and the latter refers to the necessary number of tours to collect all the waste. Notice that for the sake of deriving multiple instances, we assume various values for $m$ and derive the capacity $Q$ of the vehicle accordingly (see \Cref{sec:wasteCollApp}). We also need to make assumptions on the average speed of the vehicle to transform travel distance into travel time in the objective function and on the stop penalty value $r$. We distinguish between the average speed of the vehicle during collection ($\speedColl$) and the average speed from or to the disposal facility ($\speedDep$).

To derive the set of small instances, we consider the datasets generated in \cite{letchford2013compact} for the STSP, which consist of a collection of sparse road-network graphs. More precisely, we consider 30 datasets, each of them characterized by a total number of nodes $n$ and a probability $p$ that a node is a demand node. Note that STSP instances assume only one tour ($m=1$) and no walking distance ($\gamma = 0$). To create additional instances to our problem, we define supplementary values for both parameters. We consider $m=2$ as an additional value for the number of tours, and we define an additional value for the maximum walking distance depending on the number of nodes of each dataset: $\gamma=15$ [meters] for datasets with more than 100 nodes and $\gamma=30$ [meters] for datasets with up to 100 nodes. Hence, $m \in \{1, 2\}$ and $\gamma \in \{0, 15, 30\}$, which results into four  instances per dataset: two values for $m$ and two values for $\gamma$ ($\gamma=0$ and $\gamma \in \{15,30\}$ depending on the number of nodes). The instances are named $\text{STSP}\_n\_p-\gamma-m$, with $n \in \{50, 75, 100, 125, 150, 175, 200\}$ and $p \in \{1/3, 2/3, 3/3\}$, which yields a total of 84 instances.

\begin{figure}[h]
	\centering
	\includegraphics[width=0.5\linewidth]{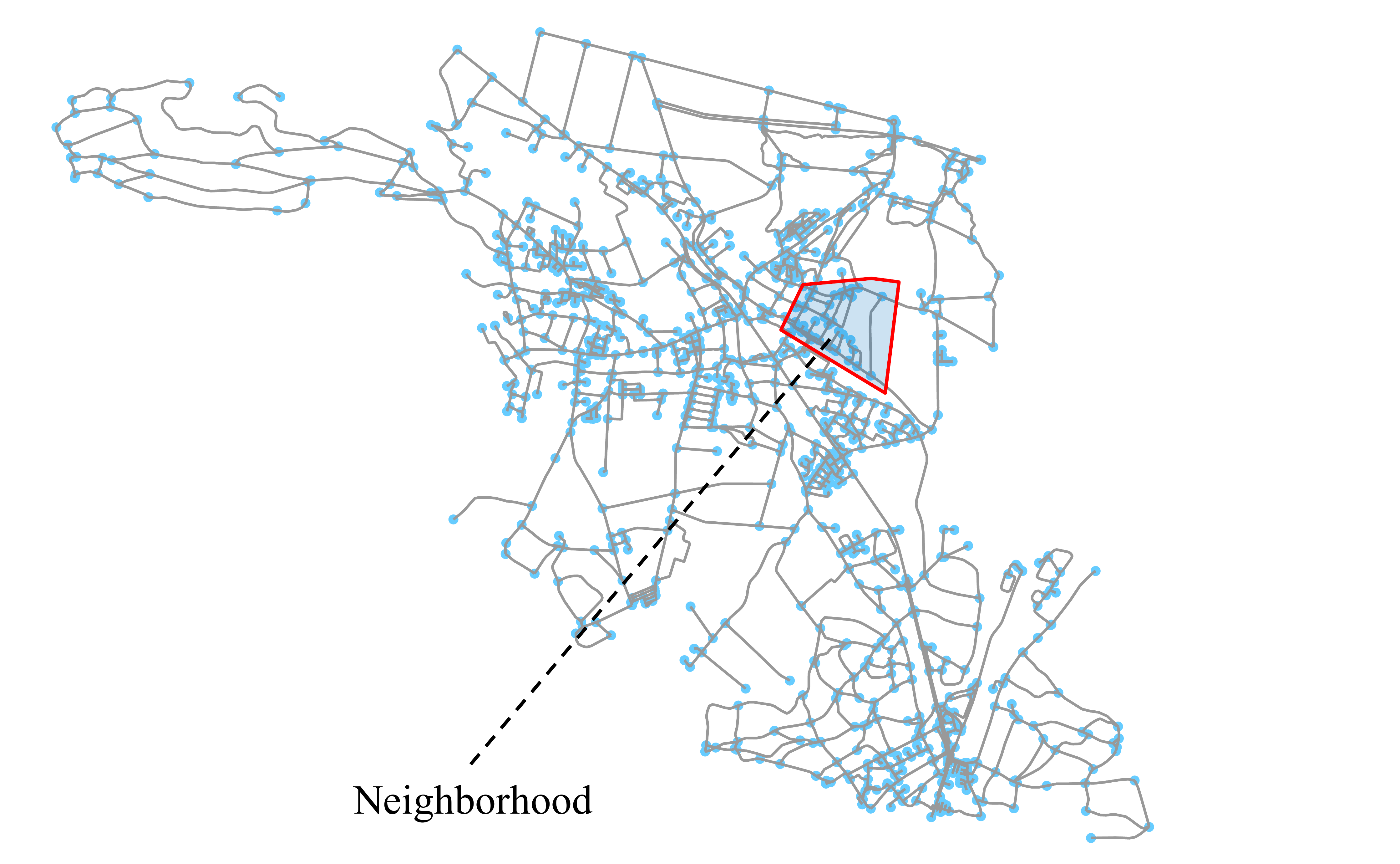}
	\caption{Municipality M3 and its neighborhood (see \Cref{sec:wasteCollApp})}
	\label{fig:munim1}
\end{figure}

The other set of instances is based on real data from four different municipalities in Switzerland, denoted by M1-M4. \Cref{fig:munim1} illustrates municipality M3. The region marked in red corresponds to the neighborhood introduced in \Cref{sec:wasteCollApp}. \Cref{tab:datasets} includes the characteristics of the municipalities with respect to area, number of inhabitants, density and number of residential buildings. Note that the datasets are labeled in increasing order according to size of their area, which is by coincidence also in decreasing order according to the population density. The underlying road-network graph is constructed according to \Cref{sec:wasteCollApp}. \Cref{tab:datasets} also contains the total number of nodes, demand nodes and arcs of the resulting graphs. In a similar way as in the small instances, we consider $m \in \{1,2,4,6\}$ and $\gamma \in \{0, 50, 100, 150\}$. This results into 16~instances for each dataset $\text{M}x$ (with $x \in \{1,2,3,4\}$), which yields 64~instances in total. They are labeled $\text{M}x-\gamma-m$.

Based on the assumed values of $\gamma$ and $m$, the \probAbbrVRP{} is reduced to other less complex problems. More precisely, when only one tour needs to be performed ($m=1$) and $\gamma = 0$, \ie there is no walking distance and consequently all residential buildings need to be visited, the problem becomes a TSP. If $m=1$ and $\gamma > 0$, it turns into a CTP. If more tours are required ($m>1$) and $\gamma = 0$, we face a SDVRP. Finally, when $m>1$ and $\gamma > 0$, the problem can be seen as an $m$-CTP.

Concerning the other parameters, we assume a larger average speed to go from and to the depot for the real-life instances, \ie $\speedColl = 2$~[meters/seconds] and $\speedDep = 14$~[meters/seconds]. We set the stop penalty value $r$ to $\stopTime = 5$~[seconds]. For the small instances, we set $\speedColl = \speedDep = 1$ [meters/seconds], as the STSP in \cite{letchford2013compact} only considers travel distances (and not travel times). Furthermore, there are no assumptions on a stop penalty value at the visited locations. For the sake of illustration, we decide to set it to $\stopTime = 5$~[meters], which can be interpreted as the average physical distance that the vehicle has to travel at a collection point to actually gather the waste.

\begin{table}[h]
	\small
	\centering
	\caption{Characteristics of the municipality datasets M2-M6 and the corresponding road-network graphs}
	\begin{tabular}{lllll}
		\toprule
		Characteristics & M1 & M2 & M3 & M4 \\ 
		\hline
		Municipalities \\
		\hspace{3mm} Area [km$^2$] & 1.52 & 5.96 & 8.81 & 36.32 \\
		\hspace{3mm} Population & 6'776 & 14'708 & 8'721 & 9'839 \\
		\hspace{3mm} Density [inhabitants/km$^2$]  & 4'455.85 & 2'466.52 & 990.39 & 270.92 \\
		\hspace{3mm} Residential buildings & 894 & 1'573 & 1'413 & 1'946  \\
		\hline
		Road-network graphs \\
		\hspace{3mm} Nodes & 652 & 1'365 & 1'353  & 3'137 \\
		\hspace{3mm} Demand nodes & 372 & 722 & 651 & 905 \\
		\hspace{3mm} Arcs & 800 & 1'651 & 1'549  & 3'383 \\
		\bottomrule
	\end{tabular}
	\label{tab:datasets}
\end{table}

The developed computer codes are implemented in Java. The instances were tested on a computer with a 3.4 GHz Intel Core i5 processor, 32 GB of RAM, operating under Windows~10. To solve the MILP formulations, we use the Gurobi 9.1.1 MIP solver via its Java API.

\subsection{Comparison of the MILP formulations} \label{sec:comp_milp_form}

The goal of this section is to test and compare against each other the MILP formulations introduced in \Cref{sec:mip}. We recall that for each network representation, we consider two different formulations with respect to the split collection feature. Note that this feature does not have an impact on the objective function value for instances with $m=1$ (\ie one tour). We set a three-hour time limit for each formulation on each instance. \Cref{tab:exact_overview} presents the number of instances solved to optimality, the number of instances for which a feasible solution was found (excluding the ones solved to optimality) and the number of instances for which no solution was found for each formulation and set of instances. Note that for each formulation and set of instances, the sum of the values of the corresponding row in the table equals the total number of instances (in brackets). We observe that the small instances are all solved by both formulations either to optimality or to feasibility within the time limit, whereas for some real-life instances no solution could be found. For both sets of instances, the road network-based (RN) formulation is able to find more often optimal solutions than its customer-based (CG) counterpart. This is also the case for the real-life instances solved to feasibility. Besides, for more than half of the real-life instances no solution could be found by CG. Concerning split collection, we observe that with this feature enabled fewer optimal and consequently more feasible solutions are obtained for the small instances with both formulations, whereas for the real-life instances, optimality is proven for one more instance with RN and fewer solutions (and fewer optimal ones) are found with CG.

\begin{table}[h]
	\small
	\centering
	\caption{Number of instances solved to optimality, feasibility and not solved by each MILP formulation (total number of instances in brackets).}
	\label{tab:exact_overview}
	\begin{tabular}{ l | c c c  | c c c }
		\toprule
		& \multicolumn{3}{c |}{\textbf{Small instances} (84)} & \multicolumn{3}{ c}{\textbf{Real-life instances} (64)}  \\
		\textbf{Formulation} & Optimal & Feasible & No solution & Optimal & Feasible & No solution \\
		\hline
		RN-wS & 56 & 28 & 0 & 14 & 47 & 3 \\ 
		CG-wS & 42 & 42 & 0 & 1 & 5 & 58 \\ 
		\hline
		RN-noS & 60 & 24 & 0  & 13 & 50 & 1  \\
		CG-noS & 54 & 30 & 0  & 5 & 26 & 33  \\
		\bottomrule
	\end{tabular}
\end{table}

We now take a closer look at the instances for which both formulations could be solved to optimality. \Cref{tab:exact_opt_time} aggregates the instances by problem type based on the values of $m$ and $\gamma$, as detailed in \Cref{sec:data_prob_inst} and includes the number of instances for which each formulation was faster than its counterpart and the average decrease in computational time. For the small instances, all TSP and SDVRP instances could be solved to optimality by both formulations, but SDVRP only when split collection is disabled, whereas fewer optimal solutions can be obtained for CTP and even fewer for $m$-CTP. RN is in general faster than CG, with an average decrease in the computational time that always surpasses 50\%. There are only two exceptions (SDVRP and $m$-CTP) for which there are more instances where CG is faster than RN, in both cases when split collection is disabled. In the case of the real-life instances, we note that only for some TSP and SDVRP instances both formulations could be solved to optimality. For these problem types, RN is always faster, and the average decrease in computational time is huge (almost 100\%) in all cases., \ie RN is extremely fast compared to CG.

\newcommand{\STAB}[1]{\begin{tabular}{@{}c@{}}#1\end{tabular}} 
\begin{table}[t]
	\small
	\centering
	\caption{Comparison of computational times by set of instances, problem type and split feature for the instances for which both formulations were solved to optimality.}
	\label{tab:exact_opt_time}
	\begin{adjustbox}{max width=\textwidth}
		\begin{tabular}{ c | c | c | c | c  | c c | c c }
			\toprule
			& \makecell{Problem \\ type} & \makecell{Total nb \\ instances} & \makecell{Split \\ collection} & \makecell{Nb inst. solved \\ by RN \& CG} & $\text{\makecell{Nb inst. \\ \textbf{RN} faster}}$ & $\text{\makecell{Avg. decrease \\ comp. time (\%)}}$ & $\text{\makecell{Nb inst. \\ \textbf{CG} faster}}$ & $\text{\makecell{Avg. decrease \\ comp. time (\%)}}$ \\
			\cmidrule[\heavyrulewidth](lr){1-9}
			\multirow{8}{*}{\STAB{\rotatebox[origin=c]{90}{Small instances}}} & TSP & 21 & \checkmark & 21 & 14 & 64.3 & 7 & 29.5 \\
			& SDVRP & 21 & \checkmark & 13 & 7 & 82.3 & 6 & 62.9 \\
			& CTP & 21 & \checkmark & 7 & 7 & 83.4 & 0 & - \\
			& $m$-CTP & 21 & \checkmark &  1 & 1 & 87.4 & 0 & - \\
			\cmidrule{2-9}
			& TSP & 21 & $\times$ & 21 & 14 & 56.9 & 7 & 43.3 \\
			& SDVRP & 21 & $\times$ & 21 & 0 & - & 21 & 77.6 \\	
			& CTP & 21 & $\times$ & 11 & 10 & 74.2 & 1 & 12.4 \\
			& $m$-CTP & 21 & $\times$ & 1 & 0 & - & 1 & 12.8 \\
			\cmidrule[\heavyrulewidth](lr){1-9}
			\multirow{8}{*}{\STAB{\rotatebox[origin=c]{90}{Real-life instances}}} & TSP & 4 & \checkmark &1 & 1 & 99.9 & 0 & - \\
			& SDVRP & 12 & \checkmark & 0 & - & - & - & - \\
			& CTP & 12 & \checkmark &0 & - & - & - & - \\
			& $m$-CTP & 36 & \checkmark & 0 & - & - & - & - \\
			\cmidrule{2-9}
			& TSP & 4 & $\times$ & 2 & 2 & 99.9 & 0 & - \\
			& SDVRP & 12 & $\times$ & 1 & 1 & 99.2 & 0 & - \\
			& CTP & 12 & $\times$ & 0 & - & - & - & - \\	
			& $m$-CTP & 36 & $\times$ & 0 & - & - & - & - \\		
			\bottomrule
		\end{tabular}
	\end{adjustbox}
\end{table}

Similar findings can be established with the instances solved to feasibility (excluding those solved to optimality) by both formulations. In this case, we analyze the number of instances for which each formulation could find better solutions than its counterpart and the average decrease in the objective function value. As we can see in \Cref{tab:exact_feas_best_sol}, RN is able to find better (or at least equally good) solutions than CG with enabled split collection for small SDVRP and $m$-CTP instances. The decrease in the objective function value is moderate: 5.4\% for SDVRP and 1.7\% for $m$-CTP. CG is able to find better solutions for more CTP instances, with an average decrease of 1.8\%. With disabled split collection, both formulations find better solutions than their counterpart for the same number of CTP and $m$-CTP instances with an average decrease that also remain moderate (around 3\%). RN outperforms CG in finding better solutions for the real-life instances, with average improvements that span from 7.5\% (CTP, split collection enabled) to 77.3\% (SDVRP, split collection disabled).  

\begin{table}[t]
	\small
	\centering
	\caption{Comparison of optimal values of the objective function by set of instances, problem type and split feature for the instances for which both formulations were solved to feasibility (excluding the ones solved to optimality).}
	\label{tab:exact_feas_best_sol}
	\begin{threeparttable}
		\begin{adjustbox}{max width=\textwidth}
			\begin{tabular}{ c | c | c | c | c  | c c | c c }
				\toprule
				& \makecell{Problem \\ type} & \makecell{Total nb \\ instances} & \makecell{Split \\ collection} & \makecell{Nb inst. solved \\ by RN \& CG} & $\text{\makecell{Nb inst. \\ \textbf{RN} better}}$ & $\text{\makecell{Avg. decrease \\ obj. value (\%)}}$ & $\text{\makecell{Nb inst. \\ \textbf{CG} better}}$ & $\text{\makecell{Avg. decrease \\ obj. value (\%)}}$ \\
				\cmidrule[\heavyrulewidth](lr){1-9}
				\multirow{8}{*}{\STAB{\rotatebox[origin=c]{90}{Small instances}}} & TSP & 21 & \checkmark & 0 & - & - & - & - \\
				& SDVRP & 21 & \checkmark & 3 & 2\tnote{a} & 5.4 & 0\tnote{a} & - \\
				& CTP & 21 & \checkmark &7 & 2\tnote{a} & 0.1 & 3\tnote{a} & 1.8 \\
				& $m$-CTP & 21 & \checkmark &18 & 16\tnote{a} & 1.7 & 1\tnote{a} & 1.2 \\
				\cmidrule{2-9}	
				& TSP & 21 & $\times$ &0 & - & - & - & - \\
				& SDVRP & 21 & $\times$ & 0 & - & - & - & - \\	
				& CTP & 21 & $\times$ &6 & 2\tnote{a} & 0.7 & 2\tnote{a} & 2.9 \\
				& $m$-CTP & 21 & $\times$ &18 & 6\tnote{a} & 0.7 & 6\tnote{a} & 3.0 \\
				\cmidrule[\heavyrulewidth](lr){1-9}
				\multirow{8}{*}{\STAB{\rotatebox[origin=c]{90}{Real-life instances}}} & TSP & 4 & \checkmark & 0 & - & - & - & - \\
				& SDVRP & 12 & \checkmark  & 0 & - & - & - & - \\
				& CTP & 12 & \checkmark & 1 & 1 & 7.5 & 0 & - \\
				& $m$-CTP & 36 & \checkmark & 0 & - & - & - & - \\
				\cmidrule{2-9}			
				& TSP & 4 & $\times$ &0 & - & - & - & - \\
				& SDVRP & 12 & $\times$ &5 & 5 & 77.3 & 0 & - \\
				& CTP & 12 & $\times$ & 1 & 1 & 44.5 & 0 & - \\	
				& $m$-CTP & 36 & $\times$ &  13 & 11 & 32.5 & 2 & 12.7 \\		
				\bottomrule
			\end{tabular}
		\end{adjustbox}
		\begin{tablenotes}\footnotesize
			\item[a] The sum of the instances for which RN was better than CG and CG was better than RN is not equal to the number \\ of instances solved by both formulations because the ones for which both RN and CG have the same performance \\ (\ie the same objective value is obtained) are not included in this table. 
		\end{tablenotes}
		
	\end{threeparttable}
	
\end{table}

These experiments illustrate that RN is preferred both from a computational and solution quality point of view. Regarding split collection, the disabling of this feature appears to contribute to both formulations in generating a larger number of optimal solutions for the small instances (especially to CG) and to CG in providing more often optimal and feasible solutions for the real-life instances. Altogether, we conclude that RN should be given preference over CG. Split collection is typically a good feature in practice and therefore, we resolve to consider RN-wS for the upcoming experiments.

We now analyze the impact of the parameters that define the problem instances ($m$ and $\gamma$) on the solvability of the models for RN-wS. To this end, Table~\ref{tab:exact_parameters} reports the number of instances that could be solved (to optimality) and the ones that could not be solved at all. The coloring of the cells allows to quickly identify the type of problem the \probAbbrVRP{} is reduced to for the considered values of $m$ and $\gamma$. Note that for the small instances, $\gamma = 15$ or $\gamma=30$ depending on the number of nodes in the associated dataset. For both datasets, all the TSP instances can be solved to optimality, in addition, a feasible solution is reported for all small instances.
For the real-life ones, SDVRP instances with $m=2$ and CTP instances with $\gamma \in \{50,100\}$ can be solved to optimality, whereas for the SDVRP ones with more tours and the CTP ones with a larger maximum walking distance only a feasible solution is obtained. For three real-life $m$-CTP instances, no solution could be found, and for the rest, none of them could be solved to optimality. Hence, it becomes apparent that as $\gamma$ and/or $m$ increase, RN-wS struggles in finding the optimal solution, and is only able to provide a feasible solution or no solution at all.

\newcolumntype{C}{>{\centering\arraybackslash}p{1.14em}} 

\begin{table}[t]
	
	\caption{Number of instances that could be solved to optimality, to feasibility (excluding the ones solved to optimality) and not solved by RN-wS with respect to $m$ (number of tours) and $\gamma$ (maximum walking distance). In brackets, the total number of instances for each combination of parameters is specified.}
	\vspace{-1.5em}
	\begin{center}
		\begin{tabular}{ll ll ll ll}
			\textcolor{tsp}{$\blacksquare$} & TSP & \textcolor{sdvrp}{$\blacksquare$} & SDVRP & \textcolor{ctp}{$\blacksquare$} & CTP & \textcolor{mctp}{$\blacksquare$} & $m$-CTP  \\
		\end{tabular}  
	\end{center}
	
	\begin{adjustbox}{max width=\textwidth}
		\begin{subtable}[h]{0.34\textwidth}
			\centering
			\footnotesize
			\begin{tabular}{ l | C C | C C | C C }
				\toprule
				& \multicolumn{2}{c|}{\footnotesize Optimal} & \multicolumn{2}{c|}{\footnotesize Feasible} & \multicolumn{2}{c}{\footnotesize No sol.} \\ 
				\diagbox[innerwidth=1.1cm,height=1.15\line]{$\gamma$}{$m$}& \footnotesize 1 & \footnotesize 2 & \footnotesize 1 & \footnotesize 2 & \footnotesize 1 & \footnotesize 2 \\
				\hline
				\footnotesize 15 / 30 & \scriptsize 14 \cellcolor{ctp} & \scriptsize 3 \cellcolor{mctp} & \scriptsize 7 \cellcolor{ctp} & \scriptsize 18 \cellcolor{mctp}& \scriptsize 0 \cellcolor{ctp} & \scriptsize 0 \cellcolor{mctp}  \\
				\footnotesize 0 & \scriptsize 21 \cellcolor{tsp} & \scriptsize 18 \cellcolor{sdvrp} & \scriptsize 0 \cellcolor{tsp} & \scriptsize 3 \cellcolor{sdvrp} & \scriptsize 0 \cellcolor{tsp} & 0 \cellcolor{sdvrp}    \\
				\hline
			\end{tabular}
			\caption{Small instances \footnotesize(21 per cell)}
			\label{tab:ex_par_small}
		\end{subtable}
		\hspace{1em}
		\begin{subtable}[h]{0.575\textwidth}
			\centering
			\footnotesize
			\begin{tabular}{l | C C C C | C C C C | C C C C}
				\toprule
				& \multicolumn{4}{c|}{\footnotesize Optimal} & \multicolumn{4}{c|}{\footnotesize Feasible}& \multicolumn{4}{c}{\footnotesize No solution} \\ 
				\diagbox[innerwidth=0.75cm,height=1.25\line]{$\gamma$}{$m$}& \footnotesize 1 & \footnotesize 2 & \footnotesize 4 & \footnotesize 6 & \footnotesize 1 & \footnotesize 2 & \footnotesize 4 & \footnotesize 6 & \footnotesize 1 & \footnotesize 2 & \footnotesize 4 & \footnotesize 6 \\
				\hline
				\footnotesize 150 & \scriptsize 0 \cellcolor{ctp} &\scriptsize 0 \cellcolor{mctp} & \scriptsize 0 \cellcolor{mctp} & \scriptsize 0 \cellcolor{mctp} & \scriptsize 4 \cellcolor{ctp} & \scriptsize 4 \cellcolor{mctp} &  \scriptsize 4 \cellcolor{mctp}& \scriptsize 3 \cellcolor{mctp} & \scriptsize 0 \cellcolor{ctp} & \scriptsize 0 \cellcolor{mctp} &  \scriptsize 0 \cellcolor{mctp}& \scriptsize 1 \cellcolor{mctp}\\
				\footnotesize 100 & \scriptsize 2 \cellcolor{ctp}  &\scriptsize 0 \cellcolor{mctp} & \scriptsize 0 \cellcolor{mctp} & \scriptsize 0 \cellcolor{mctp} & \scriptsize 2 \cellcolor{ctp} & \scriptsize 4 \cellcolor{mctp} &  \scriptsize 3 \cellcolor{mctp}& \scriptsize 4 \cellcolor{mctp} & \scriptsize 0 \cellcolor{ctp} & \scriptsize 0 \cellcolor{mctp} &  \scriptsize 1 \cellcolor{mctp}& \scriptsize 0 \cellcolor{mctp} \\
				\footnotesize 50 & \scriptsize 4 \cellcolor{ctp} & \scriptsize 0 \cellcolor{mctp} & \scriptsize 0 \cellcolor{mctp} & \scriptsize 0 \cellcolor{mctp} & \scriptsize 0 \cellcolor{ctp} & \scriptsize 4 \cellcolor{mctp} & \scriptsize 3 \cellcolor{mctp}& \scriptsize 4 \cellcolor{mctp} & \scriptsize 0 \cellcolor{ctp} & \scriptsize 0 \cellcolor{mctp} &  \scriptsize 1 \cellcolor{mctp}& \scriptsize 0 \cellcolor{mctp}\\
				\footnotesize 0 & \scriptsize 4 \cellcolor{tsp}  & \scriptsize 4 \cellcolor{sdvrp} & \scriptsize 0 \cellcolor{sdvrp} & \scriptsize 0 \cellcolor{sdvrp} & \scriptsize 0 \cellcolor{tsp}  & \scriptsize 0 \cellcolor{sdvrp} & \scriptsize 4 \cellcolor{sdvrp} & \scriptsize 4 \cellcolor{sdvrp} & \scriptsize 0 \cellcolor{tsp}  & \scriptsize 0 \cellcolor{sdvrp} & \scriptsize 0 \cellcolor{sdvrp} & \scriptsize 0 \cellcolor{sdvrp} \\
				\hline
			\end{tabular}
			\caption{Real-life instances \footnotesize(4 per cell)}
			\label{tab:ex_par_large}
		\end{subtable}
	\end{adjustbox}
	\label{tab:exact_parameters}
	
\end{table}

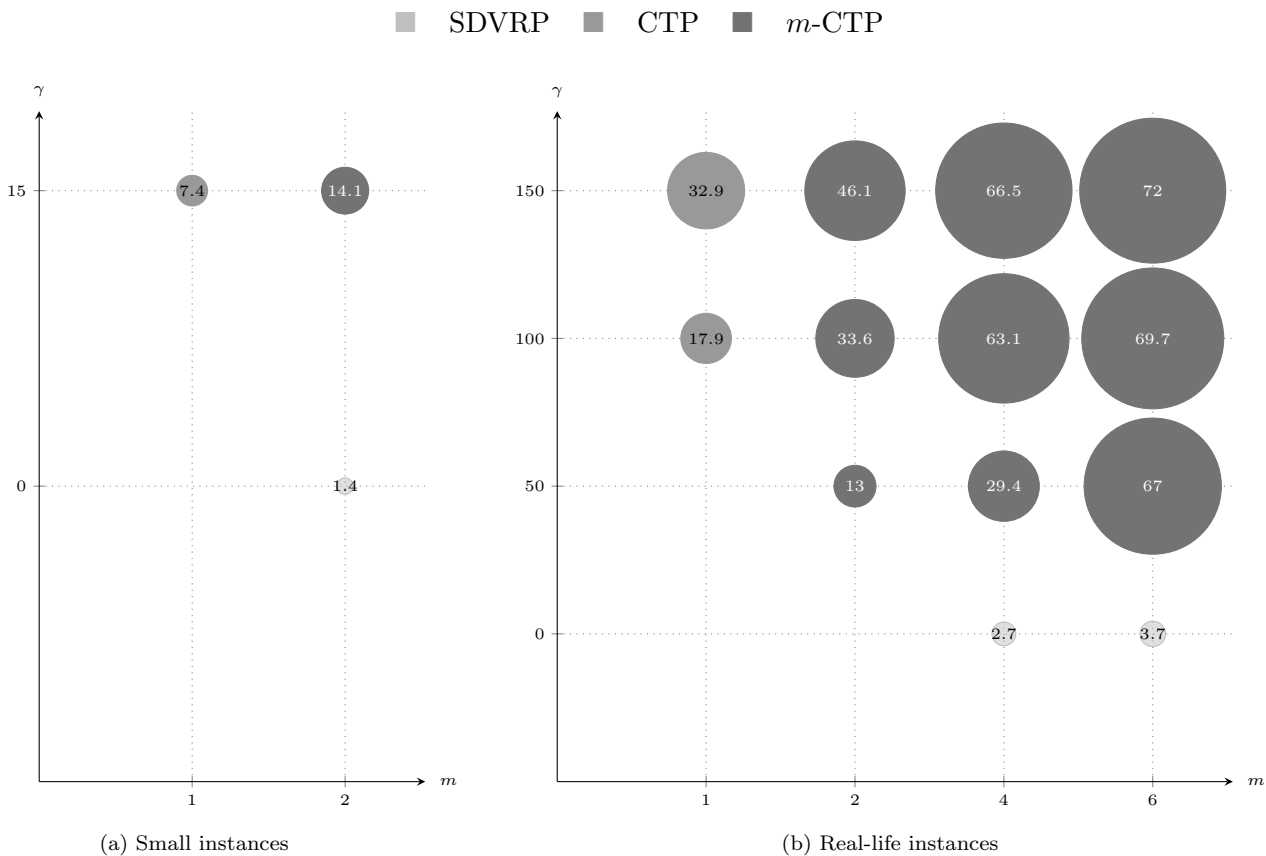
\begin{figure}[t]
	\begin{center}
		\begin{tabular}{ ll ll ll}
			\textcolor{sdvrp}{$\blacksquare$} & SDVRP & \textcolor{ctp}{$\blacksquare$} & CTP & \textcolor{mctp}{$\blacksquare$} & $m$-CTP  \\
		\end{tabular}  
	\end{center}
	\centering
	\begin{subfigure}[b]{0.3\textwidth}
		\centering
		\tiny
		\begin{tikzpicture}
			\begin{axis}[
				width=2in,
				height=3.5in,
				xlabel={$m$},
				ylabel={$\gamma$},
				grid = both, grid style={dotted,black!50},
				xtick={x,1,2},
				xmin={x}, xmax={2},
				symbolic x coords = {x, 1, 2},
				symbolic y coords = {y, 0, 15},
				ytick distance=1,ymin = y, ymax=15,
				scale only axis,
				axis lines=center,
				every axis x label/.style={at={(current axis.south east)},right=1mm},
				every axis y label/.style={at={(current axis.north west)},above=1mm},
				enlargelimits={abs=30pt,upper},
				]
				\addplot[%
				scatter=true,
				only marks,
				mark=*,
				point meta=explicit,
				fill opacity=0.5,text opacity=1,
				visualization depends on = {2.25*\thisrow{Yval} \as \perpointmarksize},
				scatter/use mapped color = {fill=sdvrp, draw = sdvrp},
				scatter/@pre marker code/.append style={
					/tikz/mark size=\perpointmarksize},
				nodes near coords*,
				nodes near coords style={text=black,font=\sffamily,anchor=center},
				] table [x={x},y={Algorithm},meta index=2] {
					x  Algorithm Val Yval 
					2 0 1.4 1.37
				};     
				
				\addplot[%
				scatter=true,
				only marks,
				mark=*,
				point meta=explicit,
				fill opacity=1,text opacity=1,
				visualization depends on = {2.25*\thisrow{Yval} \as \perpointmarksize},
				scatter/use mapped color = {fill=ctp, draw = ctp},
				scatter/@pre marker code/.append style={
					/tikz/mark size=\perpointmarksize},
				nodes near coords*,
				nodes near coords style={text=black,font=\sffamily,anchor=center},
				] table [x={x},y={Algorithm},meta index=2] {
					x  Algorithm Val Yval 
					1 15 7.4 2.57
				};     
				
				\addplot[%
				scatter=true,
				only marks,
				mark=*,
				point meta=explicit,
				fill opacity=1,text opacity=1,
				visualization depends on = {2.25*\thisrow{Yval} \as \perpointmarksize},
				scatter/use mapped color = {draw = mctp, fill=mctp},
				scatter/@pre marker code/.append style={
					/tikz/mark size=\perpointmarksize},
				nodes near coords*,
				nodes near coords style={text=white,font=\sffamily,anchor=center},
				] table [x={x},y={Algorithm},meta index=2] {
					x  Algorithm Val Yval 
					2 15 14.1 3.92
				};   
			\end{axis}
		\end{tikzpicture}
		\caption{Small instances}
		\label{fig:exact_feas_gaps_small}
	\end{subfigure}
	\hspace{3.25em}
	\begin{subfigure}[b]{0.6\textwidth}
		\centering
		\tiny
		\begin{tikzpicture}
			\begin{axis}[
				width=3.5in,
				height=3.5in,
				xlabel={$m$},
				ylabel={$\gamma$},
				grid = both, grid style={dotted,black!50},
				xtick={x,1,2,4,6},
				xmin={x}, xmax={6},
				symbolic x coords = {x, 1, 2, 4, 6},
				symbolic y coords = {y, 0, 50, 100, 150},
				ytick distance=1,ymin = y, ymax=150,
				scale only axis,
				axis lines=center,
				every axis x label/.style={at={(current axis.south east)},right=1mm},
				every axis y label/.style={at={(current axis.north west)},above=1mm},
				enlargelimits={abs=30pt,upper},
				]
				\addplot[%
				scatter=true,
				only marks,
				mark=*,
				point meta=explicit,
				fill opacity=0.5,text opacity=1,
				visualization depends on = {1.65*\thisrow{Yval} \as \perpointmarksize},
				scatter/use mapped color = {fill=sdvrp, draw = sdvrp},
				scatter/@pre marker code/.append style={
					/tikz/mark size=\perpointmarksize},
				nodes near coords*,
				nodes near coords style={text=black,font=\sffamily,anchor=center},
				] table [x={x},y={Algorithm},meta index=2] {
					x  Algorithm Val Yval 
					4 0 2.7 2.7
					6 0 3.7 2.9
				};     
				
				\addplot[%
				scatter=true,
				only marks,
				mark=*,
				point meta=explicit,
				fill opacity=1,text opacity=1,
				visualization depends on = {1.65*\thisrow{Yval} \as \perpointmarksize},
				scatter/use mapped color = {fill=ctp, draw = ctp},
				scatter/@pre marker code/.append style={
					/tikz/mark size=\perpointmarksize},
				nodes near coords*,
				nodes near coords style={text=black,font=\sffamily,anchor=center},
				] table [x={x},y={Algorithm},meta index=2] {
					x  Algorithm Val Yval 
					1 100 17.9 5.74
					1 150 32.9 8.74
				};     
				
				\addplot[%
				scatter=true,
				only marks,
				mark=*,
				point meta=explicit,
				fill opacity=1,text opacity=1,
				visualization depends on = {1.65*\thisrow{Yval} \as \perpointmarksize},
				scatter/use mapped color = {draw = mctp, fill=mctp},
				scatter/@pre marker code/.append style={
					/tikz/mark size=\perpointmarksize},
				nodes near coords*,
				nodes near coords style={text=white,font=\sffamily,anchor=center},
				] table [x={x},y={Algorithm},meta index=2] {
					x  Algorithm Val Yval 
					2 50 13 4.76
					4 50 29.4 8.04
					6 50 67 15.56
					2 100 33.6 8.88
					4 100 63.1 14.78
					6 100 69.7 16.1
					2 150 46.1 11.38
					4 150 66.5 15.46
					6 150 72 16.56
				};   
			\end{axis}
		\end{tikzpicture}
		\caption{Real-life instances}
		\label{fig:exact_feas_gaps_real}
	\end{subfigure}
	\caption{Average optimality gaps for the instances solved to feasibility (excluding the ones solved to optimality) by RN-wS with respect to $\gamma$ and $m$.}
	\label{fig:exact_feas_gaps}
\end{figure}

Finally, we evaluate the optimality gaps reported by Gurobi at the end of the run for the instances solved to feasibility with RN-wS with the bubble chart in \Cref{fig:exact_feas_gaps}. As expected, the more complex the problem type associated with the instances, the larger the average optimality gaps. This behavior is particularly recognizable in the $m$-CTP real-life instances, since the bubbles gradually become bigger as both $\gamma$ and/or $m$ increase.

\subsection{Validation and performance of the heuristic method} \label{sec:perf_heur}

This section aims to validate the heuristic method with respect to the instances solved to optimality by RN-wS and to assess its performance for other instances. As opposed to the exact method, the heuristic approach is always able to find a feasible solution. Note that, in this section we exclude instances for which the MILP formulation could not find any solution.

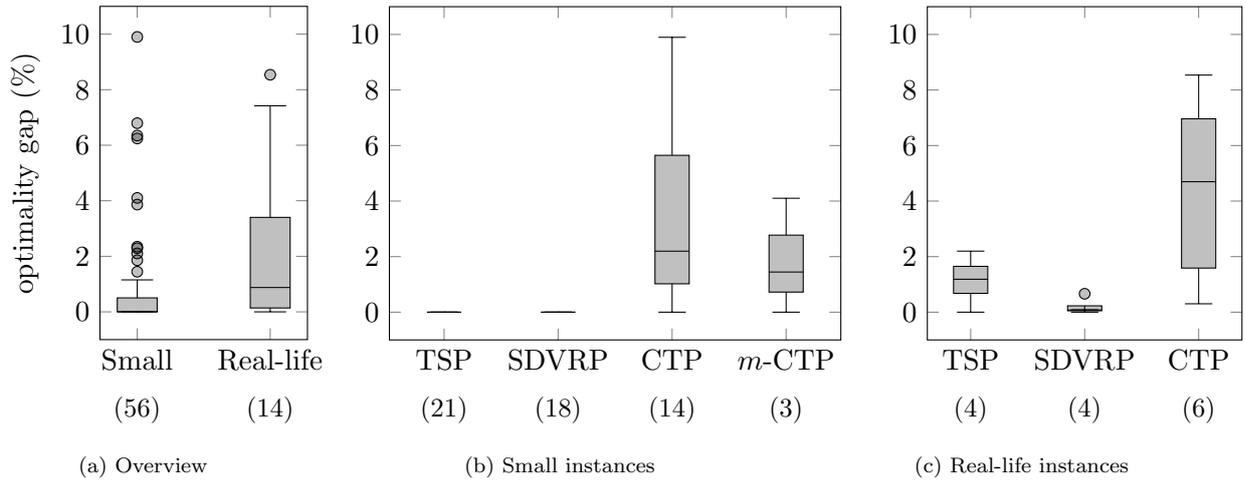
\begin{figure}[t]
	\centering
	\begin{subfigure}[b]{0.23\textwidth}
		\centering
		\begin{tikzpicture}
			\begin{axis}[
				boxplot/draw direction=y,
				ylabel={optimality gap (\%)},
				height=6cm,
				ymin=-1,
				ymax = 11,
				boxplot={
					box extend=0.3,
				},
				x=1.75cm,
				ytick={0,2,4,6,8,10},
				xtick={1,2},
				xticklabels={
					{Small\\ \small (56)},
					{Real-life\\ \small (14)},
				},
				x tick label style={
					text width=3.5cm,
					align=center
				},
				every axis plot/.append style={fill, fill opacity=0.25}
				]
				
				\addplot [mark=*, boxplot]
				table[row sep=\\,y index=0] {
					data\\
					0	\\
					0	\\
					0	\\
					0	\\
					0	\\
					0	\\
					0	\\
					0	\\
					0	\\
					0	\\
					0	\\
					0	\\
					0	\\
					0	\\
					0	\\
					0	\\
					0	\\
					0	\\
					0	\\
					0	\\
					0	\\
					0	\\
					0	\\
					0	\\
					0	\\
					0	\\
					0	\\
					0	\\
					0	\\
					0	\\
					0	\\
					0	\\
					0	\\
					0	\\
					0	\\
					0	\\
					0	\\
					0	\\
					0	\\
					0	\\
					0	\\
					0.452488688	\\
					0.653594771	\\
					0.977777778	\\
					1.150747986	\\
					1.445347787	\\
					1.856763926	\\
					2.103250478	\\
					2.288329519	\\
					2.345582486	\\
					3.861788618	\\
					4.102096627	\\
					6.242638398	\\
					6.346328196	\\
					6.793802145	\\
					9.898107715	\\
				};			
				\addplot [mark=*, boxplot]
				table[row sep=\\,y index=0] {
					data\\
					0	\\
					0	\\
					0.065478951	\\
					0.085042809	\\
					0.302081659	\\
					0.661955551	\\
					0.848899499	\\
					0.902801513	\\
					1.467410644	\\
					2.197285332	\\
					3.803329325	\\
					5.591149778	\\
					7.420013239	\\
					8.537213831	\\
				};
				
			\end{axis}
		\end{tikzpicture}
		\caption{Overview}
		\label{fig:val_heur_over}
	\end{subfigure}
	\begin{subfigure}[b]{0.4\textwidth}
		\centering
		\begin{tikzpicture}
			\begin{axis}[
				boxplot/draw direction=y,
				height=6cm,
				ymin=-1,
				ymax=11,
				boxplot={
					box extend=0.3,
				},
				x=1.5cm,
				ytick={0,2,4,6,8,10},
				xtick={1,2,3,4},
				xticklabels={
					{TSP\\ \small (21)},
					{SDVRP\\ \small (18)},
					{CTP\\ \small (14)},
					{$m$-CTP\\ \small (3)},
				},
				x tick label style={
					text width=3.5cm,
					align=center
				},
				every axis plot/.append style={fill, fill opacity=0.25}
				]
				
				\addplot [mark=*, boxplot]
				table[row sep=\\,y index=0] {
					data\\
					0	\\
					0	\\
					0	\\
					0	\\
					0	\\
					0	\\
					0	\\
					0	\\
					0	\\
					0	\\
					0	\\
					0	\\
					0	\\
					0	\\
					0	\\
					0	\\
					0	\\
					0	\\
					0	\\
					0	\\
					0	\\
				};
				
				\addplot [mark=*, boxplot]
				table[row sep=\\,y index=0] {
					data\\
					0	\\
					0	\\
					0	\\
					0	\\
					0	\\
					0	\\
					0	\\
					0	\\
					0	\\
					0	\\
					0	\\
					0	\\
					0	\\
					0	\\
					0	\\
					0	\\
					0	\\
					0	\\
				};
				\addplot [mark=*, boxplot]
				table[row sep=\\,y index=0] {
					data\\
					0	\\
					0.452488688	\\
					0.653594771	\\
					0.977777778	\\
					1.150747986	\\
					1.856763926	\\
					2.103250478	\\
					2.288329519	\\
					2.345582486	\\
					3.861788618	\\
					6.242638398	\\
					6.346328196	\\
					6.793802145	\\
					9.898107715	\\
				};
				
				\addplot [mark=*, boxplot]
				table[row sep=\\,y index=0] {
					data\\
					0	\\
					1.445347787	\\
					4.102096627	\\
				};
			\end{axis}			
		\end{tikzpicture}
		\caption{Small instances}
		\label{fig:val_heur_small}
	\end{subfigure}
	\begin{subfigure}[b]{0.3\textwidth}
		\centering
		\begin{tikzpicture}
			\begin{axis}[
				boxplot/draw direction=y,
				height=6cm,
				ymin=-1,
				ymax=11,
				boxplot={
					box extend=0.3,
				},
				x=1.5cm,
				ytick={0,2,4,6,8,10},
				xtick={1,2,3},
				xticklabels={
					{TSP\\ \small (4)},
					{SDVRP\\ \small (4)},
					{CTP\\ \small (6)},
				},
				x tick label style={
					text width=3.5cm,
					align=center
				},
				every axis plot/.append style={fill, fill opacity=0.25}
				]
				
				\addplot [mark=*, boxplot]
				table[row sep=\\,y index=0] {
					data\\
					0	\\
					0.902801513	\\
					1.467410644	\\
					2.197285332	\\
				};
				\addplot [mark=*, boxplot]
				table[row sep=\\,y index=0] {
					data\\
					0	\\
					0.065478951	\\
					0.085042809	\\
					0.661955551	\\
				};
				
				\addplot [mark=*, boxplot]
				table[row sep=\\,y index=0] {
					data\\
					0.302081659	\\
					0.848899499	\\
					3.803329325	\\
					5.591149778	\\
					7.420013239	\\
					8.537213831	\\
				};
			\end{axis}			
		\end{tikzpicture}
		\caption{Real-life instances}
		\label{fig:val_heur_real}
	\end{subfigure}
	\caption{Overview of the optimality gaps (with respect to instances solved to optimality by RN-wS) for small and real-life instances (a) and optimality gaps per problem type for small (b) and real-life instances (c). The number of considered instances for each box plot is included in brackets.}
	\label{fig:val_heur}
\end{figure}

We first consider the instances solved to optimality by RN-wS in \Cref{fig:val_heur}. \Cref{fig:val_heur_over} shows an overview of the optimality gaps for both small and real-life instances. The optimality gap is calculated as the relative difference of the objective function value obtained with the heuristic method with respect to the corresponding optimal value. The first box plot includes 56 out of the 84 small instances. For $75\%$ of them, the optimality gaps lie between $0\%$ and $0.5\%$. The median value is $0\%$, which means that the heuristic method was able to find an optimal solution for half of the instances. Note that several outliers are identified above the third quartile.
Only 14 out of the 64 real-life instances were solved to optimality. The associated gaps range from $0\%$ up to a maximum of $8.5\%$ (outlier, $\text{M}4-50-1$). In this case, for $75\%$ of the instances, the associated gaps are below $3.5\%$, and the median value is $0.9\%$.
These numbers are quite low. Therefore, we now take a closer look at each of the problem types. \Cref{fig:val_heur_small} shows that for both TSP and SDVRP small instances, the heuristic method is able to find the optimal solution. The range of the optimality gaps for CTP is larger than for $m$-CTP (higher median and maximum), but we note that most of the outliers seen in \Cref{fig:val_heur_over} for small instances correspond to CTP, and only 3 instances belong to $m$-CTP. As we can see in \Cref{fig:val_heur_real} for real-life instances, the optimality gaps for both TSP and SDVRP instances remain below $2.5\%$, whereas as expected they are larger for CTP (median value around $5\%$). As with the small instances, the outlier identified in \Cref{fig:val_heur_over} is associated with CTP.

\Cref{tab:heur_val_small} presents the value of the objective function obtained by the exact and heuristic methods for the instances inspired by \cite{letchford2013compact} that have not been adapted to our case, \ie the small TSP instances ($\gamma=0$ and $m=1$). For the heuristic method, the last time an improvement was made is as well included. Since the algorithm checks every second whether a better solution has been generated in the different routers that process the set covers from the first phase, this value can be interpreted as the time the algorithm found the best solution within the time limit (termination criterion). For this subset of small instances, the time limit is equal to 60 s. We observe that the heuristic is able to find an optimal solution for all these instances quickly, for all of them within 10 s. Thus, we could replace the termination criterion by a number of checks without improvement to speed up the heuristic approach and get a better sense of the time required to obtain the best solution. Although it is not the purpose of this experiment, we note that the computational times of the exact method is comparable to the ones reported in \citet{letchford2013compact} for these instances.

\begin{table}[h]
	\small
	\centering
	\caption{Computational results of the exact and heuristic methods for the small TSP instances (\ie $\gamma=0$ and $m=1$) solved to optimality (all of them solved to optimality by both methods).}
	\label{tab:heur_val_small}
	\begin{threeparttable}
		\begin{adjustbox}{max width=\textwidth}
			\begin{tabular}{ c | S[table-format=4.0] | S[table-format=2.2] S[table-format=1.0] | c | S[table-format=4.0] | S[table-format=2.2] S[table-format=1.0] }
				\toprule
				& & \multicolumn{2}{c|}{\textbf{Comp. time (s)}} & &  &  \multicolumn{2}{c}{\textbf{Comp. time (s)}} \\
				\makecell{Instance (STSP-$n$-$p$-$\gamma$-$m$)}  & $\text{\makecell{Objective \\ function}}$ & $\text{\makecell{Exact \\ (RN-wS)}}$  & $\text{\makecell{Heuristic \\ method\tnote{a}} }$ & \makecell{Instance (STSP-$n$-$p$-$\gamma$-$m$)}  & $\text{\makecell{Objective \\ function}}$ & $\text{\makecell{Exact \\ (RN-wS)}}$  & $\text{\makecell{Heuristic \\ method\tnote{a}}}$ \\
				\midrule
				$\text{STSP}\_{50}\_1/3-0-1$ & 789 & 1.1  & 2 & $\text{STSP}\_{125}\_3/3-0-1$  & 1582 & 26.6  & 6 \\
				$\text{STSP}\_{50}\_2/3-0-1$ & 978 & 1.7  & 2 & $\text{STSP}\_{150}\_1/3-0-1$ & 1105 & 7.2  & 3\\
				$\text{STSP}\_{50}\_3/3-0-1$ & 1107 & 3.2  & 3 &  $\text{STSP}\_{150}\_2/3-0-1$ & 1545 & 42.2  & 5\\
				$\text{STSP}\_{75}\_1/3-0-1$& 830 & 1.9  & 2 & $\text{STSP}\_{150}\_3/3-0-1$ & 1733 & 32.4  & 1 \\
				$\text{STSP}\_{75}\_2/3-0-1$ & 1029 & 6.1  & 1 & $\text{STSP}\_{175}\_1/3-0-1$ & 1272 & 8.2  & 3\\
				$\text{STSP}\_{75}\_2/3-0-1$ & 1251 & 2.0  & 3 & $\text{STSP}\_{175}\_2/3-0-1$  & 1645 & 30.6  & 5\\
				$\text{STSP}\_{100}\_1/3-0-1$ & 919 & 6.9  & 2 & $\text{STSP}\_{175}\_3/3-0-1$  & 1891 & 24.5  & 9  \\
				$\text{STSP}\_{100}\_2/3-0-1$ & 1193 & 7.4  & 1 & $\text{STSP}\_{200}\_1/3-0-1$  & 1295 & 8.9  & 1 \\			
				$\text{STSP}\_{100}\_3/3-0-1$ & 1445 & 10.0  & 5 & $\text{STSP}\_{200}\_2/3-0-1$  & 1845 & 26.7  & 6 \\			
				$\text{STSP}\_{125}\_1/3-0-1$  & 1143 & 4.5  & 2 & $\text{STSP}\_{200}\_3/3-0-1$  & 2018 & 15.9  & 7 \\
				$\text{STSP}\_{125}\_2/3-0-1$  & 1416 & 20.6  & 1 \\
				\bottomrule
			\end{tabular}
		\end{adjustbox}
		\begin{tablenotes}\footnotesize
			\item[a] This computational time corresponds to the time when reaching the best known solution after which \\ no improvement was made
		\end{tablenotes}
	\end{threeparttable}
	
\end{table}

\Cref{tab:heur_val_real} presents the 14 real-life instances solved to optimality by RN-wS. The optimality gaps are below 5\% except for 3 instances, all of them related to CTP (namely, $\text{M}2-$, $\text{M}3-$ and $\text{M}4-50-1$). These instances correspond to the most complex ones RN-wS is able to optimally solve. Note there are no $m$-CTP instances solved to optimality by RN-wS within the three-hours time limit. The instances with a zero optimality gap are M1-0-1 and M1-0-2, and correspond to the two smallest instances with only 372 candidate locations and associated problems TSP and SDVRP, respectively. Concerning computational time, the exact method is much faster in proving optimality for TSP instances ($\text{M}1-$, $\text{M}2-$, $\text{M}3-$ and $\text{M}4-0-1$) than the heuristic method in finding a solution with no improvement. For the other problem types (SDVRP and CTP), we observe that the difference in computational time is not so pronounced.

\begin{table}[t]
	\small
	\centering
	\caption{Computational results of the exact and heuristic methods for the real-life instances solved to optimality by RN-wS.}
	\label{tab:heur_val_real}
	\begin{threeparttable}
		\begin{adjustbox}{max width=\textwidth}
			\begin{tabular}{ c | S[table-format=5.1] S[table-format=5.1] | S[table-format=4.1]  S[table-format=4.0]  | S[table-format=1.1] | c | S[table-format=5.1] S[table-format=5.1]  | S[table-format=4.1] S[table-format=4.0] | S[table-format=1.1] }
				\toprule
				& \multicolumn{2}{c|}{\textbf{Objective value}} & \multicolumn{2}{c|}{\textbf{Comp. time (s)}} & & & \multicolumn{2}{c|}{\textbf{Objective value}} & \multicolumn{2}{c|}{\textbf{Comp. time (s)}} &  \\
				\makecell{Instance (M$x$-$\gamma$-$m$)}  & $\text{\makecell{Exact \\ (RN-wS)}}$ &  $\text{\makecell{Heuristic \\ method}}$ & $\text{\makecell{Exact \\ (RN-wS)}}$  & $\text{\makecell{Heuristic \\ method\tnote{a}}}$ & $\text{\makecell{Optimality \\ gap (\%)}}$ & \makecell{Instance (M$x$-$\gamma$-$m$)}  & $\text{\makecell{Exact \\ (RN-wS)}}$ &  $\text{\makecell{Heuristic \\ method}}$ & $\text{\makecell{Exact \\ (RN-wS)}}$  & $\text{\makecell{Heuristic \\ method\tnote{a}}}$ & $\text{\makecell{Optimality \\ gap (\%)}}$  \\
				\midrule
				$\text{M1}-0-1$ & 20605.9 & 20605.9 & 1.1 & 6004 & 0.0 & $\text{M3}-0-1$ & 35436.6 & 35956.6& 4.4 & 2473 & 1.4 \\ 
				$\text{M1}-0-2$ & 24033.5 &24033.5  & 801.8 & 642 & 0.0 & $\text{M3}-0-2$ & 36653 & 36677 & 860.8 & 7150 & 0.1 \\
				$\text{M1}-50-1$ & 13014.7 & 13509.6 & 1237.8 & 895  & 3.7 & $\text{M3}-50-1$ & 23382.6 & 25378.8 & 775.9 & 1721  & 7.9 \\
				$\text{M1}-100-1$ & 10526.6 & 10615.9 & 6821.6 & 6083  & 0.8 & $\text{M4}-0-1$ & 42423.5 & 42806.5& 3.9 & 8603 & 0.9 \\
				$\text{M2}-0-1$ & 34910.4 & 35677.4 & 9.5 & 6420  & 2.2 & $\text{M4}-0-2$ & 43507.5 & 43544.5 & 134.5 & 2646  & 0.1 \\
				$\text{M2}-0-2$ & 35426.9 & 35661.4 & 4862.2 & 8969  & 0.7 & $\text{M4}-50-1$ & 29986.1 & 31662.7 & 921.3 & 1370  & 5.3 \\
				$\text{M2}-50-1$ & 21768.3 & 23383.5 & 8801.3 & 7636 & 6.9 & $\text{M4}-100-1$ & 27618.4 & 27701.8 & 7921.5 & 226  & 0.3 \\		
				\bottomrule
			\end{tabular}
		\end{adjustbox}
		\begin{tablenotes}\footnotesize
			\item[a] This computational time corresponds to the time when reaching the best known solution after which \\ no improvement was made
		\end{tablenotes}
	\end{threeparttable}
\end{table}

\begin{figure}[t]
	\centering
	\begin{subfigure}[b]{0.23\textwidth}
		\centering
		\begin{tikzpicture}
			\begin{axis}[
				boxplot/draw direction=y,
				ylabel={relative difference (\%)},
				height=6cm,
				ymin=-101, 
				ymax = 21, 
				boxplot={
					box extend=0.3,
				},
				x=1.75cm,
				ytick={-100,-80,-60,-40,-20,0,20},
				xtick={1,2},
				xticklabels={
					{Small \\ \small (28)},
					{Real-life \\ \small (47)}
				},
				x tick label style={
					text width=3.5cm,
					align=center
				},
				every axis plot/.append style={fill, fill opacity=0.25}
				]
				\addplot [mark=*, boxplot]
				table[row sep=\\,y index=0] {
					data\\
					-4.37052	\\
					-4.33455	\\
					-2.95894	\\
					-1.06784	\\
					-1.06195	\\
					-0.98039	\\
					-0.35211	\\
					-0.17778	\\
					0	\\
					0	\\
					0	\\
					0	\\
					0	\\
					0.332005	\\
					0.354108	\\
					0.719942	\\
					0.782123	\\
					0.792864	\\
					1.163587	\\
					1.1883	\\
					1.270903	\\
					1.401869	\\
					1.606086	\\
					2.658662	\\
					3.583735	\\
					4.217356	\\
					5.934343	\\
					6.653426	\\
				};
				\addplot [mark=*, boxplot]
				table[row sep=\\,y index=0] {
					data\\
					-99.24379856	\\
					-98.48878822	\\
					-89.62637284	\\
					-89.03558599	\\
					-86.78132428	\\
					-81.6259214	\\
					-77.85077561	\\
					-76.41392044	\\
					-45.25770385	\\
					-37.9714071	\\
					-27.30778586	\\
					-18.92237571	\\
					-18.85103768	\\
					-17.10313718	\\
					-16.76932912	\\
					-15.9760918	\\
					-14.81454547	\\
					-13.35242572	\\
					-13.23384108	\\
					-12.56246838	\\
					-11.6822715	\\
					-10.47963658	\\
					-7.417888553	\\
					-7.277225277	\\
					-5.863939325	\\
					-5.81596696	\\
					-4.979785751	\\
					-4.965050504	\\
					-2.100496881	\\
					-1.642477249	\\
					-1.50574149	\\
					-1.228410207	\\
					-1.092265706	\\
					-1.012301001	\\
					-0.838287406	\\
					-0.727907558	\\
					-0.561769351	\\
					-0.418966747	\\
					-0.124136754	\\
					-0.101214935	\\
					3.85501E-05	\\
					0.07902616	\\
					0.429088422	\\
					0.961798124	\\
					1.330726474	\\
					1.551689064	\\
					3.742652017	\\
				};
			\end{axis}
		\end{tikzpicture}
		\caption{Overview}
		\label{fig:perf_heur_over}
	\end{subfigure}
	\begin{subfigure}[b]{0.35\textwidth}
		\centering
		\begin{tikzpicture}
			\begin{axis}[
				boxplot/draw direction=y,
				height=6cm,
				ymin=-5,
				ymax=9,
				boxplot={
					box extend=0.3,
				},
				x=1.75cm,
				xtick={1,2,3},
				ytick={-4,-2,0,2,4,6,8},
				xticklabels={
					{SDVRP\\ \small (3)},
					{CTP\\ \small (7)},
					{m-CTP\\ \small (18)},
				},
				x tick label style={
					text width=3.5cm,
					align=center
				},
				every axis plot/.append style={fill, fill opacity=0.25}
				]
				
				\addplot [mark=*, boxplot]
				table[row sep=\\,y index=0] {
					data\\
					-0.352112676	\\
					0	\\
					0	\\
				};
				\addplot [mark=*, boxplot]
				table[row sep=\\,y index=0] {
					data\\
					-2.958937198	\\
					0.354107649	\\
					0.719942405	\\
					0.782122905	\\
					1.27090301	\\
					3.583735355	\\
					5.934343434	\\
				};
				
				\addplot [mark=*, boxplot]
				table[row sep=\\,y index=0] {
					data\\
					-4.370515329	\\
					-4.334554335	\\
					-1.067839196	\\
					-1.061946903	\\
					-0.980392157	\\
					-0.177777778	\\
					0	\\
					0	\\
					0	\\
					0.332005312	\\
					0.792864222	\\
					1.163586585	\\
					1.188299817	\\
					1.401869159	\\
					1.606086221	\\
					2.658662093	\\
					4.217356042	\\
					6.653426018	\\
				};
			\end{axis}
		\end{tikzpicture}
		\caption{Small instances}
		\label{fig:perf_heur_small_new}
	\end{subfigure}
	\begin{subfigure}[b]{0.35\textwidth}
		\centering
		\begin{tikzpicture}
			\begin{axis}[
				boxplot/draw direction=y,
				height=6cm,
				ymin=-101, 
				ymax=21,
				boxplot={
					box extend=0.3,
				},
				x=1.75cm,
				ytick={-100,-80,-60,-40,-20,0,20},
				xtick={1,2,3},
				xticklabels={
					{SDVRP\\ \small (8)},
					{CTP\\ \small (6)},
					{m-CTP\\ \small (33)},
				},
				x tick label style={
					text width=3.5cm,
					align=center
				},
				every axis plot/.append style={fill, fill opacity=0.25}
				]
				
				\addplot [mark=*, boxplot]
				table[row sep=\\,y index=0] {
					data\\
					-2.100496881	\\
					-1.092265706	\\
					-0.727907558	\\
					-0.561769351	\\
					-0.418966747	\\
					-0.124136754	\\
					-0.101214935	\\
					3.85501E-05	\\
				};
				\addplot [mark=*, boxplot]
				table[row sep=\\,y index=0] {
					data\\
					-16.76932912	\\
					-1.228410207	\\
					-1.012301001	\\
					-0.838287406	\\
					0.429088422	\\
					1.551689064	\\
				};
				
				\addplot [mark=*, boxplot]
				table[row sep=\\,y index=0] {
					data\\
					-99.24379856	\\
					-98.48878822	\\
					-89.62637284	\\
					-89.03558599	\\
					-86.78132428	\\
					-81.6259214	\\
					-77.85077561	\\
					-76.41392044	\\
					-45.25770385	\\
					-37.9714071	\\
					-27.30778586	\\
					-18.92237571	\\
					-18.85103768	\\
					-17.10313718	\\
					-15.9760918	\\
					-14.81454547	\\
					-13.35242572	\\
					-13.23384108	\\
					-12.56246838	\\
					-11.6822715	\\
					-10.47963658	\\
					-7.417888553	\\
					-7.277225277	\\
					-5.863939325	\\
					-5.81596696	\\
					-4.979785751	\\
					-4.965050504	\\
					-1.642477249	\\
					-1.50574149	\\
					0.07902616	\\
					0.961798124	\\
					1.330726474	\\
					3.742652017	\\
				};
			\end{axis}
		\end{tikzpicture}
		\caption{Real-life instances}
		\label{fig:perf_heur_real_new}
	\end{subfigure}
	\caption{Overview of the relative differences (for instances solved to feasibility by RN-wS, excluding the ones solved to optimality) for small and real-life instances (a) and relative differences per problem type for small (b) and real-life instances (c). The number of considered instances for each box plot is included in brackets.}
	\label{fig:perf_heur_new}
\end{figure}
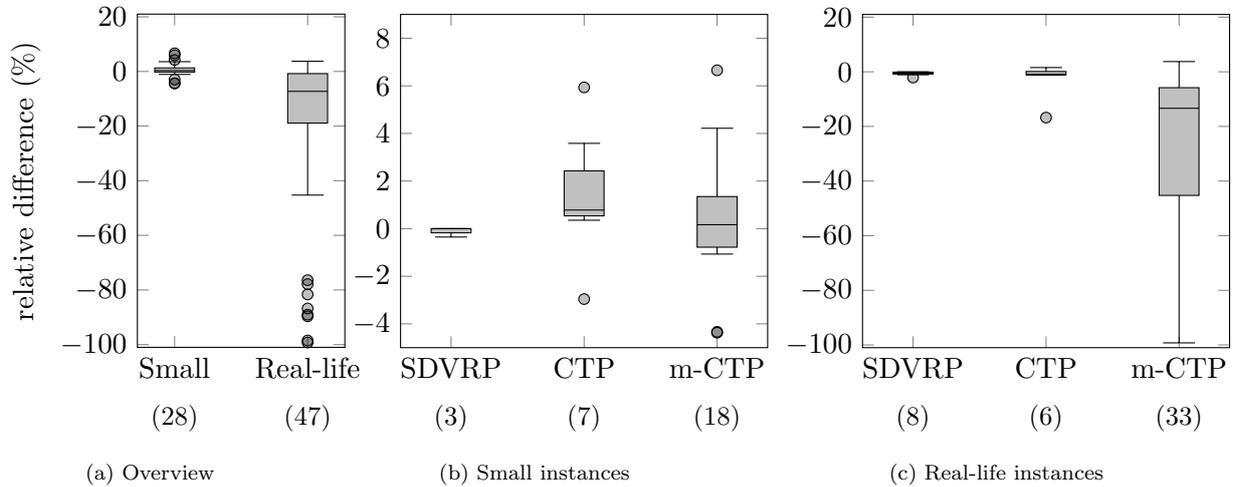

Lastly, \Cref{fig:perf_heur_new} shows the relative differences with respect to the instances solved to feasibility by RN-wS. Note that in this case both the exact and heuristic approaches were run for the same amount of time. The relative difference is calculated as $(Z_{\text{heur}}-Z_{\text{RN-wS}})/Z_{\text{RN-wS}}$, where $Z_{\text{heur}}$ and $Z_{\text{RN-wS}}$ refer to the objective function value obtained with the heuristic and the exact method, respectively. Hence, a negative relative difference indicates that the heuristic method was able to find a better solution than the exact one. As we can see in \Cref{fig:perf_heur_over}, for 28 out of 84 instances the relative differences range from -4.4\% to 6.7\%, with a median value equal to 0.3\%. This means that the heuristic outperforms the exact method for approximately half of the instances. In the real-life case, 47 out of 64 instances were solved to feasibility by RN-wS. Several outliers with relative differences lower than -70\% are reported. For the remaining instances, the relative differences go from -45.3\% to 3.7\%, with 75\% of the instances having a relative difference below -0.8\%. Thus, the heuristic provides a better solution for most of the instances. As we can see in \Cref{fig:perf_heur_small_new} for the small instances, the heuristic method could find a better or equal quality solution for the SDVRP instances, whereas for the CTP and $m$-CTP the median of the relative difference is equal to 0.8\% and 0.3\%, respectively, which indicates that RN-wS could find in general slightly better solutions than the heuristic. For the real-life instances, we can see in \Cref{fig:perf_heur_real_new} that the same behavior for SDVRP instances is observed. For CTP instances, the heuristic method is able to find better solutions for more than 50\% of the instances (median value equal to -0.9\%). The major improvement is observed for $m$-CTP, since for only 4 instances (out of 33) the relative difference is positive (between 0.1\% and 3.7\%). Indeed, 75\% of the instances report a relative difference below -5.8\% and the median value is equal to -13.4\%.

\subsection{Practical aspects} \label{sec:practical_aspects}

In this experiment, we analyze the savings in the total collection time with respect to door-to-door collection (no walking distance) for the real-life instances. As presented in \Cref{sec:perf_heur}, the heuristic performed well in finding good solutions for real-life instances and is therefore considered as the method to use in practice. To this end, we consider the solution generated with the heuristic method (within the three-hour time limit) for each instance with $\gamma > 0$. We then compare it against the heuristic solution of the instance associated with the same municipality and number of tours but $\gamma=0$, i.e. door-to-door collection.

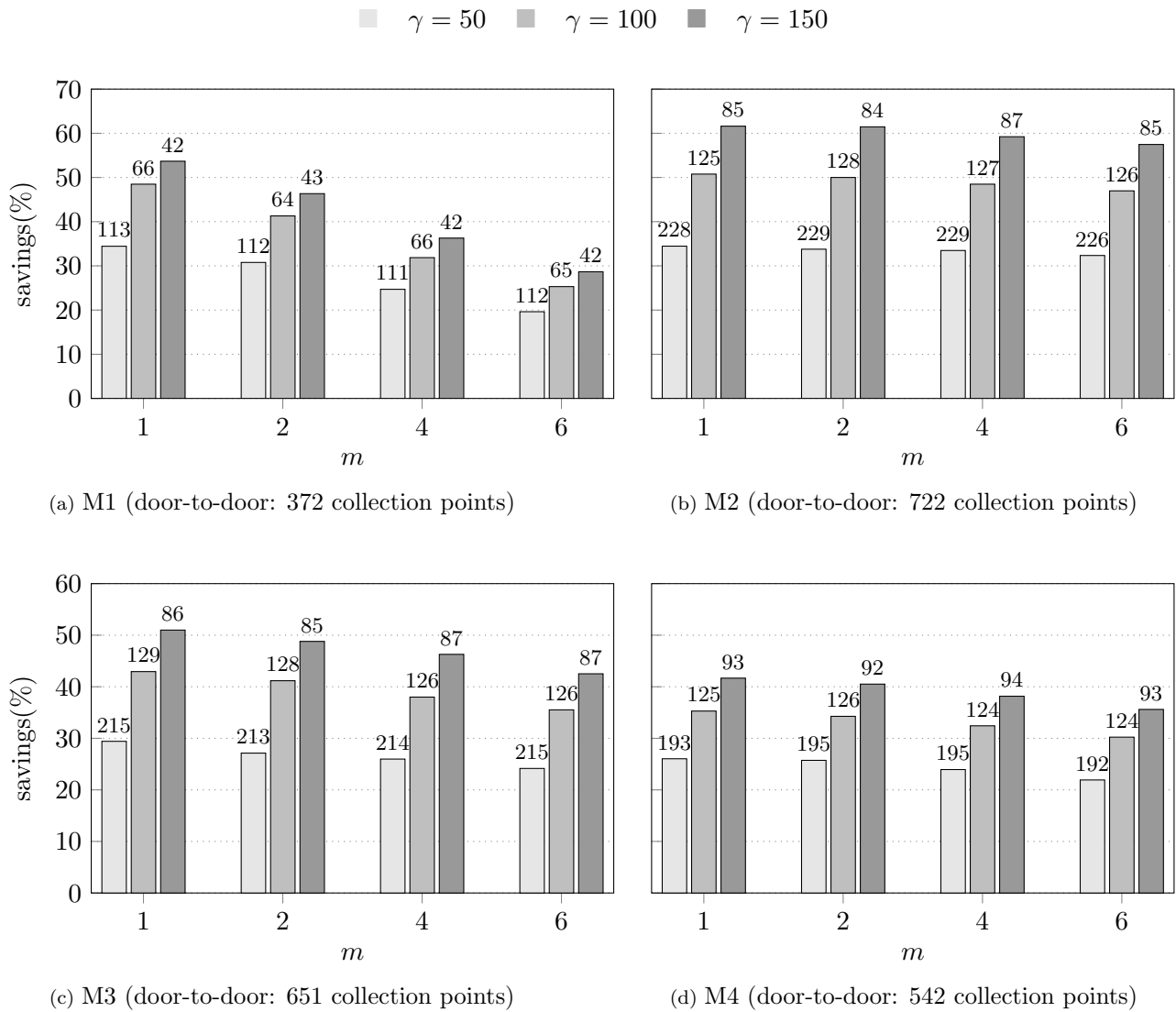
\begin{figure*}[t]
	\begin{center}
		\begin{tabular}{ll ll ll ll}
			\textcolor{tsp}{$\blacksquare$} & $\gamma=50$ & \textcolor{sdvrp}{$\blacksquare$} & $\gamma=100$ & \textcolor{ctp}{$\blacksquare$} & $\gamma=150$ \\
		\end{tabular}  
	\end{center}
	\centering
	\begin{subfigure}[b]{0.475\textwidth}
		\begin{tikzpicture}
			\begin{axis}[height=0.75\textwidth,
				ybar,
				bar width=10pt,
				xlabel={$m$},
				ylabel={savings(\%)},
				xtick=data,
				ymin=0, ymax=70,
				xticklabels from table={\tableA}{m},
				grid style={dotted,black!50},
				xtick pos=left,
				ytick pos=left,
				ytick={0, 10, 20, 30, 40, 50, 60, 70},
				ymajorgrids,
				x=2cm,
				enlarge x limits={abs=0.75cm},
				every node near coord/.append style={font=\footnotesize}
				]
				\addplot[fill=tsp, 
				visualization depends on={value \thisrow{laby50}\as\mylabel},
				nodes near coords = {\mylabel{}},
				] table [col sep = comma, x expr=\coordindex, y=y50]{Plots/M1.csv};
				\addplot[fill=sdvrp, 
				visualization depends on={value \thisrow{laby100}\as\mylabel},
				nodes near coords = {\mylabel{}},
				] table [col sep = comma, x expr=\coordindex, y=y100]{Plots/M1.csv};
				\addplot[fill=ctp, 
				visualization depends on={value \thisrow{laby150}\as\mylabel},
				nodes near coords = {\mylabel{}},
				] table [col sep = comma, x expr=\coordindex, y=y150]{Plots/M1.csv};
			\end{axis}
		\end{tikzpicture}
		\caption[]{\small M1 (door-to-door: 372 collection points)}    
		\label{fig:per_sav_M1}
	\end{subfigure}
	\hfill
	\begin{subfigure}[b]{0.475\textwidth}  
		\centering 
		\begin{tikzpicture}
			\begin{axis}[height=0.75\textwidth,
				ybar,
				bar width=10pt,
				xlabel={$m$},
				xtick=data,
				ymin=0, ymax=70,
				xticklabels from table={\tableB}{m},
				grid style={dotted,black!50},
				xtick pos=left,
				ytick pos=left,
				ytick={0, 10, 20, 30, 40, 50, 60, 70},
				yticklabels=,
				ymajorgrids,
				x=2cm,
				enlarge x limits={abs=0.75cm},
				every node near coord/.append style={font=\footnotesize}
				]
				\addplot[fill=tsp, 
				visualization depends on={value \thisrow{laby50}\as\mylabel},
				nodes near coords = {\mylabel{}},
				] table [col sep = comma, x expr=\coordindex, y=y50]{Plots/M2.csv};
				\addplot[fill=sdvrp, 
				visualization depends on={value \thisrow{laby100}\as\mylabel},
				nodes near coords = {\mylabel{}},
				] table [col sep = comma, x expr=\coordindex, y=y100]{Plots/M2.csv};
				\addplot[fill=ctp, 
				visualization depends on={value \thisrow{laby150}\as\mylabel},
				nodes near coords = {\mylabel{}},
				] table [col sep = comma, x expr=\coordindex, y=y150]{Plots/M2.csv};	
			\end{axis}
		\end{tikzpicture}
		\caption[]{\small M2 (door-to-door: 722 collection points)}      
		\label{fig:per_sav_M2}
	\end{subfigure}
	\vskip\baselineskip
	\begin{subfigure}[b]{0.475\textwidth}   
		\centering 
		\begin{tikzpicture}
			\begin{axis}[height=0.75\textwidth,
				ybar,
				bar width=10pt,
				xlabel={$m$},
				ylabel={savings(\%)},
				xtick=data,
				ymin=0, ymax=60,
				xticklabels from table={\tableA}{m},
				ymajorgrids,
				grid style={dotted,black!50},
				xtick pos=left,
				ytick pos=left,
				x=2cm,
				ytick={0, 10, 20, 30, 40, 50, 60},
				enlarge x limits={abs=0.75cm},
				every node near coord/.append style={font=\footnotesize}
				]
				\addplot[fill=tsp, 
				visualization depends on={value \thisrow{laby50}\as\mylabel},
				nodes near coords = {\mylabel{}},
				] table [col sep = comma, x expr=\coordindex, y=y50]{Plots/M3.csv};
				\addplot[fill=sdvrp, 
				visualization depends on={value \thisrow{laby100}\as\mylabel},
				nodes near coords = {\mylabel{}},
				] table [col sep = comma, x expr=\coordindex, y=y100]{Plots/M3.csv};
				\addplot[fill=ctp, 
				visualization depends on={value \thisrow{laby150}\as\mylabel},
				nodes near coords = {\mylabel{}},
				] table [col sep = comma, x expr=\coordindex, y=y150]{Plots/M3.csv};	
			\end{axis}
		\end{tikzpicture}
		\caption[]{\small M3 (door-to-door: 651 collection points)}  
		\label{fig:per_sav_M3}
	\end{subfigure}
	\hfill
	\begin{subfigure}[b]{0.475\textwidth}   
		\centering 
		\begin{tikzpicture}
			\begin{axis}[height=0.75\textwidth,
				ybar,
				bar width=10pt,
				xlabel={$m$},
				xtick=data,
				ymin=0, ymax=60,
				xticklabels from table={\tableA}{m},
				ymajorgrids,
				grid style={dotted,black!50},
				xtick pos=left,
				ytick pos=left,
				x=2cm,
				ytick={0, 10, 20, 30, 40, 50, 60},
				yticklabels=,
				enlarge x limits={abs=0.75cm},
				every node near coord/.append style={font=\footnotesize}
				]
				\addplot[fill=tsp, 
				visualization depends on={value \thisrow{laby50}\as\mylabel},
				nodes near coords = {\mylabel{}},
				] table [col sep = comma, x expr=\coordindex, y=y50]{Plots/M4.csv};
				\addplot[fill=sdvrp, 
				visualization depends on={value \thisrow{laby100}\as\mylabel},
				nodes near coords = {\mylabel{}},
				] table [col sep = comma, x expr=\coordindex, y=y100]{Plots/M4.csv};
				\addplot[fill=ctp, 
				visualization depends on={value \thisrow{laby150}\as\mylabel},
				nodes near coords = {\mylabel{}},
				] table [col sep = comma, x expr=\coordindex, y=y150]{Plots/M4.csv};	
			\end{axis}
		\end{tikzpicture}
		\caption[]{\small M4 (door-to-door: 542 collection points)}  
		\label{fig:per_sav_M4}
	\end{subfigure}
	\caption[] {\small Percentage of savings with respect to door-to-door collection ($\gamma = 0$) solutions of the real-life instances obtained with the heuristic method.} 
	\label{fig:percentage_savings2}
\end{figure*}

\Cref{fig:percentage_savings2} displays the savings (as a percentage) with respect to door-to-door collection for each municipality, walking distance ($\gamma > 0$) and number of tours ($m$). In all cases we observe savings of at least $20\%$ up to a maximum of over $60\%$ (M2-150-1 and M2-150-2). As expected, the larger the walking distance, the higher the savings. We observe that increasing $\gamma$ to $50$ [meters] comes with a larger gain in collection time (around $30\%$) than increasing $\gamma$ further to $100$ [meters] or $150$ [meters] (additional gain with respect to the previous value of $\gamma$ around $20\%$ and $10\%$, respectively). The same can be observed with the number of collection points visited by the vehicles (labels on top of the bars). Going from door-to-door to $\gamma = 50$ [meters] reduces the number of points by around $60\%$ while increasing to $\gamma = 100$ [meters] or $150$ [meters] removes only half or a third of the points, respectively.

We also observe that the number of tours has almost no impact on the savings for M$2$-M$4$ instances. For M$1$ however, \Cref{fig:per_sav_M1} shows that the savings decrease with increasing value of $m$. This might be due to a larger geographical distance of the network to the depot $\depot{}$ and consequently a higher constant driving cost to dump the waste in each tour.

Taking a closer look at the characteristics of the municipalities (see \Cref{tab:datasets}), we note that for smaller and denser municipalities (\eg M1 and M2) the savings with respect to door-to-door collection are higher than for larger and less densely populated ones (\eg M3 and M4). The highest savings could be achieved for M2, which is the municipality with the highest population.

To visualize the effect of different walking distances on the number of collection points visited by the vehicle, \Cref{fig:fig_M1_solutions} in \ref{sec:app_pract} presents the network representation of the solutions obtained by the heuristic method for all instances ($\gamma \in \{0,50,100,150\}, m \in \{1,2,4,6\}$) of the dataset M1. The tours are marked in different colors, which helps to easily identify the various clusters of collection points. For visualization reasons, the depot is not shown in the images. We observe that with increasing walking distance (images from bottom to top) the number of visited nodes in the graph decrease which as discussed above correlates with the gain in collection time. Furthermore, the images show that for each column of instances (same value of $m$) the clustering of nodes is similar based on the areas in which the nodes lie. In the second column ($m=2$), for instance, the two clusters of all instances split the street network into a south-west and a north-east section which in real life are separated by a water canal. Therefore, it makes sense to define separate clusters for these two areas, which is as well visible in the solutions with $m >2$.

\section{Conclusion}
\label{sec:conclusion}

In this paper, we formulate and solve the \probAbbrVRP{}, a particular version of the $m$-CTP where the constraints on the length and number of nodes of each tour are replaced by vehicle capacity constraints. Furthermore, the coverage of demand nodes is determined by exogenously given preference lists that enforce that each demand node is covered at the first node in its preference list visited by a vehicle. We develop two compact MILP formulations for each considered graph representation: a road-network (RN) graph and a customer-based (CG) graph as typically used in VRP. In both cases, we adapt the formulation to allow for split collection, \ie the demand of the nodes that must be covered can be split and served in different vehicles. To solve practically relevant instances, we propose a two-phased heuristic that first generates sets of nodes to be visited (set covers) and then determines the tours. In the first phase, we obtain a collection of set covers by addressing multiple SCP on the nodes to visit and their alternatives, i.e. additional nodes that can cover the same demand nodes. In the second phase, the SDVRP on each considered set cover is reduced to a CVRP with an a priori splitting strategy so that it can be solved with the state-of-the-art metaheuristic for CVRP (\citealp{vidal2020hybrid}).

The MILP formulations and the heuristic method are tested on a set of small instances inspired by \citet{letchford2013compact} and a set of real-life instances from four municipalities in Switzerland. The computational experiments confirm the advantages of the road-network representation such that the associated formulation (RN) outperforms its counterpart (CG) and provides a more intuitive characterization of the actual network. In particular, RN performs well for instances with a small number of tours ($m$) as seen in Tables \ref{tab:exhaustive_results_3}-\ref{tab:exhaustive_results_4} in \ref{sec:app_results} for the real-life instances. We know that split collection is the state of practice and therefore conclude to use the RN-wS formulation for further analysis. We further observe that for practically relevant problem instances ($m$-CTP instances with $\gamma>0,m>1$), RN-wS is in general not able to prove optimality anymore, and might only provide a feasible solution or no solution at all.

The proposed two-phased heuristic method provides good solutions for the \probAbbrVRP{} and can handle the large instances that general MILP solvers fail to solve. This method was able to find optimal or close to optimal solutions with optimality gaps below $0.5\%$ and $3.5\%$ for $75\%$ of the small and real-life instances solved to optimality by RN-wS, respectively. In addition, it found better solutions for many of the real-life instances for which the exact method failed at proving optimality within the given time limit. From a practical perspective, we observe that the larger the walking distance $\gamma$, the higher the savings in travel time with respect to door-to-door collection ($\gamma=0$). Small and dense municipalities report higher savings than large and less densely populated ones.

The presented MILP formulation and heuristic method could be further extended to accommodate other waste collection systems. For instance, an interesting concept results from introducing intermediate disposal facilities (\citealp{ramos2020new, markov2016integrating}) and replacing the homogeneous fleet of vehicles by an heterogeneous one consisting of small vehicles, potentially electric, that bring the waste to the intermediate facilities and large ones that empty them and bring the waste to the depot. Additionally, the link between the two phases of the heuristic could be reinforced via a feedback loop that, for example, increases the probability of selecting collection points which are present in good solutions. To better represent reality, uncertain waste productions and travel times could be modeled by a set of discrete scenarios with the goal to minimize the worst or average cost over all of them.

\section*{Acknowledgements}
\noindent
The authors gratefully acknowledge the support of Innosuisse under grant 36157.1 IP-EE. They also appreciate the involvement of Schwendimann AG in conducting practical experiments and providing the associated data.


\bibliographystyle{elsarticle-harv} 
\bibliography{ref}
\newpage


\appendix
\section{}

\subsection{Exhaustive results for compact MILP formulations and heuristic}
\label{sec:app_results}

\renewcommand{\thefootnote}{\fnsymbol{footnote}}
\renewcommand*{\thefootnote}{\fnsymbol{footnote}}

\begin{table}[h]
	\small
	\centering
	\caption{Exhaustive results for the compact MILP formulations and the heuristic.}
	\label{tab:exhaustive_results_1}
	
	\begin{adjustbox}{max width=\textwidth}
		\begin{threeparttable}
			
			\begin{tabular}{ c | l | c  | c c | c c | c c | c c  | c c | c }
				
				\toprule
				& & & \multicolumn{2}{c|}{\textbf{CG-noS}} &  \multicolumn{2}{c|}{\textbf{RN-noS}} & \multicolumn{2}{c|}{\textbf{CG-wS}} &  \multicolumn{2}{c|}{\textbf{RN-wS}} & \multicolumn{2}{c|}{\textbf{Heuristic}}  & \\
				& \makecell{Instance} & \makecell{Candidate \\ locations} & $\text{\makecell{Obj.\\ fun.}}$ & $\text{\makecell{Comp. \\ time (s)}}$ & $\text{\makecell{Obj. \\ fun.}}$ & $\text{\makecell{Comp. \\ time (s)}}$ & $\text{\makecell{Obj. \\ fun.}}$ & $\text{\makecell{Comp. \\ time(s)}}$ & $\text{\makecell{Obj. \\ fun.}}$ & $\text{\makecell{Comp. \\ time (s)}}$ & $\text{\makecell{Obj. \\ fun.}}$ & $\text{\makecell{Comp. \\ time (s)}}$ & \makecell{Gap \% \\ to RN-wS}\tnote{1} \\
				\midrule
				\multirow{6}{*}{\STAB{\rotatebox[origin=c]{90}{Small instances}}}
				& STSP\_50\_1/3-0-1 & 15 & \textbf{789}\footnotemark[1] & 0.2 & \textbf{789}\footnotemark[1] & 1.4 & \textbf{789}\footnotemark[1] & 0.4 & \textbf{789}\footnotemark[1] & 1.1 & \textbf{789} & 2 & 0\\
				& STSP\_50\_1/3-0-2 & 15 & \textbf{1164}\footnotemark[1] & 0.6 & \textbf{1164}\footnotemark[1] & 7.5 & \textbf{1164}\footnotemark[1] & 7.9 & \textbf{1164}\footnotemark[1] & 17.1 & \textbf{1164} & 2 & 0\\
				& STSP\_50\_1/3-30-1 & 37 & \textbf{764}\footnotemark[1] & 8.7 & \textbf{764}\footnotemark[1] & 3.3 & \textbf{764}\footnotemark[1] & 8.0 & \textbf{764}\footnotemark[1] & 2.7 & \textbf{764} & 3 & 0\\
				& STSP\_50\_1/3-30-2 & 37 & \textbf{1024}\footnotemark[1] & 105.6 & \textbf{1024}\footnotemark[1] & 121.0 & \textbf{1024}\footnotemark[1] & 1014.3 & \textbf{1024}\footnotemark[1] & 127.8 & \textbf{1024} & 4 & 0\\
				& STSP\_50\_2/3-0-1 & 32 & \textbf{978}\footnotemark[1] & 1.4 & \textbf{978}\footnotemark[1] & 1.8 & \textbf{978}\footnotemark[1] & 2.7 & \textbf{978}\footnotemark[1] & 1.7 & \textbf{978} & 2 & 0\\
				& STSP\_50\_2/3-0-2 & 32 & \textbf{1362}\footnotemark[1] & 12.4 & \textbf{1362}\footnotemark[1] & 32.4 & \textbf{1362}\footnotemark[1] & 606.4 & \textbf{1362}\footnotemark[1] & 73.0 & \textbf{1362} & 3 & 0\\
				& STSP\_50\_2/3-30-1 & 49 & \textbf{839}\footnotemark[1] & 38.8 & \textbf{839}\footnotemark[1] & 8.3 & \textbf{839}\footnotemark[1] & 172.8 & \textbf{839}\footnotemark[1] & 4.6 & 896 & 249 & 6.8\\
				& STSP\_50\_2/3-30-2 & 49 & \textbf{1097} & 10800 & \textbf{1097}\footnotemark[1] & 3532.3 & 1099 & 10800 & \textbf{1097}\footnotemark[1] & 2483.1 & 1142 & 2 & 4.1\\
				& STSP\_50\_3/3-0-1 & 49 & \textbf{1107}\footnotemark[1] & 4.9 & \textbf{1107}\footnotemark[1] & 3.5 & \textbf{1107}\footnotemark[1] & 12.1 & \textbf{1107}\footnotemark[1] & 3.2 & \textbf{1107} & 3 & 0\\
				& STSP\_50\_3/3-0-2 & 49 & \textbf{1446}\footnotemark[1] & 5.9 & \textbf{1446}\footnotemark[1] & 16.8 & \textbf{1446}\footnotemark[1] & 177.0 & \textbf{1446}\footnotemark[1] & 18.9 & \textbf{1446} & 4 & 0\\
				& STSP\_50\_3/3-30-1 & 49 & \textbf{849}\footnotemark[1] & 261.4 & \textbf{849}\footnotemark[1] & 6.7 & \textbf{849}\footnotemark[1] & 198.7 & \textbf{849}\footnotemark[1] & 6.8 & 902 & 9 & 6.2\\
				& STSP\_50\_3/3-30-2 & 49 & \textbf{1107} & 10800 & \textbf{1107}\footnotemark[1] & 3365.9 & 1109 & 10800 & \textbf{1107}\footnotemark[1] & 1953.0 & 1123 & 4 & 1.4\\
				& STSP\_75\_1/3-0-1 & 24 & \textbf{830}\footnotemark[1] & 1.3 & \textbf{830}\footnotemark[1] & 1.9 & \textbf{830}\footnotemark[1] & 1.4 & \textbf{830}\footnotemark[1] & 1.9 & \textbf{830} & 2 & 0\\
				& STSP\_75\_1/3-0-2 & 24 & \textbf{1186}\footnotemark[1] & 1.4 & \textbf{1186}\footnotemark[1] & 52.2 & \textbf{1186}\footnotemark[1] & 36.1 & \textbf{1186}\footnotemark[1] & 104.0 & \textbf{1186} & 3 & 0\\
				& STSP\_75\_1/3-30-1 & 62 & \textbf{663}\footnotemark[1] & 1055.4 & \textbf{663}\footnotemark[1] & 136.8 & \textbf{663}\footnotemark[1] & 1044.5 & \textbf{663}\footnotemark[1] & 108.2 & 666 & 5 & 0.5\\
				& STSP\_75\_1/3-30-2 & 62 & \textbf{928} & 10800 & \textbf{928} & 10800 & \textbf{928} & 10800 & \textbf{928} & 10800 & \textbf{928} & 28 & 0\\
				& STSP\_75\_2/3-0-1 & 49 & \textbf{1029}\footnotemark[1] & 7.9 & \textbf{978}\footnotemark[1] & 1.8 & \textbf{978}\footnotemark[1] & 2.7 & \textbf{978}\footnotemark[1] & 1.7 & \textbf{978} & 2 & 0\\
				& STSP\_75\_2/3-0-2 & 32 & \textbf{1362}\footnotemark[1] & 12.4 & \textbf{1362}\footnotemark[1] & 32.4 & \textbf{1362}\footnotemark[1] & 606.4 & \textbf{1362}\footnotemark[1] & 73.0 & \textbf{1362} & 3 & 0\\
				& STSP\_75\_2/3-30-1 & 49 & \textbf{839}\footnotemark[1] & 38.8 & \textbf{839}\footnotemark[1] & 8.3 & \textbf{839}\footnotemark[1] & 172.8 & \textbf{839}\footnotemark[1] & 4.6 & 896 & 249 & 6.8\\
				& STSP\_75\_2/3-30-2 & 49 & \textbf{1097} & 10800 & \textbf{1097}\footnotemark[1] & 3532.3 & 1099 & 10800 & \textbf{1097}\footnotemark[1] & 2483.1 & 1142 & 2 & 4.1\\
				& STSP\_75\_3/3-0-1 & 49 & \textbf{1107}\footnotemark[1] & 4.9 & \textbf{1107}\footnotemark[1] & 3.5 & \textbf{1107}\footnotemark[1] & 12.1 & \textbf{1107}\footnotemark[1] & 3.2 & \textbf{1107} & 3 & 0\\
				& STSP\_75\_3/3-0-2 & 49 & \textbf{1446}\footnotemark[1] & 5.9 & \textbf{1446}\footnotemark[1] & 16.8 & \textbf{1446}\footnotemark[1] & 177.0 & \textbf{1446}\footnotemark[1] & 18.9 & \textbf{1446} & 4 & 0\\
				& STSP\_75\_3/3-30-1 & 49 & \textbf{849}\footnotemark[1] & 261.4 & \textbf{849}\footnotemark[1] & 6.7 & \textbf{849}\footnotemark[1] & 198.7 & \textbf{849}\footnotemark[1] & 6.8 & 902 & 9 & 6.2\\
				& STSP\_75\_3/3-30-2 & 49 & \textbf{1107} & 10800 & \textbf{1107}\footnotemark[1] & 3365.9 & 1109 & 10800 & \textbf{1107}\footnotemark[1] & 1953.0 & 1123 & 4 & 1.4\\
				& STSP\_100\_1/3-0-1& 32 & \textbf{919}\footnotemark[1] & 3.9 & \textbf{919}\footnotemark[1] & 7.3 & \textbf{919}\footnotemark[1] & 5.8 & \textbf{919}\footnotemark[1] & 6.9 & \textbf{919} & 2 & 0\\
				& STSP\_100\_1/3-0-2 & 32 & \textbf{1339}\footnotemark[1] & 5.2 & \textbf{1339}\footnotemark[1] & 52.2 & \textbf{1339}\footnotemark[1] & 193.9 & \textbf{1339}\footnotemark[1] & 224.2 & \textbf{1339} & 3 & 0\\
				& STSP\_100\_1/3-30-1 & 83 & \textbf{687} & 10800 & \textbf{687}\footnotemark[1] & 490.5 & \textbf{687} & 10800 & \textbf{687}\footnotemark[1] & 133.1 & 755 & 6 & 9.9\\
				& STSP\_100\_1/3-30-2 & 83 & 1010 & 10800 & 1010 & 10800 & 1016 & 10800 & \textbf{1007} & 10800 & 1074 & 3 & 6.7\\
				& STSP\_100\_2/3-0-1 & 65 & \textbf{1193}\footnotemark[1] & 10.8 & \textbf{1193}\footnotemark[1] & 6.5 & \textbf{1193}\footnotemark[1] & 9.7 & \textbf{1193}\footnotemark[1] & 7.4 & \textbf{1193} & 1 & 0\\
				& STSP\_100\_2/3-0-2 & 65 & \textbf{1668}\footnotemark[1] & 54.9 & \textbf{1668}\footnotemark[1] & 312.8 & \textbf{1668}\footnotemark[1] & 744.4 & \textbf{1668}\footnotemark[1] & 428.8 & \textbf{1668} & 1 & 0\\
				& STSP\_100\_2/3-30-1 & 97 & 806 & 10800 & \textbf{792}\footnotemark[1] & 2472.5 & \textbf{792} & 10800 & \textbf{792} & 10800 & 839 & 351 & 5.9\\
				& STSP\_100\_2/3-30-2 & 97 & 1037 & 10800 & \textbf{1035} & 10800 & 1037 & 10800 & \textbf{1035} & 10800 & \textbf{1035} & 1298 & 0\\
				& STSP\_100\_3/3-0-1 & 99 & \textbf{1445}\footnotemark[1] & 41.2 & \textbf{1445}\footnotemark[1] & 7.7 & \textbf{1445}\footnotemark[1] & 50.5 & \textbf{1445}\footnotemark[1] & 10.0 & \textbf{1445} & 5 & 0\\
				& STSP\_100\_3/3-0-2 & 99 & \textbf{2038}\footnotemark[1] & 102.3 & \textbf{2038}\footnotemark[1] & 205.3 & 2094 & 10800 & \textbf{2038}\footnotemark[1] & 364.2 & \textbf{2038} & 1 & 0\\
				& STSP\_100\_3/3-30-1 & 99 & 896 & 10800 & \textbf{895} & 10800 & 896 & 10800 & \textbf{895} & 10800 & 902 & 320 & 0.8\\
				& STSP\_100\_3/3-30-2 & 99 & 1127 & 10800 & 1125 & 10800 & 1174 & 10800 & 1125 & 10800 & \textbf{1123} & 2898 & -0.2\\
				& STSP\_125\_1/3-0-1 & 40 & \textbf{1143}\footnotemark[1] & 3.2 & \textbf{1143}\footnotemark[1] & 4.1 & \textbf{1143}\footnotemark[1] & 2.5 & \textbf{1143}\footnotemark[1] & 4.5 & \textbf{1143} & 2 & 0\\
				& STSP\_125\_1/3-0-2 & 40 & \textbf{1478}\footnotemark[1] & 13.8 & \textbf{1478}\footnotemark[1] & 88.0 & \textbf{1478}\footnotemark[1] & 100.9 & \textbf{1478}\footnotemark[1] & 677.4 & \textbf{1478} & 3 & 0\\
				& STSP\_125\_1/3-15-1 & 88 & \textbf{1125}\footnotemark[1] & 881.6 & \textbf{1125}\footnotemark[1] & 417.3 & \textbf{1125}\footnotemark[1] & 1699.1 & \textbf{1125}\footnotemark[1] & 299.5 & 1136 & 9 & 1.0\\
				& STSP\_125\_1/3-15-2 & 88 & \textbf{1265} & 10800 & \textbf{1265} & 10800 & 1301 & 10800 & \textbf{1265} & 10800 & \textbf{1265} & 28 & 0\\
				\bottomrule
				
			\end{tabular}
			
			\begin{tablenotes}\footnotesize
				\item[1] The gap is calculated as the relative difference of the objective function value obtained with the heuristic method with respect to the corresponding value of the RN-wS formulation.
			\end{tablenotes}
			
		\end{threeparttable}
	\end{adjustbox}
	
\end{table}

\begin{table}[H]
	\small
	\centering
	\caption{Exhaustive results for the compact MILP formulations and the heuristic.}
	\label{tab:exhaustive_results_2}
	\begin{adjustbox}{max width=\textwidth}
		\begin{threeparttable}
			
			\begin{tabular}{ c | l | c  | c c | c c | c c | c c  | c c | c }
				\toprule
				& & & \multicolumn{2}{c|}{\textbf{CG-noS}} &  \multicolumn{2}{c|}{\textbf{RN-noS}} & \multicolumn{2}{c|}{\textbf{CG-wS}} &  \multicolumn{2}{c|}{\textbf{RN-wS}} & \multicolumn{2}{c|}{\textbf{Heuristic}}  & \\
				& \makecell{Instance} & \makecell{Candidate \\ locations} & $\text{\makecell{Obj.\\ fun.}}$ & $\text{\makecell{Comp. \\ time (s)}}$ & $\text{\makecell{Obj. \\ fun.}}$ & $\text{\makecell{Comp. \\ time (s)}}$ & $\text{\makecell{Obj. \\ fun.}}$ & $\text{\makecell{Comp. \\ time(s)}}$ & $\text{\makecell{Obj. \\ fun.}}$ & $\text{\makecell{Comp. \\ time (s)}}$ & $\text{\makecell{Obj. \\ fun.}}$ & $\text{\makecell{Comp. \\ time (s)}}$ & \makecell{Gap \% \\ to RN-wS}\tnote{1} \\
				\midrule
				\multirow{6}{*}{\STAB{\rotatebox[origin=c]{90}{Small instances}}}
				& STSP\_125\_2/3-0-1 & 82 & \textbf{1416}\footnotemark[1] & 19.6 & \textbf{1416}\footnotemark[1] & 16 & \textbf{1416}\footnotemark[1] & 22.8 & \textbf{1416}\footnotemark[1] & 20.6 & \textbf{1416} & 1 & 0\\
				& STSP\_125\_2/3-0-2 & 82 & \textbf{1943}\footnotemark[1] & 184.3 & \textbf{1943}\footnotemark[1] & 2062.9 & \textbf{1943} & 10800 & \textbf{1943}\footnotemark[1] & 4741.2 & \textbf{1943} & 6 & 0\\
				& STSP\_125\_2/3-15-1 & 117 & \textbf{1311}\footnotemark[1] & 5497.7 & \textbf{1311}\footnotemark[1] & 1056.6 & \textbf{1311} & 10800 & \textbf{1311}\footnotemark[1] & 626.0 & 1341 & 4331 & 2.3\\
				& STSP\_125\_2/3-15-2 & 117 & 1463 & 10800 & \textbf{1461} & 10800 & 1501 & 10800 & \textbf{1461} & 10800 & 1478 & 199 & 1.2\\
				& STSP\_125\_3/3-0-1 & 124 & \textbf{1582}\footnotemark[1] & 81.7 & \textbf{1582}\footnotemark[1] & 26.7 & \textbf{1582}\footnotemark[1] & 231.9 & \textbf{1582}\footnotemark[1] & 26.6 & \textbf{1582} & 6 & 0\\
				& STSP\_125\_3/3-0-2 & 124 & \textbf{2295}\footnotemark[1] & 2423.5 & \textbf{2295}\footnotemark[1] & 9033.8 & \textbf{2295} & 10800 & \textbf{2295}\footnotemark[1] & 10659.8 & \textbf{2295} & 7 & 0\\
				& STSP\_125\_3/3-15-1 & 124 & \textbf{1377}\footnotemark[1] & 5776.3 & \textbf{1377}\footnotemark[1] & 739.3 & \textbf{1377}\footnotemark[1] & 5335.9 & \textbf{1377}\footnotemark[1] & 739.4 & 1386 & 5006 & 0.7\\
				& STSP\_125\_3/3-15-2 & 124 & \textbf{1497} & 10800 & 1526 & 10800 & 1516 & 10800 & 1498 & 10800 & 1519 & 6151 & 1.4\\
				& STSP\_150\_1/3-0-1 & 49 & \textbf{1105}\footnotemark[1] & 6.9 & \textbf{1105}\footnotemark[1] & 6.6 & \textbf{1105}\footnotemark[1] & 6.7 & \textbf{1105}\footnotemark[1] & 7.2 & \textbf{1105} & 3 & 0\\
				& STSP\_150\_1/3-0-2 & 49 & \textbf{1502}\footnotemark[1] & 147.7 & \textbf{1502}\footnotemark[1] & 945.4 & \textbf{1502}\footnotemark[1] & 945.6 & \textbf{1502}\footnotemark[1] & 3429.5 & \textbf{1502} & 5 & 0\\
				& STSP\_150\_1/3-15-1 & 112 & \textbf{984}\footnotemark[1] & 10621.1 & \textbf{984}\footnotemark[1] & 3291.4 & 987 & 10800 & \textbf{984}\footnotemark[1] & 1505.4 & 1022 & 6 & 3.9\\
				& STSP\_150\_1/3-15-2 & 112 & 1183 & 10800 & \textbf{1166} & 10800 & 1191 & 10800 & \textbf{1166} & 10800 & 1197 & 12 & 2.7\\
				& STSP\_150\_2/3-0-1 & 99 & \textbf{1545}\footnotemark[1] & 79.6 & \textbf{1545}\footnotemark[1] & 42.9 & \textbf{1545}\footnotemark[1] & 139.6 & \textbf{1545}\footnotemark[1] & 42.2 & \textbf{1545} & 5 & 0\\
				& STSP\_150\_2/3-0-2 & 99 & \textbf{2097}\footnotemark[1] & 108.5 & \textbf{2097}\footnotemark[1] & 801.4 & \textbf{2097}\footnotemark[1] & 10002.4 & \textbf{2097}\footnotemark[1] & 1299.8 & \textbf{2097} & 3 & 0\\
				& STSP\_150\_2/3-15-1 & 144 & \textbf{1279}\footnotemark[1] & 9032.4 & \textbf{1279}\footnotemark[1] & 10313.7 & \textbf{1279} & 10800 & \textbf{1279}\footnotemark[1] & 7613.0 & 1309 & 1225 & 2.3\\
				& STSP\_150\_2/3-15-2 & 144 & 1396 & 10800 & \textbf{1393} & 10800 & 1442 & 10800 & 1428 & 10800 & 1414 & 234 & -1.0\\
				& STSP\_150\_3/3-0-1 & 149 & \textbf{1733}\footnotemark[1] & 142.0 & \textbf{1733}\footnotemark[1] & 29.2 & \textbf{1733}\footnotemark[1] & 332.8 & \textbf{1733}\footnotemark[1] & 32.4 & \textbf{1733} & 1 & 0\\
				& STSP\_150\_3/3-0-2 & 149 & \textbf{2552}\footnotemark[1] & 5267.3 & \textbf{2552}\footnotemark[1] & 9342.8 & 2619 & 10800 & \textbf{2552} & 10800 & \textbf{2552} & 8 & 0\\
				& STSP\_150\_3/3-15-1 & 149 & \textbf{1412} & 10800 & \textbf{1412} & 10800 & \textbf{1412} & 10800 & \textbf{1412} & 10800 & 1417 & 7326 & 0.4\\
				& STSP\_150\_3/3-15-2 & 149 & \textbf{1466} & 10800 & 1506 & 10800 & 1575 & 10800 & 1533 & 10800 & \textbf{1466} & 88 & -4.4\\
				& STSP\_175\_1/3-0-1 & 57 & \textbf{1272}\footnotemark[1] & 5.0 & \textbf{1272}\footnotemark[1] & 7.9 & \textbf{1272}\footnotemark[1] & 8.4 & \textbf{1272}\footnotemark[1] & 8.2 & \textbf{1272} & 3 & 0\\
				& STSP\_175\_1/3-0-2 & 57 & \textbf{1674}\footnotemark[1] & 12.1 & \textbf{1674}\footnotemark[1] & 313.2 & \textbf{1674}\footnotemark[1] & 463.2 & \textbf{1674}\footnotemark[1] & 3564.8 & \textbf{1674} & 7 & 0\\
				& STSP\_175\_1/3-15-1 & 140 & \textbf{1103}\footnotemark[1] & 4878.5 & \textbf{1103}\footnotemark[1] & 3279.6 & \textbf{1103}\footnotemark[1] & 8249.0 & \textbf{1103}\footnotemark[1] & 2505.2 & 1173 & 26 & 6.3\\
				& STSP\_175\_1/3-15-2 & 140 & \textbf{1233} & 10800 & \textbf{1233} & 10800 & 1287 & 10800 & \textbf{1233} & 10800 & 1285 & 2270 & 4.2\\
				& STSP\_175\_2/3-0-1 & 115 & \textbf{1645}\footnotemark[1] & 78.9 & \textbf{1645}\footnotemark[1] & 30.4 & \textbf{1645}\footnotemark[1] & 83.3 & \textbf{1645}\footnotemark[1] & 30.6 & \textbf{1645} & 5 & 0\\
				& STSP\_175\_2/3-0-2 & 115 & \textbf{2303}\footnotemark[1] & 297.0 & \textbf{2303}\footnotemark[1] & 1792.1 & 2398 & 10800 & \textbf{2303}\footnotemark[1] & 2506.0 & \textbf{2303} & 11 & 0\\
				& STSP\_175\_2/3-15-1 & 169 & \textbf{1348} & 10800 & \textbf{1348} & 10800 & 1356 & 10800 & 1389 & 10800 & 1399 & 445 & 0.7\\
				& STSP\_175\_2/3-15-2 & 169 & \textbf{1476} & 10800 & 1514 & 10800 & 1488 & 10800 & 1506 & 10800 & 1511 & 4150 & 0.3\\
				& STSP\_175\_3/3-0-1 & 174 & \textbf{1891}\footnotemark[1] & 229.2 & \textbf{1891}\footnotemark[1] & 24.4 & \textbf{1891}\footnotemark[1] & 826.6 & \textbf{1891}\footnotemark[1] & 24.5 & \textbf{1891} & 9 & 0\\
				& STSP\_175\_3/3-0-2 & 174 & \textbf{2830}\footnotemark[1] & 4007.9 & \textbf{2830}\footnotemark[1] & 9684.2 & 3092 & 10800 & 2840 & 10800 & \textbf{2830} & 28 & -0.4\\
				& STSP\_175\_3/3-15-1 & 174 & 1507 & 10800 & \textbf{1488} & 10800 & 1497 & 10800 & 1495 & 10800 & 1514 & 3304 & 1.3\\
				& STSP\_175\_3/3-15-2 & 174 & \textbf{1571} & 10800 & 1631 & 10800 & 1611 & 10800 & 1592 & 10800 & 1575 & 1799 & -1.1\\
				& STSP\_200\_1/3-0-1 & 65 & \textbf{1295}\footnotemark[1] & 3.8 & \textbf{1295}\footnotemark[1] & 9.1 & \textbf{1295}\footnotemark[1] & 7.0 & \textbf{1295}\footnotemark[1] & 8.9 & \textbf{1295} & 1 & 0\\
				& STSP\_200\_1/3-0-2 & 65 & \textbf{1764}\footnotemark[1] & 169.2 & \textbf{1764}\footnotemark[1] & 1729.5 & \textbf{1764} & 10800 & \textbf{1764} & 10800 & \textbf{1764} & 13 & 0\\
				& STSP\_200\_1/3-15-1 & 154 & \textbf{1046} & 10800 & \textbf{1046}\footnotemark[1] & 5873.5 & \textbf{1046} & 10800 & \textbf{1046}\footnotemark[1] & 1077.3 & 1068 & 7742 & 2.1\\
				& STSP\_200\_1/3-15-2 & 154 & \textbf{1183} & 10800 & \textbf{1183} & 10800 & 1190 & 10800 & \textbf{1183} & 10800 & 1202 & 5365 & 1.6\\
				& STSP\_200\_2/3-0-1 & 132 & \textbf{1845}\footnotemark[1] & 89.0 & \textbf{1845}\footnotemark[1] & 26.2 & \textbf{1845}\footnotemark[1] & 153.9 & \textbf{1845}\footnotemark[1] & 26.7 & \textbf{1845} & 6 & 0\\
				& STSP\_200\_2/3-0-2 & 132 & \textbf{2543}\footnotemark[1] & 163.1 & \textbf{2543}\footnotemark[1] & 1599.8 & \textbf{2543}\footnotemark[1] & 10436.0 & \textbf{2543}\footnotemark[1] & 2322.9 & \textbf{2543} & 21 & 0\\
				& STSP\_200\_2/3-15-1 & 198 & \textbf{1443} & 10800 & 1477 & 10800 & 1446 & 10800 & 1451 & 10800 & 1503 & 6322 & 3.6\\
				& STSP\_200\_2/3-15-2 & 198 & \textbf{1512} & 10800 & 1559 & 10800 & 1639 & 10800 & 1638 & 10800 & 1567 & 5631 & -4.3\\
				& STSP\_200\_3/3-0-1 & 199 & \textbf{2018}\footnotemark[1] & 296.5 & \textbf{2018}\footnotemark[1] & 16.5 & \textbf{2018}\footnotemark[1] & 618.5 & \textbf{2018}\footnotemark[1] & 15.9 & \textbf{2018} & 7 & 0\\
				& STSP\_200\_3/3-0-2 & 199 & \textbf{3046}\footnotemark[1] & 370.4 & \textbf{3046} & 2719.8 & 3327 & 10800 & \textbf{3046}\footnotemark[1] & 3934.6 & \textbf{3046} & 16 & 0\\
				& STSP\_200\_3/3-15-1 & 199 & \textbf{1580} & 10800 & 1636 & 10800 & 1610 & 10800 & 1656 & 10800 & 1607 & 2719 & -3.0\\
				& STSP\_200\_3/3-15-2 & 199 & \textbf{1653} & 10800 & 1721 & 10800 & 1726 & 10800 & 1695 & 10800 & 1677 & 5482 & -1.1\\
				\bottomrule
			\end{tabular}
			
			\begin{tablenotes}\footnotesize
				\item[1] The gap is calculated as the relative difference of the objective function value obtained with the heuristic method with respect to the corresponding value of the RN-wS formulation.
			\end{tablenotes}
		\end{threeparttable}
	\end{adjustbox}
	
\end{table}

\begin{table}[H]
	\small
	\centering
	\caption{Exhaustive results for the compact MILP formulations and the heuristic.}
	\label{tab:exhaustive_results_3}
	\begin{adjustbox}{max width=\textwidth}
		\begin{threeparttable}
			
			\begin{tabular}{ c | l | c  | c c | c c | c c | c c  | c c | c }
				\toprule
				& & & \multicolumn{2}{c|}{\textbf{CG-noS}} &  \multicolumn{2}{c|}{\textbf{RN-noS}} & \multicolumn{2}{c|}{\textbf{CG-wS}} &  \multicolumn{2}{c|}{\textbf{RN-wS}} & \multicolumn{2}{c|}{\textbf{Heuristic}}  & \\
				& \makecell{Instance} & \makecell{Candidate \\ locations} & $\text{\makecell{Obj.\\ fun.}}$ & $\text{\makecell{Comp. \\ time (s)}}$ & $\text{\makecell{Obj. \\ fun.}}$ & $\text{\makecell{Comp. \\ time (s)}}$ & $\text{\makecell{Obj. \\ fun.}}$ & $\text{\makecell{Comp. \\ time(s)}}$ & $\text{\makecell{Obj. \\ fun.}}$ & $\text{\makecell{Comp. \\ time (s)}}$ & $\text{\makecell{Obj. \\ fun.}}$ & $\text{\makecell{Comp. \\ time (s)}}$ & \makecell{Gap \% \\ to RN-wS}\tnote{1} \\
				\midrule
				\multirow{6}{*}{\STAB{\rotatebox[origin=c]{90}{Real-life instances}}}
				& M1-0-1 & 372 & \textbf{20605.9}\footnotemark[1] & 5482.5 & \textbf{20605.9}\footnotemark[1] & 1.8 & \textbf{20605.9}\footnotemark[1] & 1333.1 & \textbf{20605.9}\footnotemark[1] & 1.1 & \textbf{20605.9} & 6004 & 0\\
				& M1-0-2 & 372 & \textbf{24033.5}\footnotemark[1] & 1574.8 & \textbf{24033.5}\footnotemark[1] & 13.2 & - & - & \textbf{24033.5}\footnotemark[1] & 801.8 & \textbf{24033.5} & 642 & 0\\
				& M1-0-4 & 372 & \textbf{31128.4}\footnotemark[1] & 840.3 & \textbf{31128.4} & 10800 & - & - & \textbf{31128.4} & 10800 & \textbf{31128.4} & 2841 & 0\\
				& M1-0-6 & 372 & \textbf{38385.7}\footnotemark[1] & 9929.3 & 38433.4 & 10800 & - & - & 38433.4 & 10800 & \textbf{38385.7} & 8231 & -0.1\\
				& M1-50-1 & 592 & - & - & \textbf{13014.7}\footnotemark[1] & 959.2 & 15490.8 & 10800 & \textbf{13014.7}\footnotemark[1]  & 1237.8 & 13509.6 & 895 & 3.8\\
				& M1-50-2 & 592 & 18000.4 & 10800 & \textbf{16421.2} & 10800 & - & - & 16478.5 & 10800 & 16637 & 4050 & 1.0\\
				& M1-50-4 & 592 & 24975.9 & 10800 & 24278.6 & 10800 & - & - & 23799.6 & 10800 & \textbf{23441.2} & 2461 & -1.5\\
				& M1-50-6 & 592 & 170751.9 & 10800 & 34630.5 & 10800 & - & - & 38052.8 & 10800 & \textbf{30852.3} & 7402 & -18.9\\
				& M1-100-1 & 629 & - & - & \textbf{10526.6}\footnotemark[1] & 5475.5 & 29042.1 & 10800 & \textbf{10526.6}\footnotemark[1] & 6821.6 & 10615.9 & 6083 & 0.8\\
				& M1-100-2 & 629 & 16418.8 & 10800 & 14652.7 & 10800 & - & - & 14342.4 & 10800 & \textbf{14106.8} & 5089 & -1.6\\
				& M1-100-4 & 629 & 22068.6 & 10800 & 22195.2 & 10800 & - & - & 22874.8 & 10800 & \textbf{21210.1} & 3828 & -7.3\\
				& M1-100-6 & 629 & 31838.9 & 10800 & 42322 & 10800 & - & - & 30969.2 & 10800 & \textbf{28671.9} & 2262 & -7.4\\
				& M1-150-1 & 636 & 17340.4 & 10800 & 9632 & 10800 & 10442.6 & 10800 & 9663.7 & 10800 & \textbf{9545} & 4038 & -1.2\\
				& M1-150-2 & 636 & 14780.9 & 10800 & 13320.4 & 10800 & - & - & 13692.7 & 10800 & \textbf{12896.3} & 116 & -5.8\\
				& M1-150-4 & 636 & 27393.4 & 10800 & 20782.9 & 10800 & - & - & 22457.8 & 10800 & \textbf{19834.2} & 4580 & -11.7\\
				& M1-150-6 & 636 & 39969.6 & 10800 & 34780.3 & 10800 & - & - & 31604.2 & 10800 & \textbf{27384.3} & 3340 & -13.3\\
				& M2-0-1 & 722 & 100788.6 & 10800 & \textbf{34910.4}\footnotemark[1] & 12.2 & - & - & \textbf{34910.4}\footnotemark[1] & 9.5 & 35677.4 & 6420 & 2.2\\
				& M2-0-2 & 722 & 447901 & 10800 & \textbf{35426.9}\footnotemark[1] & 3581.7 & - & - & \textbf{35426.9}\footnotemark[1]  & 4862.2 & 35661.4 & 8969 & 0.7\\
				& M2-0-4 & 722 & 347966 & 10800 & 36702.4 & 10800 & - & - & 36546 & 10800 & \textbf{36509} & 7774 & -0.1\\
				& M2-0-6 & 722 & - & - & 38080.6 & 10800 & - & - & 38118.9 & 10800 & \textbf{37702.6} & 9399 & -1.1\\
				& M2-50-1 & 1158 & - & - & \textbf{21768.3} & 10800 & - & - & \textbf{21768.3}\footnotemark[1] & 8801.3 & 23383.5 & 7636 & 7.4\\
				& M2-50-2 & 1158 & - & - & 24024.6 & 10800 & - & - & \textbf{23599.8} & 10800 & 23618.4 & 6396 & 0.08\\
				& M2-50-4 & 1158 & - & - & 44976 & 10800 & - & - & - & - & \textbf{24286.6} & 2378 & -\\
				& M2-50-6 & 1158 & - & - &  149764.6 & 10800 & - & - & 232639.6 & 10800 & 25507.6 & 405 & -89\\
				& M2-100-1 & 1262 & - & - & \textbf{17541.9} & 10800 & - & - & 17745.7 & 10800 & 17566.1 & 5934 & -1.0\\
				& M2-100-2 & 1262 & - & - & 21069.1 & 10800 & - & - & 19917.6 & 10800 & \textbf{17830.3} & 1349 & -10.5\\
				& M2-100-4 & 1262 & - & - & 37424.4 & 10800 & - & - & - & - & 18808.4 & 671 & -\\
				& M2-100-6 & 1262 & - & - & 133081.3 & 10800 & - & - & 151242.8 & 10800 & \textbf{19992.3} & 840 & -86.8\\
				& M2-150-1 & 1301 & - & - & 13981.8 & 10800 & - & - & 16455.2 & 10800 & \textbf{13695.8} & 6116 & -16.8\\
				& M2-150-2 & 1301 & - & - & 14832 & 10800 & - & - & 16937.4 & 10800 & \textbf{13744.5} & 4023 & -18.9\\
				& M2-150-4 & 1301 & - & - & 24760.9 & 10800 & - & - & 143524.6 & 10800 & \textbf{14888.7} & 326 & -89.6\\
				& M2-150-6 & 1301 & - & - & 84633.9 & 10800 & - & - & 2120542.1 & 10800 & \textbf{16035.6} & 4478 & -99.2\\
				& M3-0-1 & 651 & 40132.5 & 10800 & \textbf{35436.6}\footnotemark[1] & 5.3 & 517467.6 & 10800 & \textbf{35436.6}\footnotemark[1] & 4.4 & 35956.6 & 2473 & 1.5\\
				& M3-0-2 & 651 & 37789.7 & 10800 & \textbf{36653}\footnotemark[1] & 740.7 & - & - & \textbf{36653}\footnotemark[1] & 860.8 & 36677 & 7150 & 0.07\\
				& M3-0-4 & 651 & 389344.4 & 10800 & 39822.8 & 10800 & - & - & 39838.9 & 10800 & \textbf{39548.9} & 4646 & -0.7\\
				& M3-0-6 & 651 & 520917.8 & 10800 & 42609 & 10800 & - & - & 42757.6 & 10800 & \textbf{42578.4} & 8524 & -0.4\\
				& M3-50-1 & 1083 & - & - & \textbf{23382.6}\footnotemark[1] & 726.5 & - & - & \textbf{23382.6}\footnotemark[1] & 775.9 & 25378.8 & 1721 & 8.5\\
				& M3-50-2 & 1083 & - & - & \textbf{25348.4} & 10800 & - & - & 26376.6 & 10800 & 26727.6 & 4658 & 1.3\\
				& M3-50-4 & 1083 & - & - & 59895.6 & 10800 & - & - & 33742.7 & 10800 & \textbf{29277.2} & 8116 & -13.2\\
				& M3-50-6 & 1083 & - & - & 203165.1 & 10800 & - & - & 175755 & 10800 & \textbf{32293.4} & 4011 & -81.6\\
				& M3-100-1 & 1200 & - & - & 20373.9 & 10800 & - & - & \textbf{20203.8} & 10800 & 20517.3 & 5668 & 1.6\\
				& M3-100-2 & 1200 & - & - & 170434.9 & 10800 & - & - & 29682.8 & 10800 & \textbf{21577.1} & 1334 & -27.3\\
				& M3-100-4 & 1200 & - & - & 42747.9 & 10800 & - & - & 1622680.6 & 10800 & \textbf{24522.1} & 3909 & -98.5\\
				& M3-100-6 & 1200 & - & - & 89812.1 & 10800 & - & - & 116399.8 & 10800 & \textbf{27454.1} & 7891 & -76.4\\
				& M3-150-1 & 1261 & - & - & 17593 & 10800 & - & - & \textbf{17546.5} & 10800 & 17621.8 & 3598 & 0.4\\
				& M3-150-2 & 1261 & - & - & 23463 & 10800 & - & - & 19950.1 & 10800 & \textbf{18780.2} & 1099 & -5.9\\
				& M3-150-4 & 1261 & - & - & 35304.5 & 10800 & - & - & 34249.1 & 10800 & \textbf{21244.2} & 3481 & -38\\
				& M3-150-6 & 1261 & - & - & 88844 & 10800 & - & - & 110504.1 & 10800 & \textbf{24475.8} & 8187 & -77.9\\
				& M4-0-1 & 542 & \textbf{42423.5}\footnotemark[1] & 2512.2 & \textbf{42423.5}\footnotemark[1] & 3.7 & \textbf{42423.5} & 10800 & \textbf{42423.5}\footnotemark[1] & 3.9 & 42806.5 & 8603 & 0.9\\
				& M4-0-2 & 542 & \textbf{43507.5} & 10800 & \textbf{43507.5}\footnotemark[1] & 87.8 & - & - & \textbf{43507.5}\footnotemark[1] & 134.5 & 43544.5 & 2646 & 0.09\\
				& M4-0-4 & 542 & 51613.4 & 10800 & \textbf{46045.5} & 10800 & - & - & 46483.5 & 10800 & 46222.4 & 3628 & -0.6\\
				& M4-0-6 & 542 & 527916.5 & 10800 & 50021.6 & 10800 & - & - & 50467.1 & 10800 & \textbf{49407} & 4914 & -2.1\\
				& M4-50-1 & 1029 & 80541.7 & 10800 & \textbf{29986.1}\footnotemark[1] & 831.1 & - & - & \textbf{29986.1}\footnotemark[1] & 921.3 & 31662.7 & 1370 & 5.6\\
				& M4-50-2 & 1029 & 54139.7 & 10800 & \textbf{31162.7} & 10800 & - & - & 31175.2 & 10800 & 32342 & 5033 & 3.7\\
				& M4-50-4 & 1029 & 102941.4 & 10800 & 43659.6 & 10800 & - & - & 36993.4 & 10800 & \textbf{35156.6} & 4879 & -5.0\\
				& M4-50-6 & 1029 & 181828.1 & 10800 & 43642.2 & 10800 & - & - & 46535.8 & 10800 & \textbf{38576.7} & 5102 & -17.1\\
				\bottomrule
			\end{tabular}
			
			\begin{tablenotes}\footnotesize
				\item[1] The gap is calculated as the relative difference of the objective function value obtained with the heuristic method with respect to the corresponding value of the RN-wS formulation.
			\end{tablenotes}
		\end{threeparttable}
	\end{adjustbox}
	
\end{table}

\begin{table}[H]
	\small
	\centering
	\caption{Exhaustive results for the compact MILP formulations and the heuristic.}
	\label{tab:exhaustive_results_4}
	\begin{adjustbox}{max width=\textwidth}
		\begin{threeparttable}
			
			\begin{tabular}{ c | l | c  | c c | c c | c c | c c  | c c | c }
				\toprule
				& & & \multicolumn{2}{c|}{\textbf{CG-noS}} &  \multicolumn{2}{c|}{\textbf{RN-noS}} & \multicolumn{2}{c|}{\textbf{CG-wS}} &  \multicolumn{2}{c|}{\textbf{RN-wS}} & \multicolumn{2}{c|}{\textbf{Heuristic}}  & \\
				& \makecell{Instance} & \makecell{Candidate \\ locations} & $\text{\makecell{Obj.\\ fun.}}$ & $\text{\makecell{Comp. \\ time (s)}}$ & $\text{\makecell{Obj. \\ fun.}}$ & $\text{\makecell{Comp. \\ time (s)}}$ & $\text{\makecell{Obj. \\ fun.}}$ & $\text{\makecell{Comp. \\ time(s)}}$ & $\text{\makecell{Obj. \\ fun.}}$ & $\text{\makecell{Comp. \\ time (s)}}$ & $\text{\makecell{Obj. \\ fun.}}$ & $\text{\makecell{Comp. \\ time (s)}}$ & \makecell{Gap \% \\ to RN-wS}\tnote{1} \\
				\midrule
				\multirow{6}{*}{\STAB{\rotatebox[origin=c]{90}{Real-life instances}}}
				& M4-100-1 & 1197 & - & - & \textbf{27618.4}\footnotemark[1] & 10512.9 & - & - & \textbf{27618.4}\footnotemark[1] & 7921.5 & 27701.8 & 226 & 0.3\\
				& M4-100-2 & 1197 & 164966.9 & 10800 & 28957.8 & 10800 & - & - & 30124.8 & 10800 & \textbf{28624.6} & 9570 & -5.0\\
				& M4-100-4 & 1197 & 183809.4 & 10800 & - & - & - & - & 37173.9 & 10800 & \textbf{31234.9} & 7809 & -16.0\\
				& M4-100-6 & 1197 & - & - & 71002.9 & 10800 & - & - & 62974 & 10800 & \textbf{34473.4} & 457 & -45.3\\
				& M4-150-1 & 1319 & - & - & 25129.5 & 10800 & - & - & 25179.8 & 10800 & \textbf{24968.7} & 7674 & -0.8\\
				& M4-150-2 & 1319 & - & - & 30337.4 & 10800 & - & - & 29625 & 10800 & \textbf{25903.4} & 7360 & -12.6\\
				& M4-150-4 & 1319 & - & - & 39339.3 & 10800 & - & - & 33552.4 & 10800 & \textbf{28581.8} & 2597 & -14.8\\
				& M4-150-6 & 1319 & - & - & 52729.7 & 10800 & - & - & - & - & \textbf{31816.5} & 9606 & -\\
				\bottomrule
			\end{tabular}
			
			\begin{tablenotes}\footnotesize
				\item[1] The gap is calculated as the relative difference of the objective function value obtained with the heuristic method with respect to the corresponding value of the RN-wS formulation.
			\end{tablenotes}
		\end{threeparttable}
	\end{adjustbox}
	
\end{table}
\renewcommand*{\thefootnote}{\arabic{footnote}}

\newpage

\subsection{Practical aspects}
\label{sec:app_pract}

\begin{figure}[H]
	\begin{subfigure}{.24\textwidth}
		\centering
		\includegraphics[width=.9\textwidth]{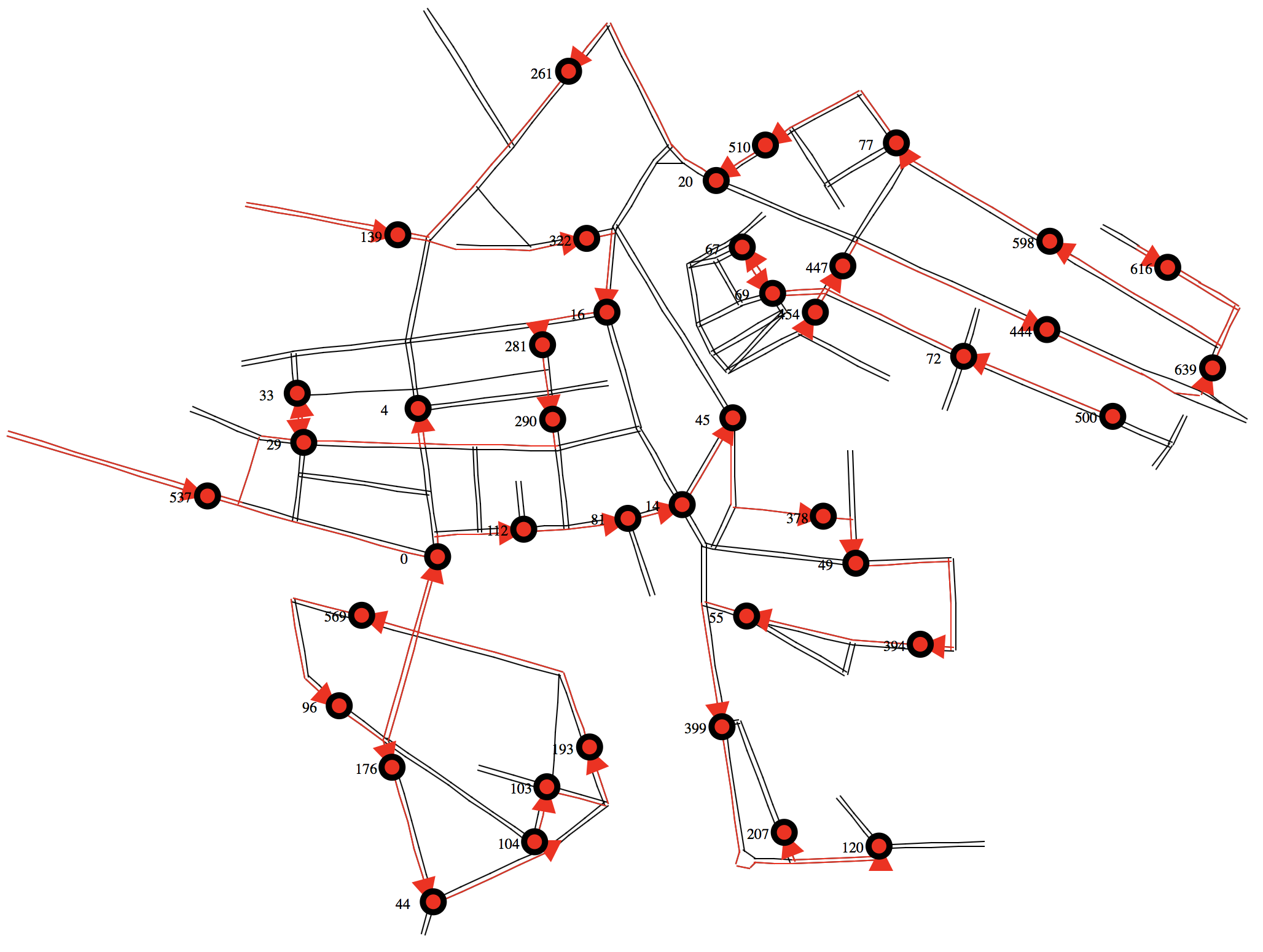}
		\caption{$\gamma=150, m=1$}
		\label{fig:sfig1}
	\end{subfigure}%
	\begin{subfigure}{.24\textwidth}
		\centering
		\includegraphics[width=.9\textwidth]{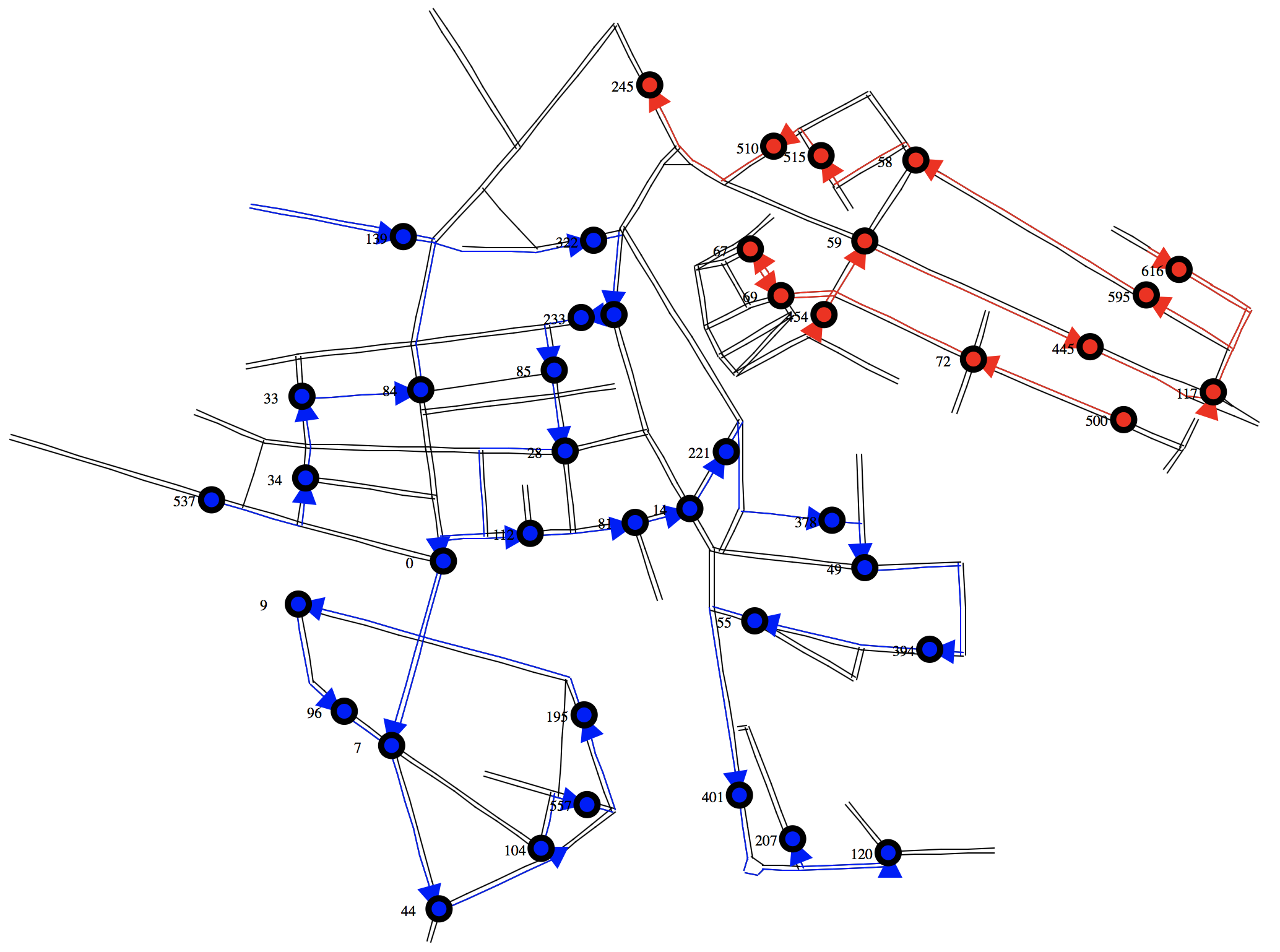}
		\caption{$\gamma=150, m=2$}
		\label{fig:sfig2}
	\end{subfigure}
	\begin{subfigure}{.24\textwidth}
		\centering
		\includegraphics[width=.9\textwidth]{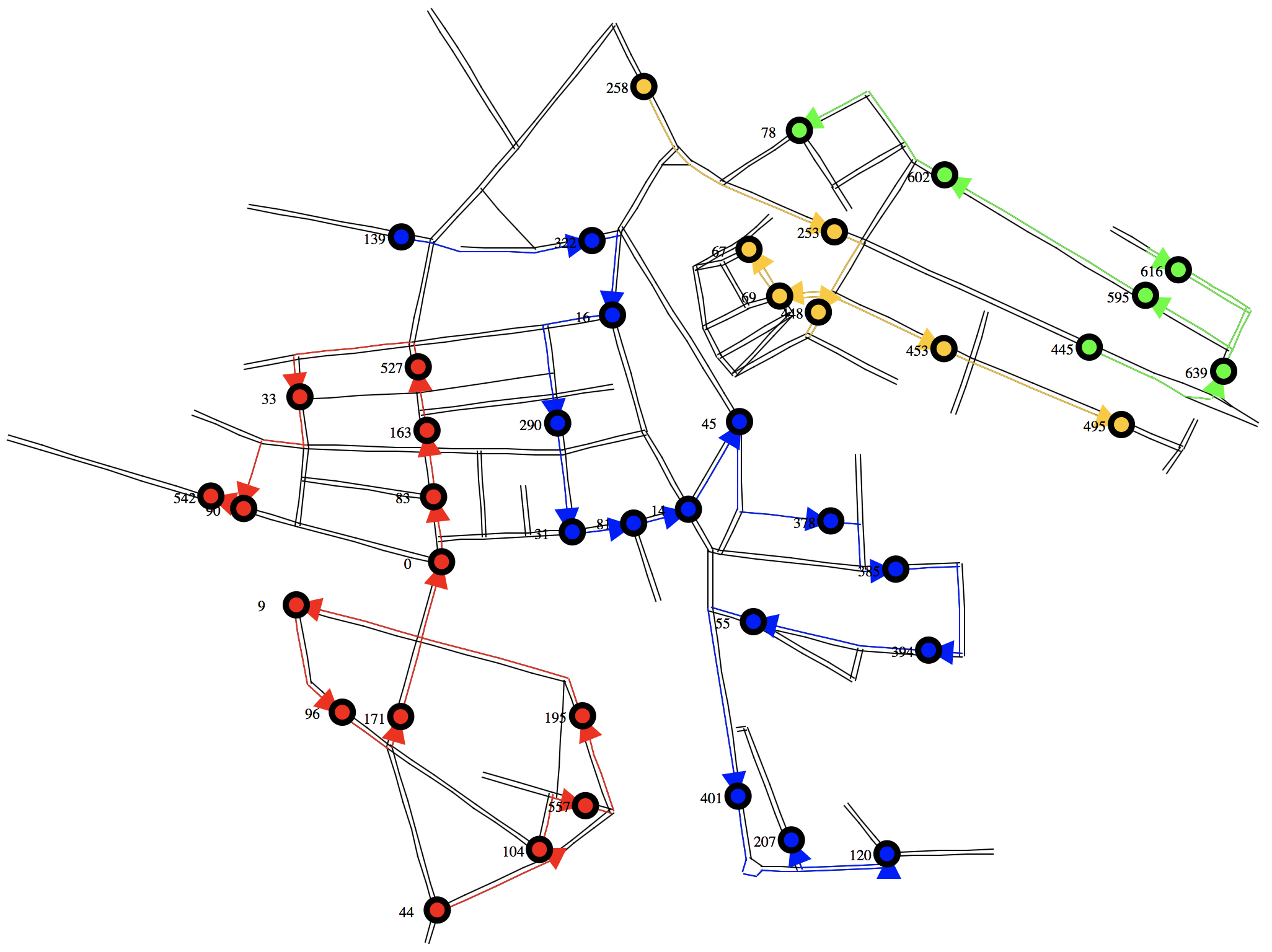}
		\caption{$\gamma=150, m=4$}
		\label{fig:sfig2}
	\end{subfigure}
	\begin{subfigure}{.24\textwidth}
		\centering
		\includegraphics[width=.9\textwidth]{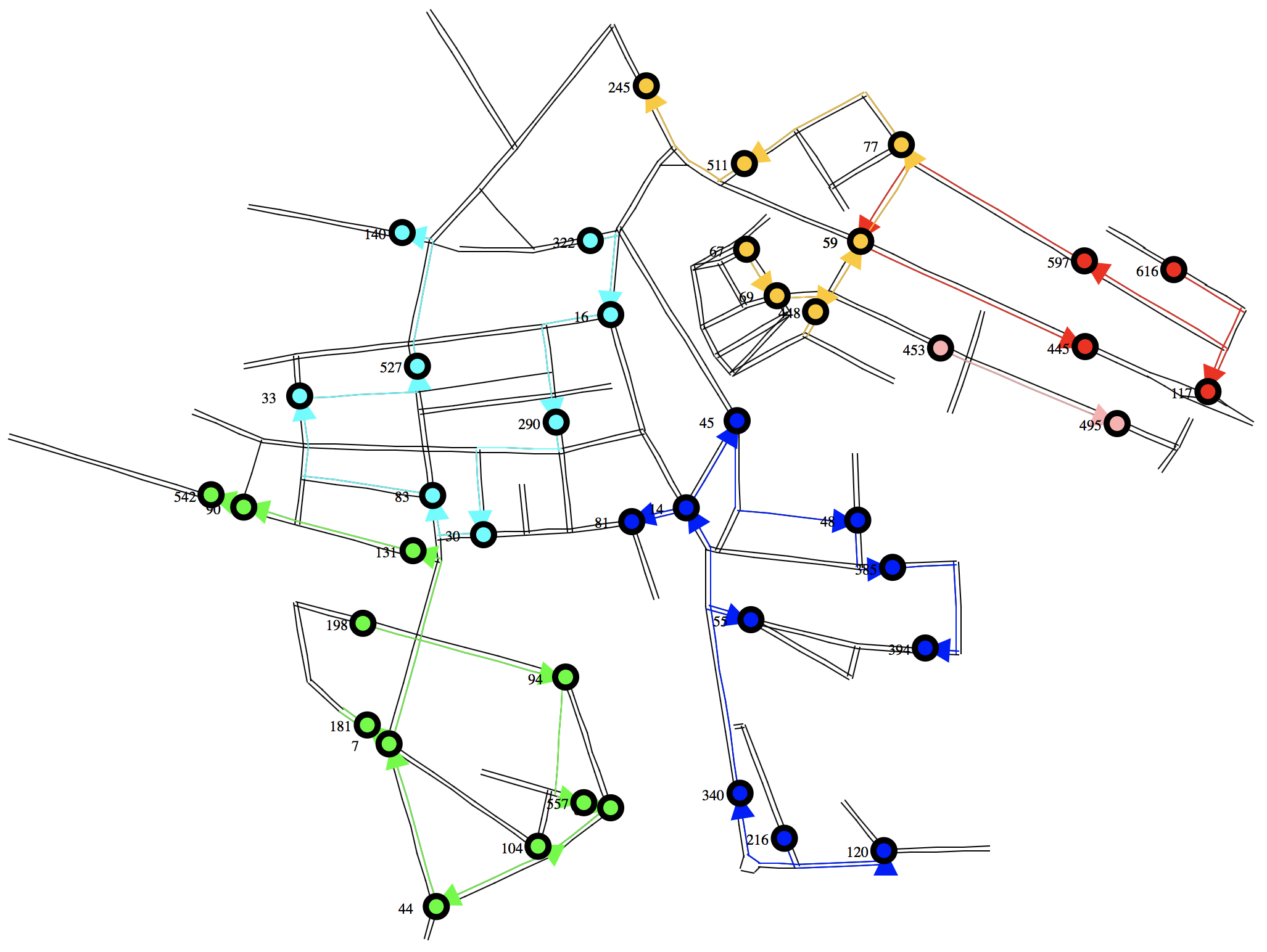}
		\caption{$\gamma=150, m=6$}
		\label{fig:sfig2}
	\end{subfigure}
	\begin{subfigure}{.24\textwidth}
		\centering
		\includegraphics[width=.9\textwidth]{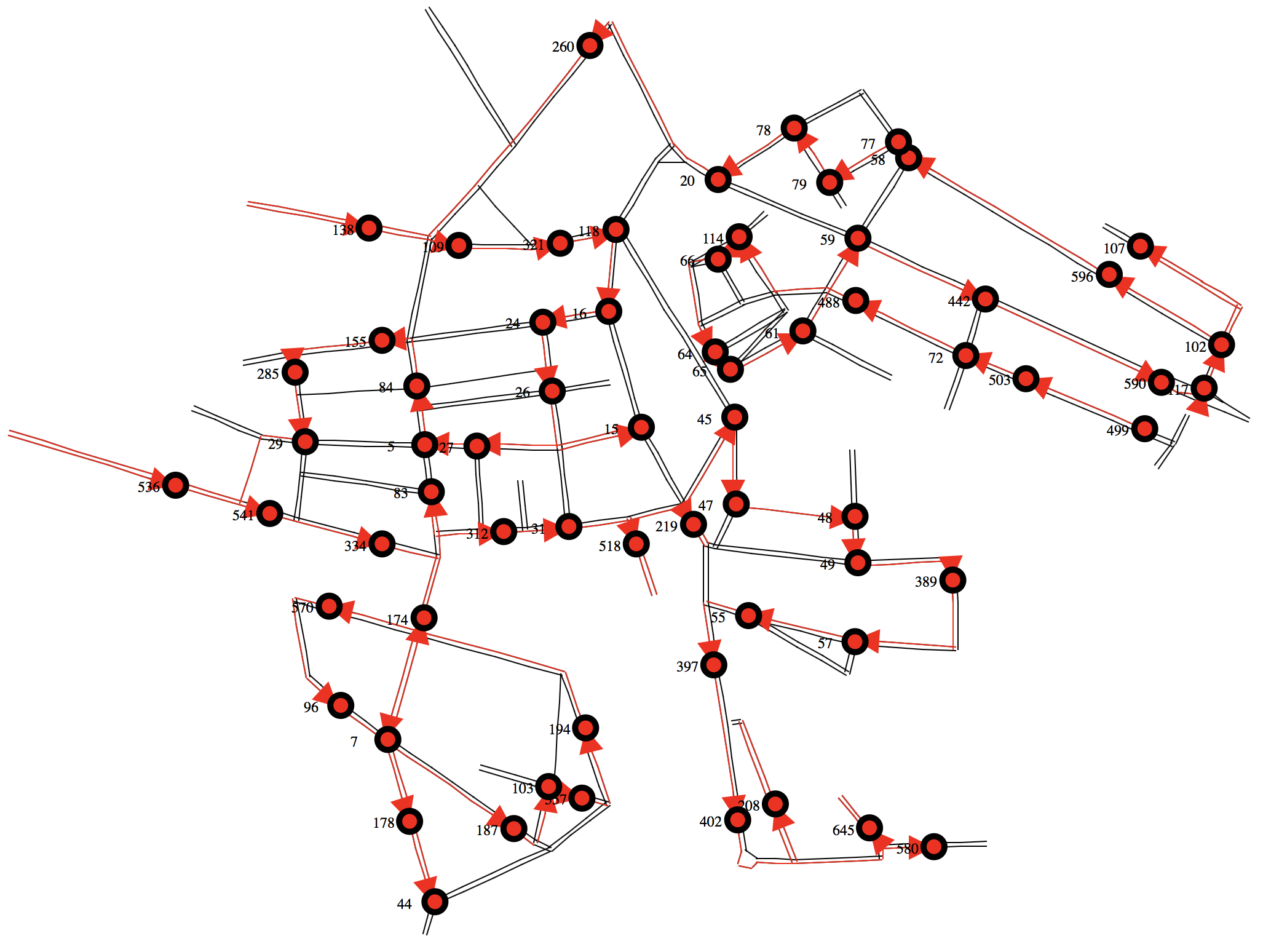}
		\caption{$\gamma=100, m=1$}
		\label{fig:sfig1}
	\end{subfigure}%
	\begin{subfigure}{.24\textwidth}
		\centering
		\includegraphics[width=.9\textwidth]{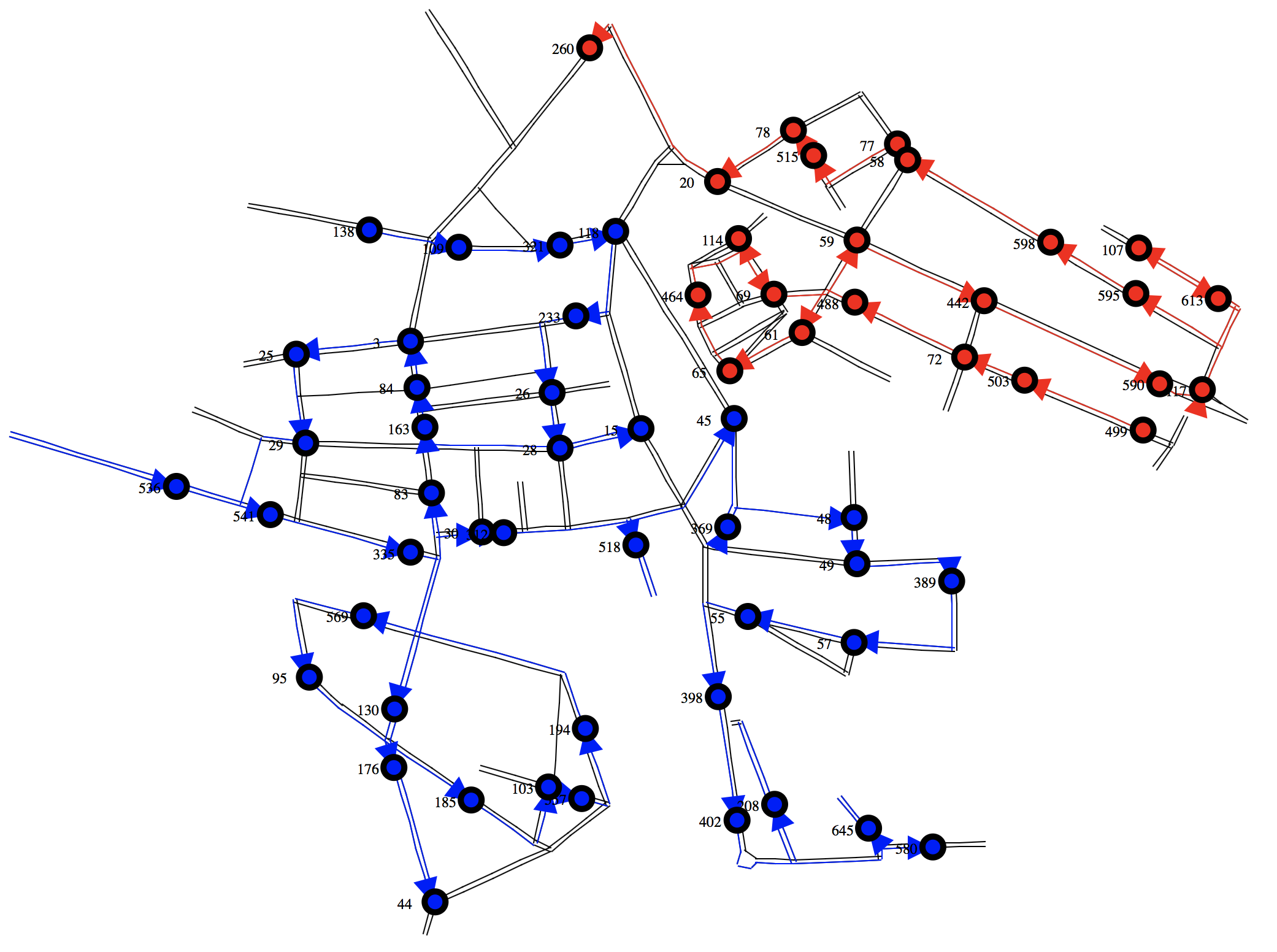}
		\caption{$\gamma=100, m=2$}
		\label{fig:sfig2}
	\end{subfigure}
	\begin{subfigure}{.24\textwidth}
		\centering
		\includegraphics[width=.9\textwidth]{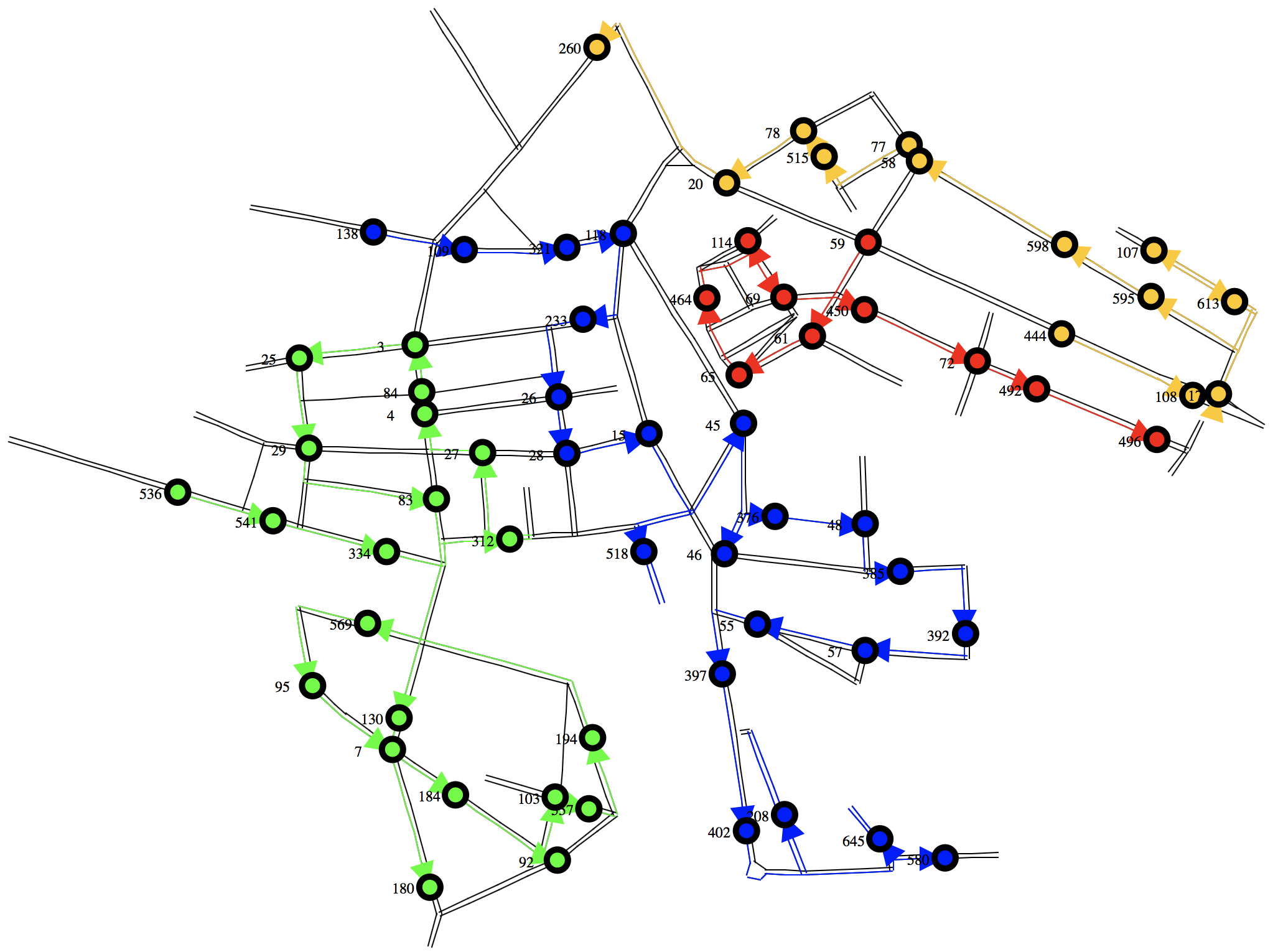}
		\caption{$\gamma=100, m=4$}
		\label{fig:sfig2}
	\end{subfigure}
	\begin{subfigure}{.24\textwidth}
		\centering
		\includegraphics[width=.9\textwidth]{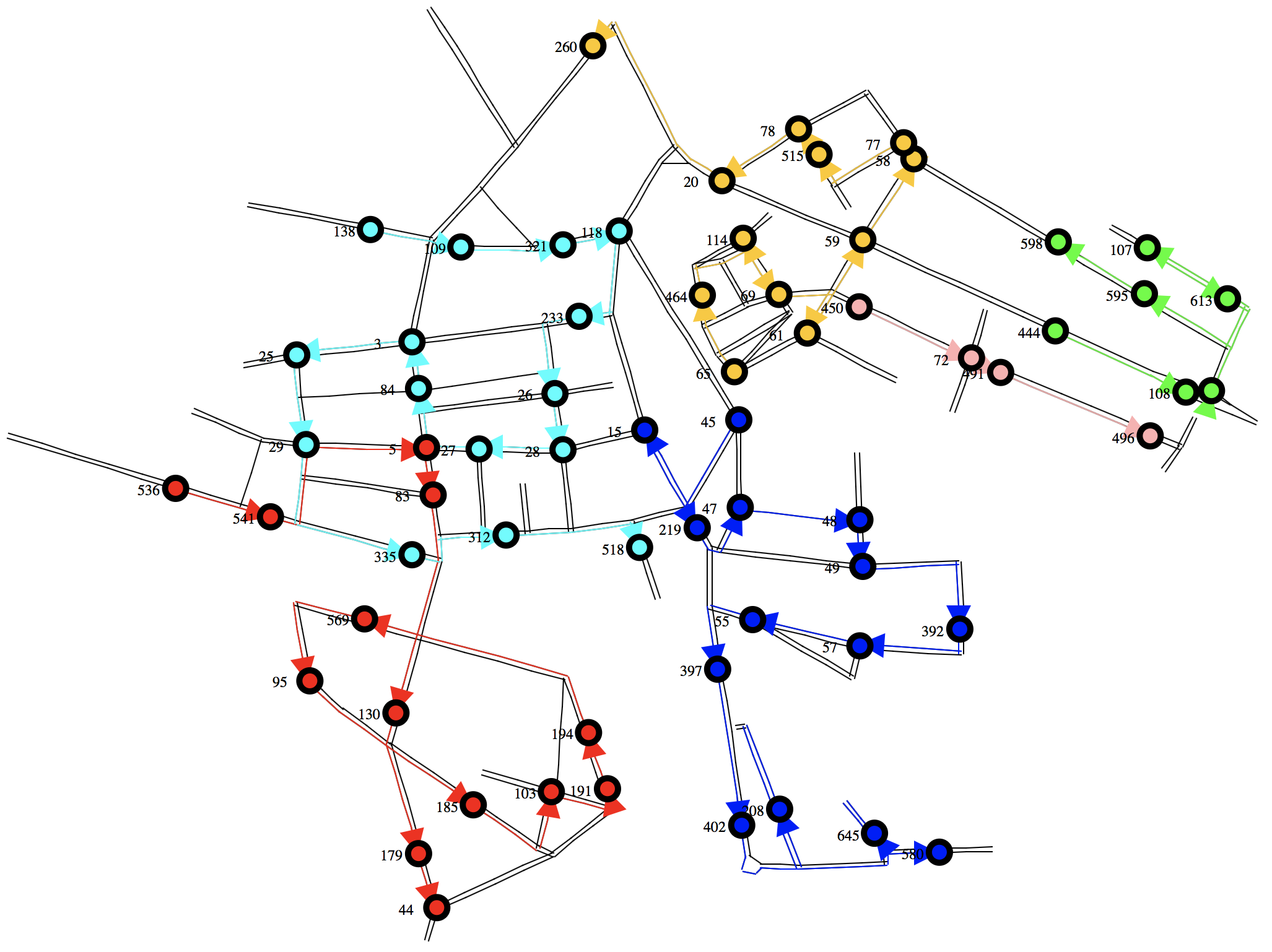}
		\caption{$\gamma=100, m=6$}
		\label{fig:sfig2}
	\end{subfigure}
	\begin{subfigure}{.24\textwidth}
		\centering
		\includegraphics[width=.9\textwidth]{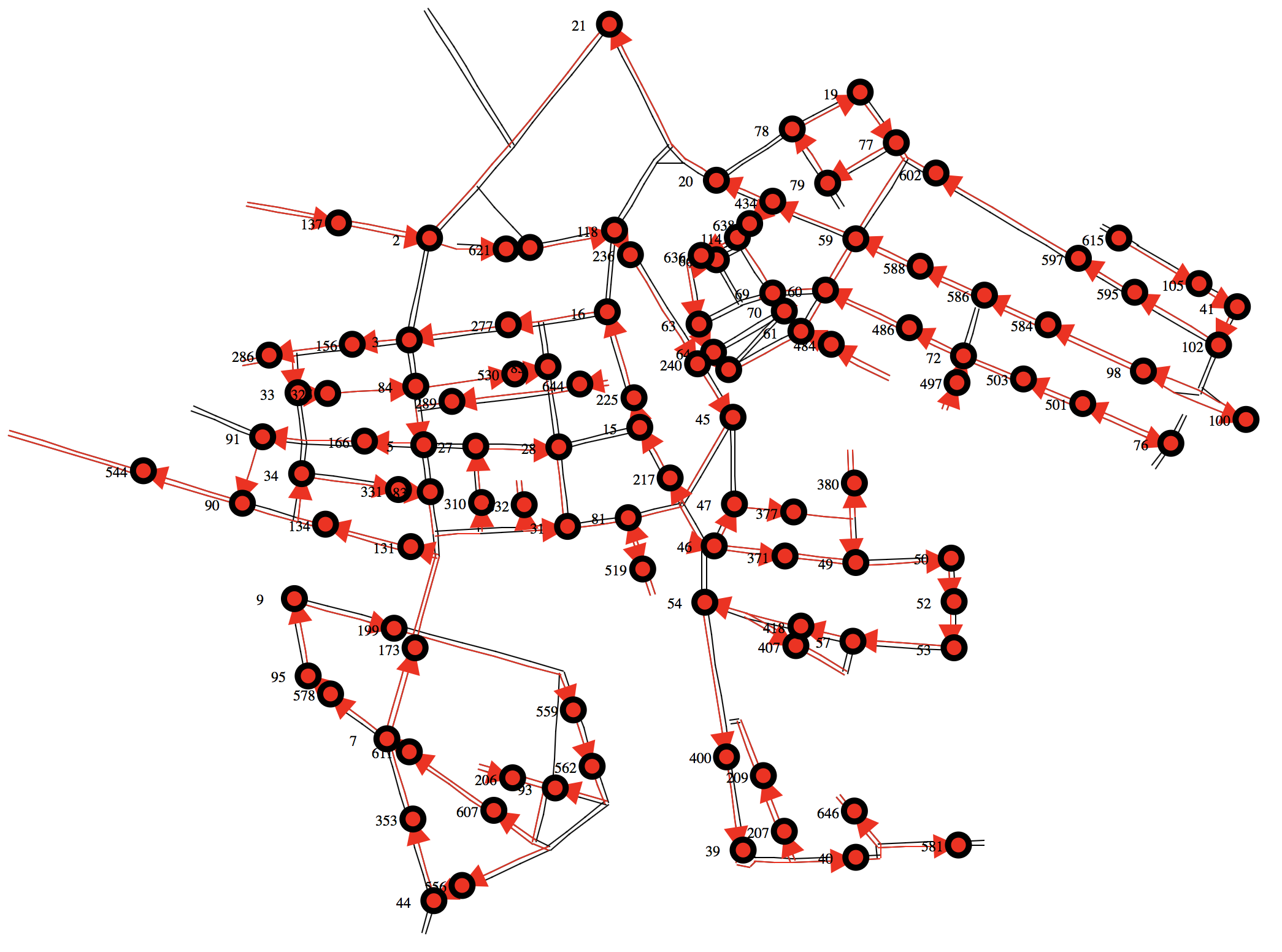}
		\caption{$\gamma=50, m=1$}
		\label{fig:sfig1}
	\end{subfigure}%
	\begin{subfigure}{.24\textwidth}
		\centering
		\includegraphics[width=.9\textwidth]{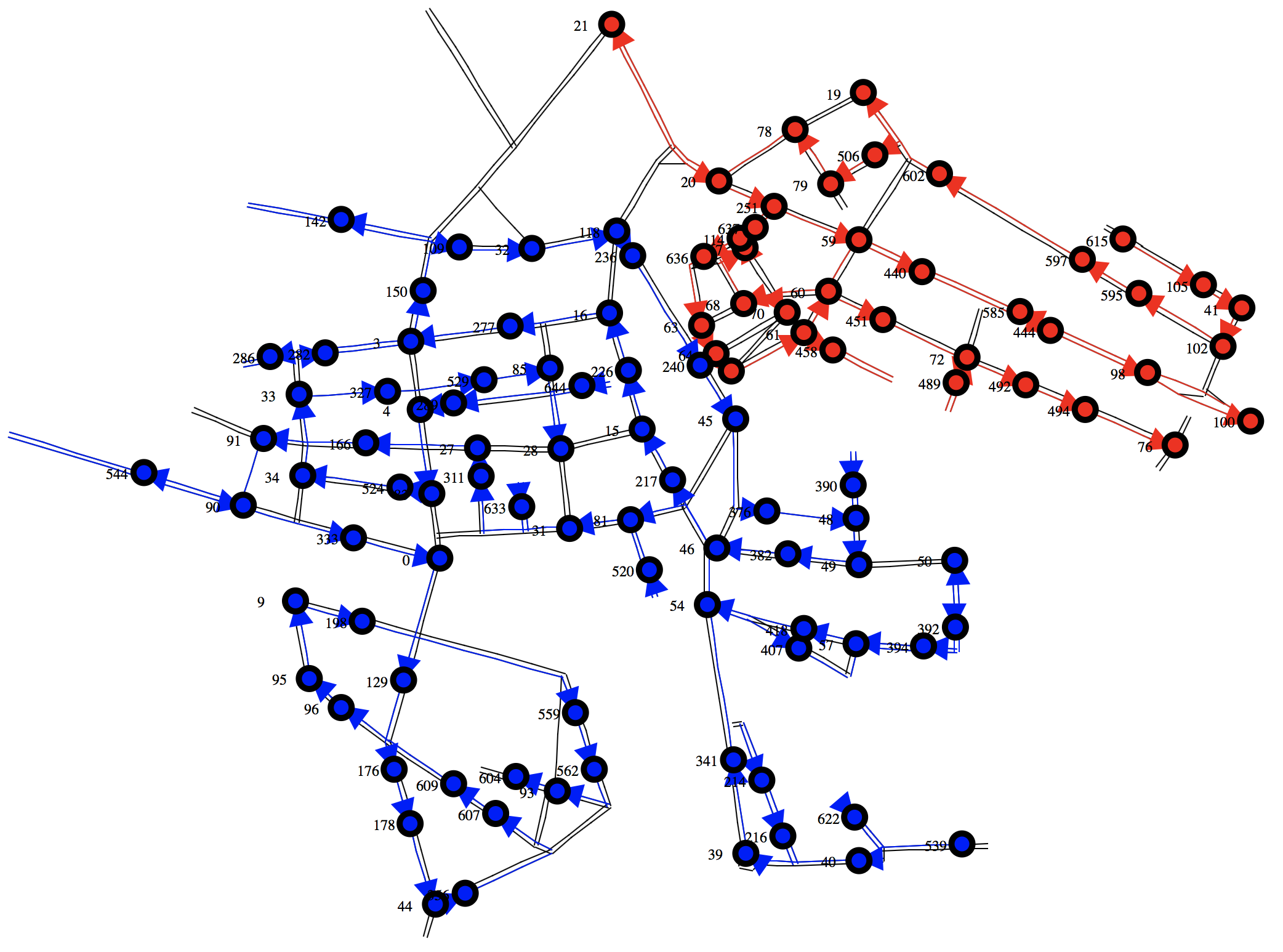}
		\caption{$\gamma=50, m=2$}
		\label{fig:sfig2}
	\end{subfigure}
	\begin{subfigure}{.24\textwidth}
		\centering
		\includegraphics[width=.9\textwidth]{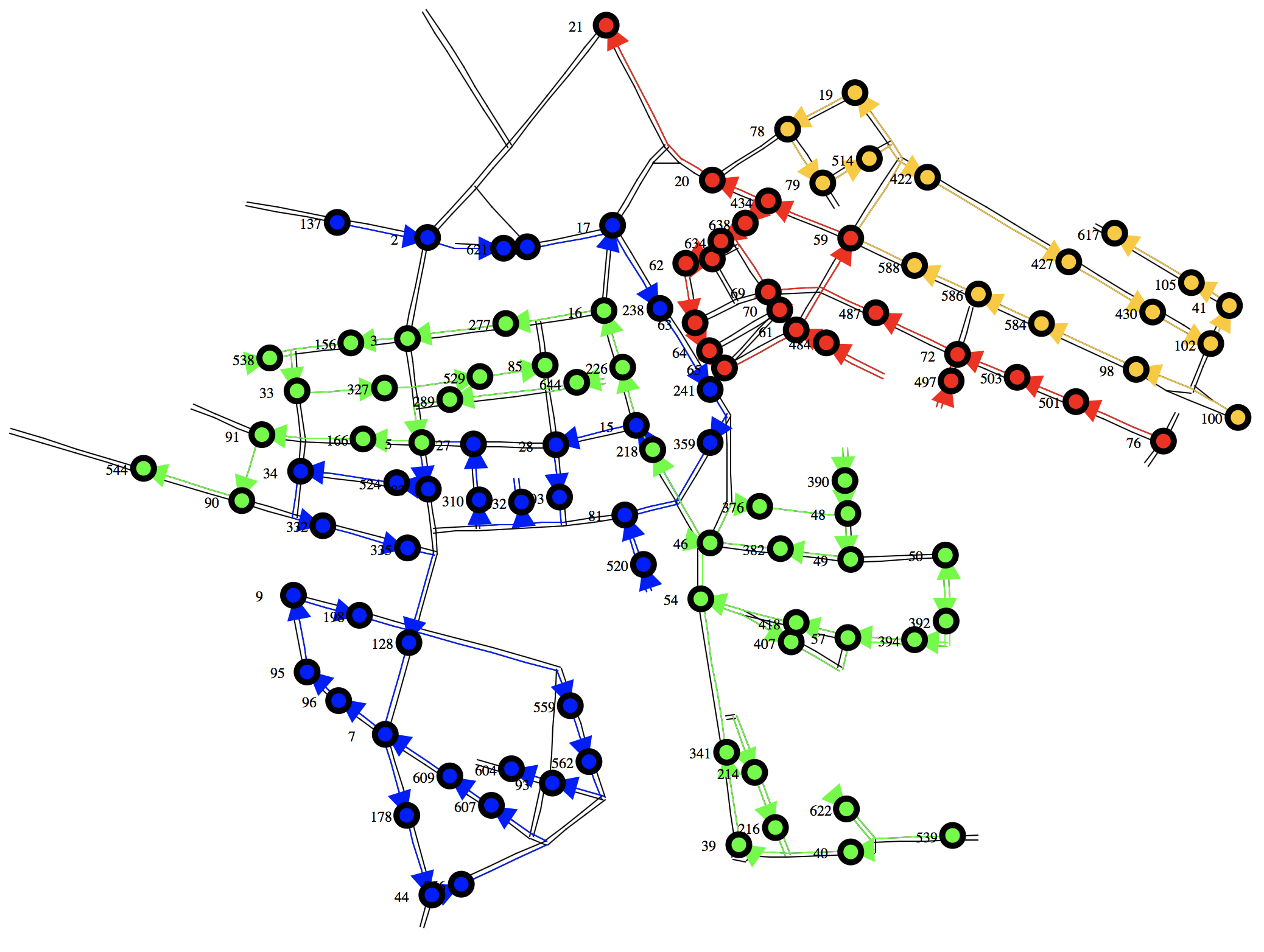}
		\caption{$\gamma=50, m=4$}
		\label{fig:sfig2}
	\end{subfigure}
	\begin{subfigure}{.24\textwidth}
		\centering
		\includegraphics[width=.9\textwidth]{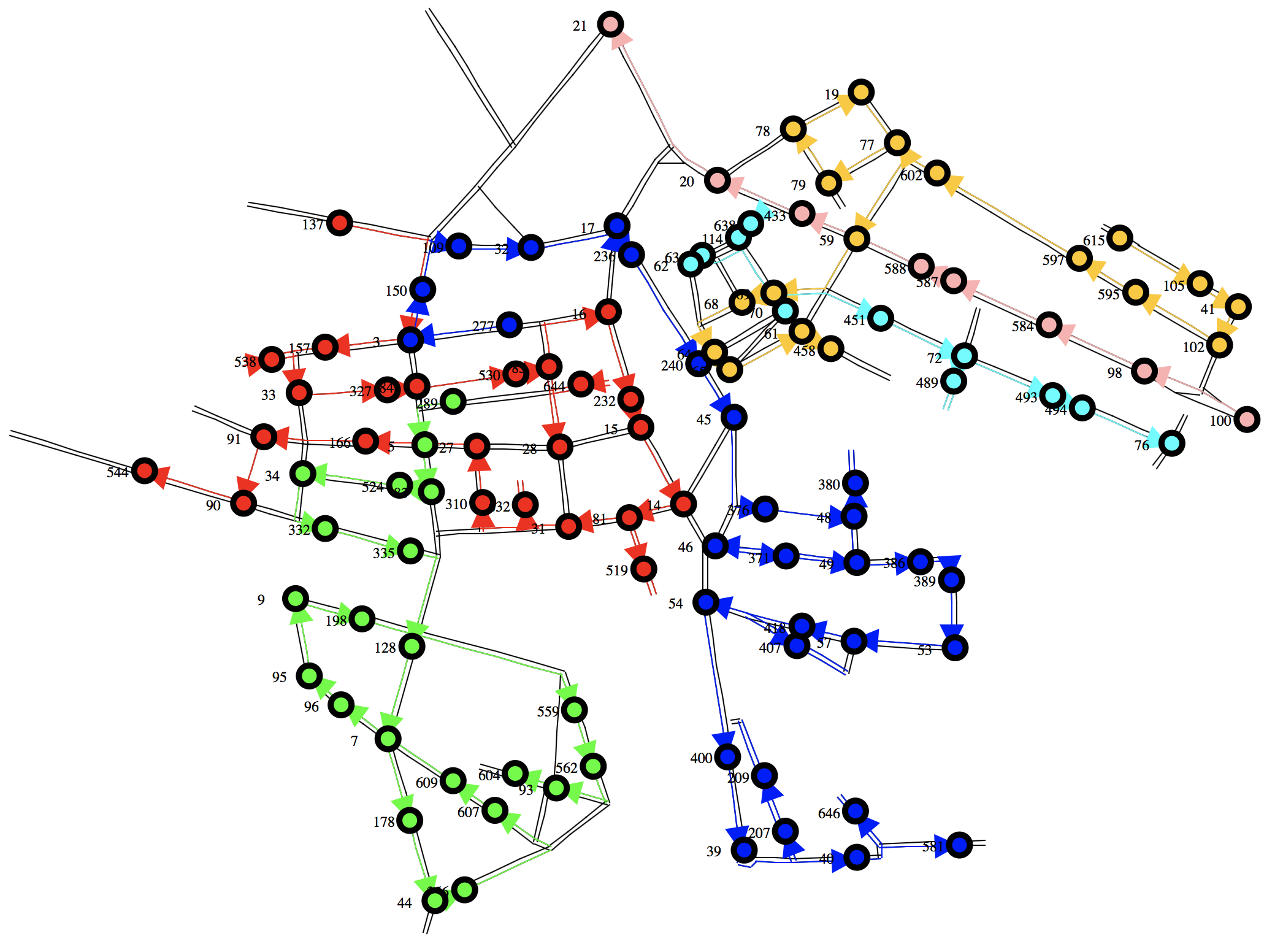}
		\caption{$\gamma=50, m=6$}
		\label{fig:sfig2}
	\end{subfigure}
	\begin{subfigure}{.24\textwidth}
		\centering
		\includegraphics[width=.9\textwidth]{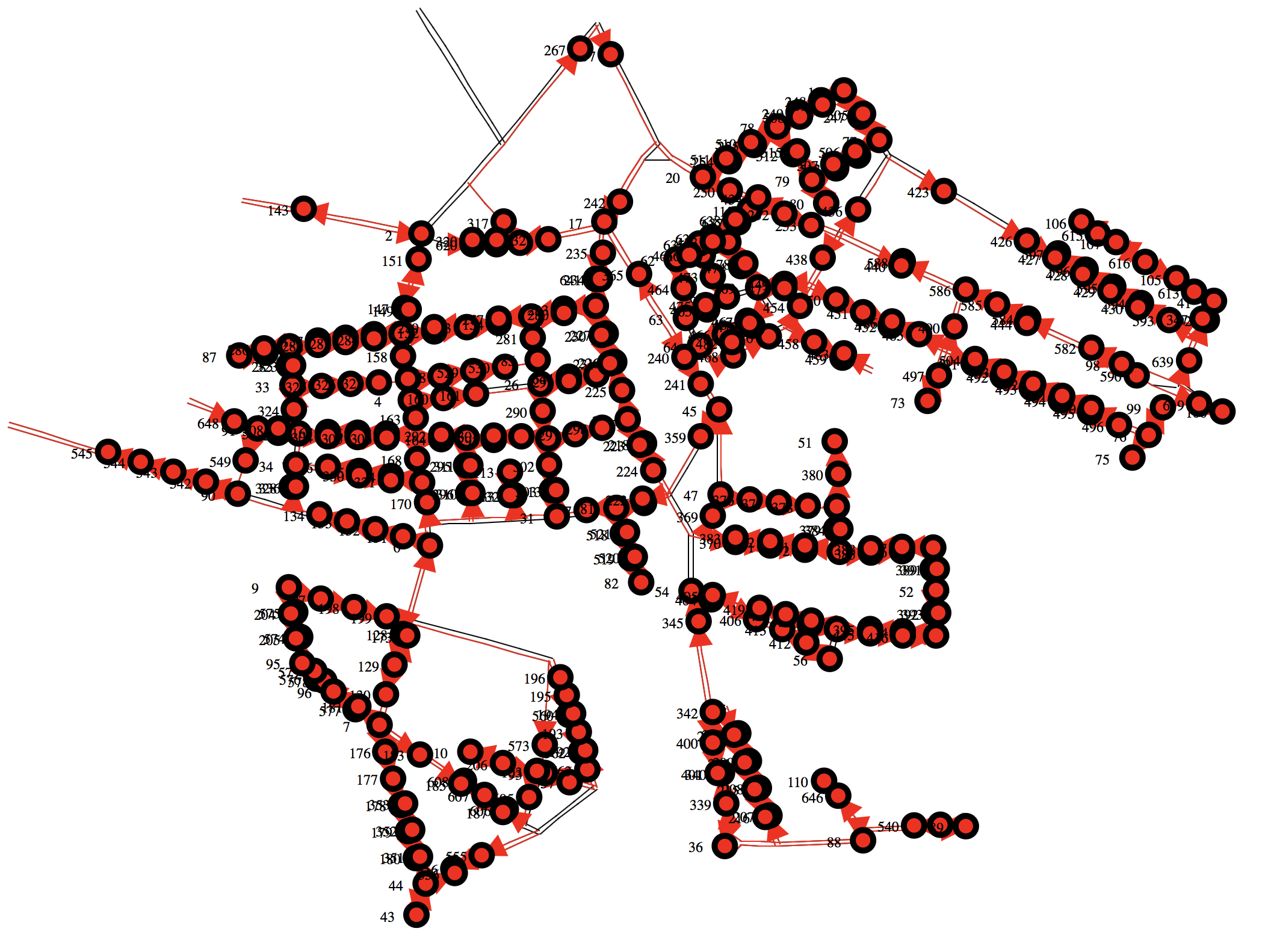}
		\caption{$\gamma=0, m=1$}
		\label{fig:sfig1}
	\end{subfigure}%
	\begin{subfigure}{.24\textwidth}
		\centering
		\includegraphics[width=.9\textwidth]{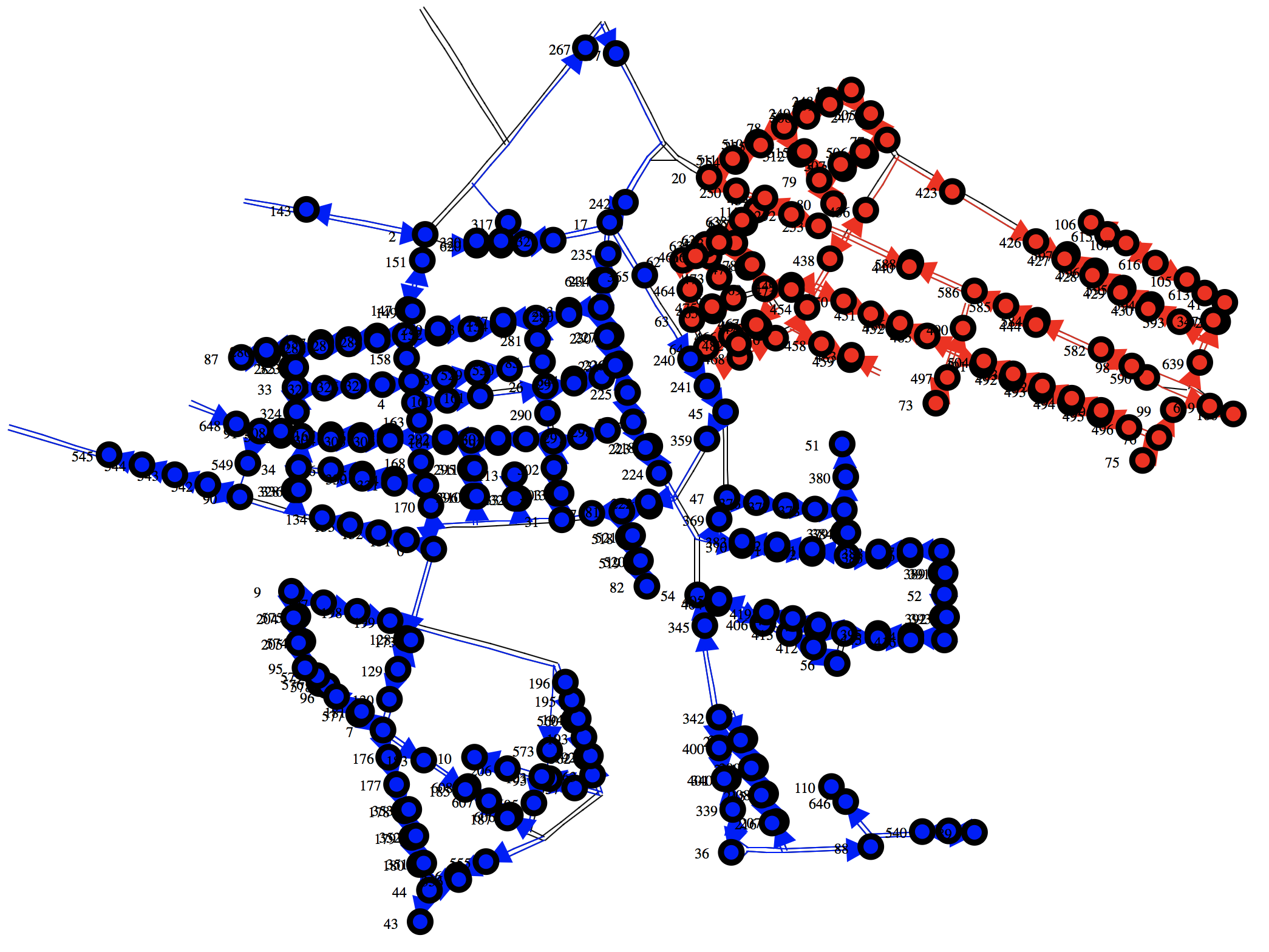}
		\caption{$\gamma=0, m=2$}
		\label{fig:sfig2}
	\end{subfigure}
	\begin{subfigure}{.24\textwidth}
		\centering
		\includegraphics[width=.9\textwidth]{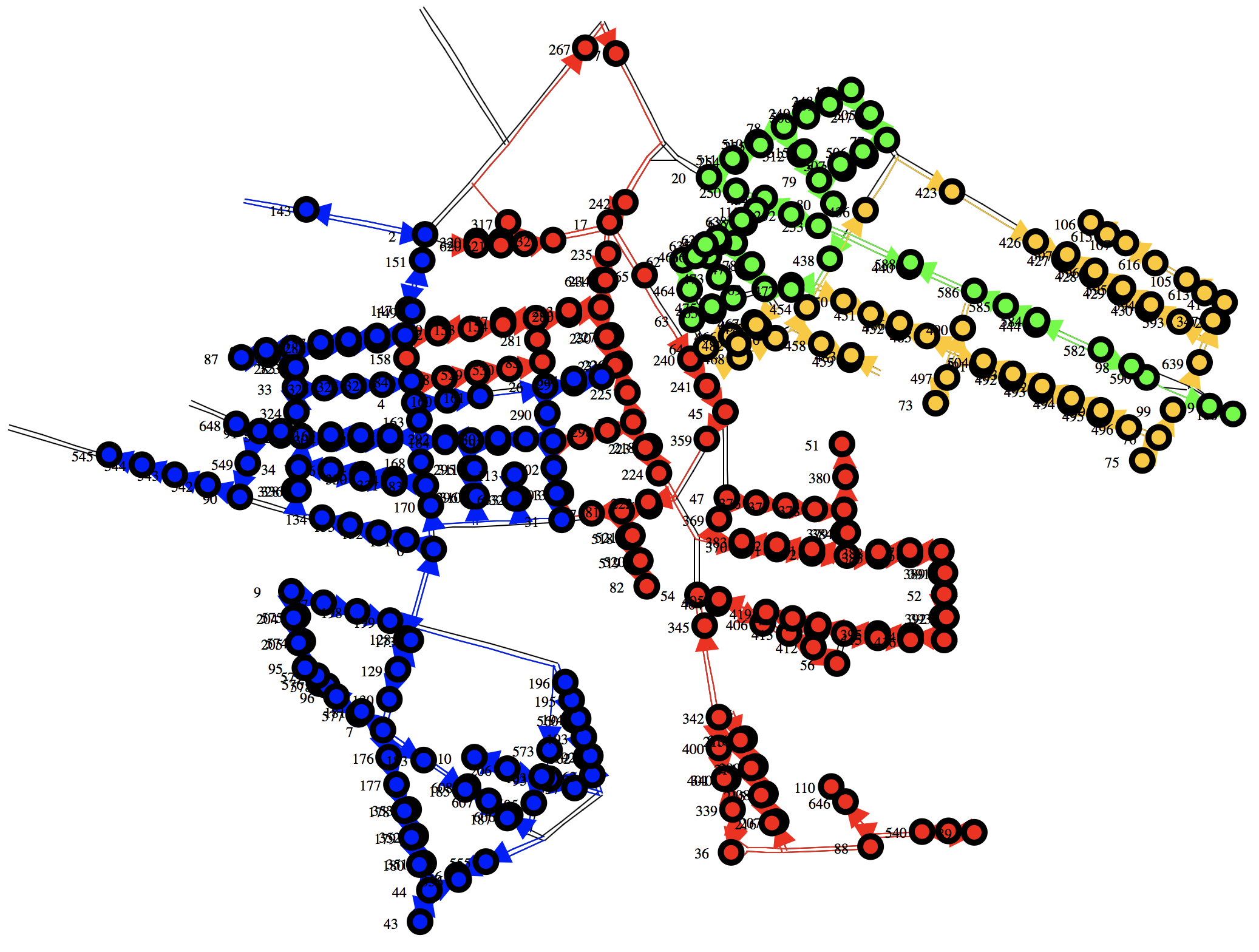}
		\caption{$\gamma=0, m=4$}
		\label{fig:sfig2}
	\end{subfigure}
	\begin{subfigure}{.24\textwidth}
		\centering
		\includegraphics[width=.9\textwidth]{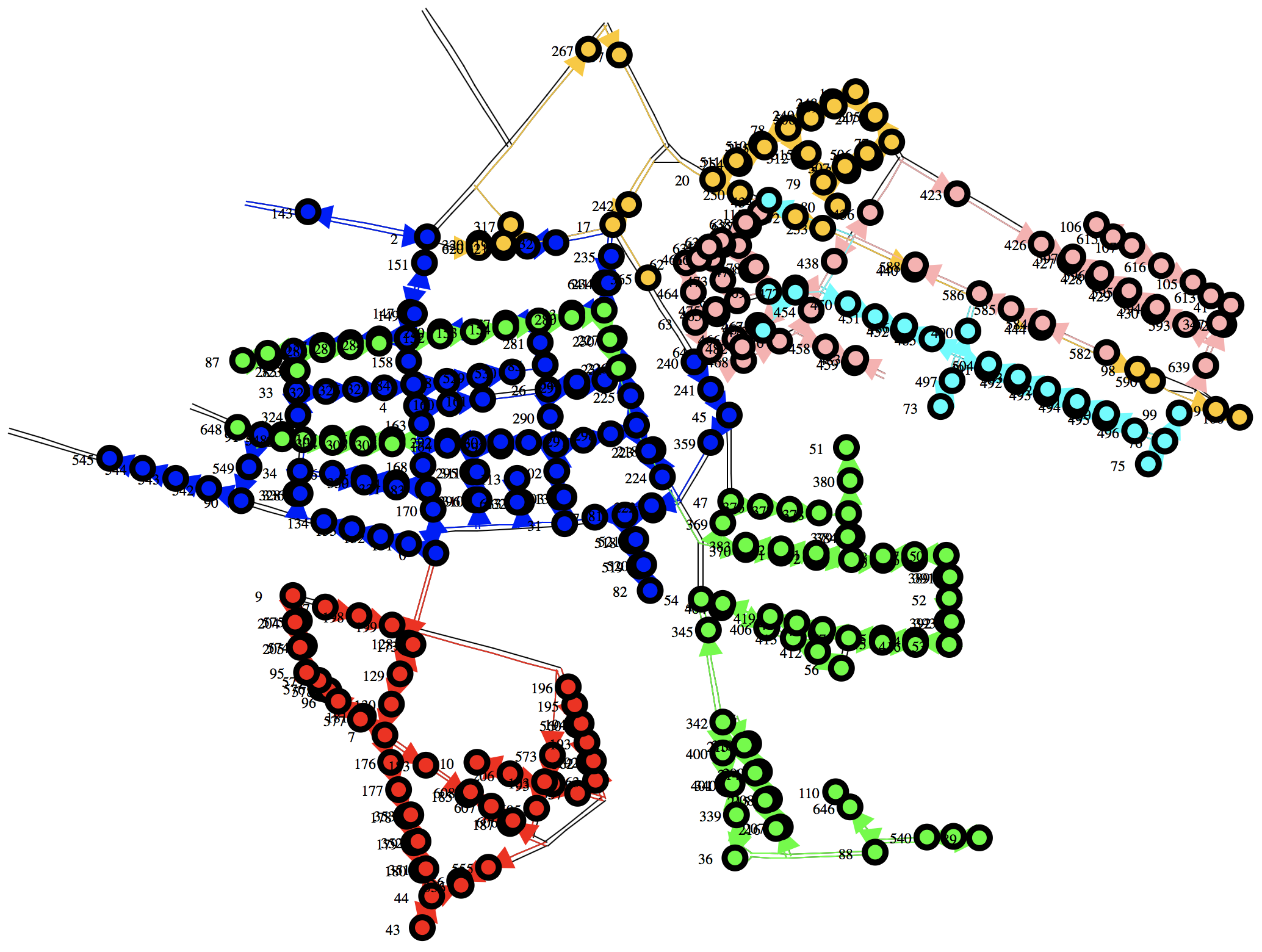}
		\caption{$\gamma=0, m=6$}
		\label{fig:sfig2}
	\end{subfigure}
	\caption{Solutions of the heuristic approach for instances of data set M1}
	\label{fig:fig_M1_solutions}
\end{figure}




\end{document}